`


\input amstex
\expandafter\ifx\csname mathdefs.tex\endcsname\relax
  \expandafter\gdef\csname mathdefs.tex\endcsname{}
\else \message{Hey!  Apparently you were trying to
  \string twice.   This does not make sense.} 
\errmessage{Please edit your file (probably \jobname.tex) and remove
any duplicate ``\string\input'' lines} \fi




\catcode`\X=12\catcode`\@=11

\def\n@wcount{\alloc@0\count\countdef\insc@unt}
\def\n@wwrite{\alloc@7\write\chardef\sixt@@n}
\def\n@wread{\alloc@6\read\chardef\sixt@@n}
\def\r@s@t{\relax}\def\v@idline{\par}\def\@mputate#1/{#1}
\def\l@c@l#1X{\firstpart.#1}\def\gl@b@l#1X{#1}\def\t@d@l#1X{{}}

\def\crossrefs#1{\ifx\all#1\let\tr@ce=\all\else\def\tr@ce{#1,}\fi
   \n@wwrite\cit@tionsout\openout\cit@tionsout=\jobname.cit 
   \write\cit@tionsout{\tr@ce}\expandafter\setfl@gs\tr@ce,}
\def\setfl@gs#1,{\def\@{#1}\ifx\@\empty\let\next=\relax
   \else\let\next=\setfl@gs\expandafter\xdef
   \csname#1tr@cetrue\endcsname{}\fi\next}
\def\m@ketag#1#2{\expandafter\n@wcount\csname#2tagno\endcsname
     \csname#2tagno\endcsname=0\let\tail=\all\xdef\all{\tail#2,}
   \ifx#1\l@c@l\let\tail=\r@s@t\xdef\r@s@t{\csname#2tagno\endcsname=0\tail}\fi
   \expandafter\gdef\csname#2cite\endcsname##1{\expandafter
     \ifx\csname#2tag##1\endcsname\relax?\else\csname#2tag##1\endcsname\fi
     \expandafter\ifx\csname#2tr@cetrue\endcsname\relax\else
     \write\cit@tionsout{#2tag ##1 cited on page \folio.}\fi}
   \expandafter\gdef\csname#2page\endcsname##1{\expandafter
     \ifx\csname#2page##1\endcsname\relax?\else\csname#2page##1\endcsname\fi
     \expandafter\ifx\csname#2tr@cetrue\endcsname\relax\else
     \write\cit@tionsout{#2tag ##1 cited on page \folio.}\fi}
   \expandafter\gdef\csname#2tag\endcsname##1{\expandafter
      \ifx\csname#2check##1\endcsname\relax
      \expandafter\xdef\csname#2check##1\endcsname{}%
      \else\immediate\write16{Warning: #2tag ##1 used more than once.}\fi
      \multit@g{#1}{#2}##1/X%
      \write\t@gsout{#2tag ##1 assigned number \csname#2tag##1\endcsname\space
      on page \number\count0.}%
   \csname#2tag##1\endcsname}}
\def\multit@g#1#2#3/#4X{\def\t@mp{#4}\ifx\t@mp\empty%
      \global\advance\csname#2tagno\endcsname by 1 
      \expandafter\xdef\csname#2tag#3\endcsname
      {#1\number\csname#2tagno\endcsnameX}%
   \else\expandafter\ifx\csname#2last#3\endcsname\relax
      \expandafter\n@wcount\csname#2last#3\endcsname
      \global\advance\csname#2tagno\endcsname by 1 
      \expandafter\xdef\csname#2tag#3\endcsname
      {#1\number\csname#2tagno\endcsnameX}
      \write\t@gsout{#2tag #3 assigned number \csname#2tag#3\endcsname\space
      on page \number\count0.}\fi
   \global\advance\csname#2last#3\endcsname by 1
   \def\t@mp{\expandafter\xdef\csname#2tag#3/}%
   \expandafter\t@mp\@mputate#4\endcsname
   {\csname#2tag#3\endcsname\lastpart{\csname#2last#3\endcsname}}\fi}
\def\t@gs#1{\def\all{}\m@ketag#1e\m@ketag#1s\m@ketag\t@d@l p
   \m@ketag\gl@b@l r \n@wread\t@gsin
   \openin\t@gsin=\jobname.tgs \re@der \closein\t@gsin
   \n@wwrite\t@gsout\openout\t@gsout=\jobname.tgs }
\outer\def\localtags{\t@gs\l@c@l}
\outer\def\globaltags{\t@gs\gl@b@l}
\outer\def\newlocaltag#1{\m@ketag\l@c@l{#1}}
\outer\def\newglobaltag#1{\m@ketag\gl@b@l{#1}}

\newif\ifpr@ 
\def\m@kecs #1tag #2 assigned number #3 on page #4.%
   {\expandafter\gdef\csname#1tag#2\endcsname{#3}
   \expandafter\gdef\csname#1page#2\endcsname{#4}
   \ifpr@\expandafter\xdef\csname#1check#2\endcsname{}\fi}
\def\re@der{\ifeof\t@gsin\let\next=\relax\else
   \read\t@gsin to\t@gline\ifx\t@gline\v@idline\else
   \expandafter\m@kecs \t@gline\fi\let \next=\re@der\fi\next}
\def\pretags#1{\pr@true\pret@gs#1,,}
\def\pret@gs#1,{\def\@{#1}\ifx\@\empty\let\n@xtfile=\relax
   \else\let\n@xtfile=\pret@gs \openin\t@gsin=#1.tgs \message{#1} \re@der 
   \closein\t@gsin\fi \n@xtfile}

\newcount\sectno\sectno=0\newcount\subsectno\subsectno=0
\newif\ifultr@local \def\ultralocal{\ultr@localtrue}
\def\firstpart{\number\sectno}
\def\lastpart#1{\ifcase#1 \or a\or b\or c\or d\or e\or f\or g\or h\or 
   i\or k\or l\or m\or n\or o\or p\or q\or r\or s\or t\or u\or v\or w\or 
   x\or y\or z \fi}

\def\resetall{\global\advance\sectno by 1\subsectno=0
   \gdef\firstpart{\number\sectno}\r@s@t}
\def\resetsub{\global\advance\subsectno by 1
   \gdef\firstpart{\number\sectno.\number\subsectno}\r@s@t}
\def\newsection#1\par{\resetall\vskip0pt plus.3\vsize\penalty-250
   \vskip0pt plus-.3\vsize\bigskip\bigskip
   \message{#1}\leftline{\bf#1}\nobreak\bigskip}
\def\subsection#1\par{\ifultr@local\resetsub\fi
   \vskip0pt plus.2\vsize\penalty-250\vskip0pt plus-.2\vsize
   \bigskip\smallskip\message{#1}\leftline{\bf#1}\nobreak\medskip}

\def\t@gsoff#1,{\def\@{#1}\ifx\@\empty\let\next=\relax\else\let\next=\t@gsoff
   \def\@@{p}\ifx\@\@@\else
   \expandafter\gdef\csname#1cite\endcsname##1{\zeigen{##1}}
   \expandafter\gdef\csname#1page\endcsname##1{?}
   \expandafter\gdef\csname#1tag\endcsname##1{\zeigen{##1}}\fi\fi\next}
\def\verbatimtags{\ifx\all\relax\else\expandafter\t@gsoff\all,\fi}
\def\zeigen#1{\hbox{$\langle$}#1\hbox{$\rangle$}}

\def\(#1){\edef\dot@g{\ifmmode\ifinner(\hbox{\noexpand\etag{#1}})
   \else\noexpand\eqno(\hbox{\noexpand\etag{#1}})\fi
   \else(\noexpand\ecite{#1})\fi}\dot@g}

\newif\ifbr@ck
\def\eat#1{}
\def\[#1]{\br@cktrue[\br@cket#1'X]}
\def\br@cket#1'#2X{\def\temp{#2}\ifx\temp\empty\let\next\eat
   \else\let\next\br@cket\fi
   \ifbr@ck\br@ckfalse\br@ck@t#1,X\else\br@cktrue#1\fi\next#2X}
\def\br@ck@t#1,#2X{\def\temp{#2}\ifx\temp\empty\let\neext\eat
   \else\let\neext\br@ck@t\def\temp{,}\fi
   \def\teemp{#1}\ifx\teemp\empty\else\rcite{#1}\fi\temp\neext#2X}
\def\resetbr@cket{\gdef\[##1]{[\rtag{##1}]}}
\def\references{\resetbr@cket\newsection References\par}

\newtoks\symb@ls\newtoks\s@mb@ls\newtoks\p@gelist\n@wcount\ftn@mber
    \ftn@mber=1\newif\ifftn@mbers\ftn@mbersfalse\newif\ifbyp@ge\byp@gefalse
\def\defm@rk{\ifftn@mbers\n@mberm@rk\else\symb@lm@rk\fi}
\def\n@mberm@rk{\xdef\m@rk{{\the\ftn@mber}}%
    \global\advance\ftn@mber by 1 }
\def\rot@te#1{\let\temp=#1\global#1=\expandafter\r@t@te\the\temp,X}
\def\r@t@te#1,#2X{{#2#1}\xdef\m@rk{{#1}}}
\def\b@@st#1{{$^{#1}$}}\def\str@p#1{#1}
\def\symb@lm@rk{\ifbyp@ge\rot@te\p@gelist\ifnum\expandafter\str@p\m@rk=1 
    \s@mb@ls=\symb@ls\fi\write\f@nsout{\number\count0}\fi \rot@te\s@mb@ls}
\def\byp@ge{\byp@getrue\n@wwrite\f@nsin\openin\f@nsin=\jobname.fns 
    \n@wcount\currentp@ge\currentp@ge=0\p@gelist={0}
    \re@dfns\closein\f@nsin\rot@te\p@gelist
    \n@wread\f@nsout\openout\f@nsout=\jobname.fns }
\def\m@kelist#1X#2{{#1,#2}}
\def\re@dfns{\ifeof\f@nsin\let\next=\relax\else\read\f@nsin to \f@nline
    \ifx\f@nline\v@idline\else\let\t@mplist=\p@gelist
    \ifnum\currentp@ge=\f@nline
    \global\p@gelist=\expandafter\m@kelist\the\t@mplistX0
    \else\currentp@ge=\f@nline
    \global\p@gelist=\expandafter\m@kelist\the\t@mplistX1\fi\fi
    \let\next=\re@dfns\fi\next}
\def\symbols#1{\symb@ls={#1}\s@mb@ls=\symb@ls} 
\def\bigsymbol{\textstyle}
\symbols{\bigsymbol\ast,\dagger,\ddagger,\sharp,\flat,\natural,\star}
\def\ftnumbers{\ftn@mberstrue} \def\ftsymbols{\ftn@mbersfalse}
\def\paginal{\byp@ge} \def\resetftnumbers{\ftn@mber=1}
\def\ftnote#1{\defm@rk\expandafter\expandafter\expandafter\footnote
    \expandafter\b@@st\m@rk{#1}}

\long\def\jump#1\endjump{}
\def\ssum{\mathop{\lower .1em\hbox{$\textstyle\Sigma$}}\nolimits}

\def\qed{\nobreak\kern 1em \vrule height .5em width .5em depth 0em}
\def\newneq{\hbox{\rlap{\hbox to 1\wd9{\hss$=$\hss}}\raise .1em 
   \hbox to 1\wd9{\hss$\scriptscriptstyle/$\hss}}}
\def\subsetne{\setbox9 = \hbox{$\subset$}\mathrel{\hbox{\rlap
   {\lower .4em \newneq}\raise .13em \hbox{$\subset$}}}}
\def\supsetne{\setbox9 = \hbox{$\subset$}\mathrel{\hbox{\rlap
   {\lower .4em \newneq}\raise .13em \hbox{$\supset$}}}}

\def\vbar{\mathchoice{\vrule height6.3ptdepth-.5ptwidth.8pt\kern-.8pt}
   {\vrule height6.3ptdepth-.5ptwidth.8pt\kern-.8pt}
   {\vrule height4.1ptdepth-.35ptwidth.6pt\kern-.6pt}
   {\vrule height3.1ptdepth-.25ptwidth.5pt\kern-.5pt}}
\def\f@dge{\mathchoice{}{}{\mkern.5mu}{\mkern.8mu}}
\def\b@c#1#2{{\rm \mkern#2mu\vbar\mkern-#2mu#1}}
\def\b@b#1{{\rm I\mkern-3.5mu #1}}
\def\b@a#1#2{{\rm #1\mkern-#2mu\f@dge #1}}
\def\bb#1{{\count4=`#1 \advance\count4by-64 \ifcase\count4\or\b@a A{11.5}\or
   \b@b B\or\b@c C{5}\or\b@b D\or\b@b E\or\b@b F \or\b@c G{5}\or\b@b H\or
   \b@b I\or\b@c J{3}\or\b@b K\or\b@b L \or\b@b M\or\b@b N\or\b@c O{5} \or
   \b@b P\or\b@c Q{5}\or\b@b R\or\b@a S{8}\or\b@a T{10.5}\or\b@c U{5}\or
   \b@a V{12}\or\b@a W{16.5}\or\b@a X{11}\or\b@a Y{11.7}\or\b@a Z{7.5}\fi}}

\catcode`\X=11 \catcode`\@=12

\expandafter\ifx\csname citeadd.tex\endcsname\relax
\expandafter\gdef\csname citeadd.tex\endcsname{}
\else \message{Hey!  Apparently you were trying to
\string twice.   This does not make sense.} 
\errmessage{Please edit your file (probably \jobname.tex) and remove
any duplicate ``\string\input'' lines} \fi

\sectno=-2   
\localtags
\NoBlackBoxes
\ifx\shlhetal\undefinedcontrolsequence\let\shlhetal\relax\fi
\define\mr{\medskip\roster}
\define\sn{\smallskip\noindent}
\define\mn{\medskip\noindent}
\define\bn{\bigskip\noindent}
\define\ub{\underbar}
\define\wilog{\text{without loss of generality}}
\define\ermn{\endroster\medskip\noindent}

\define\dbcu{\dsize\bigcup}
\define\nl{\newline}
\documentstyle {amsppt}
\topmatter
\title {\it Categoricity for Abstract Classes with Amalgamation} \endtitle
\rightheadtext{Categoricity}
\author {Saharon Shelah \thanks {\null\newline
Partially supported by the United States-Israel Binational Science 
Foundation and I thank Alice Leonhardt for the beautiful typing.\newline
Publ. No. 394 \newline
Done 6,9/88 \newline
Latest Revision - 98/Sept/23} \endthanks} \endauthor
\affil {Institute of Mathematics \\
The Hebrew University \\
Jerusalem, Israel
\medskip
Department of Mathematics \\
Rutgers University \\
New Brunswick, NJ  USA} \endaffil

\abstract {Let ${\frak K}$ be an abstract elementary class with
amalgamation, and Lowenheim Skolem number LS$({\frak K})$.  We prove that
for a suitable Hanf number $\chi_0$ if $\chi_0 < \lambda_0 \le \lambda_1$,
and ${\frak K}$ is categorical in $\lambda^+_1$ then it is categorical in
$\lambda_0$.} \endabstract
\endtopmatter
\document  

\expandafter\ifx\csname alice2jlem.tex\endcsname\relax
  \expandafter\gdef\csname alice2jlem.tex\endcsname{}
\else \message{Hey!  Apparently you were trying to
\string  twice.   This does not make sense.}
\errmessage{Please edit your file (probably \jobname.tex) and remove
any duplicate ``\string\input'' lines} \fi

\expandafter\ifx\csname bib4plain.tex\endcsname\relax
  \expandafter\gdef\csname bib4plain.tex\endcsname{}
\else \message{Hey!  Apparently you were trying to \string twice.   This does not make sense.}
\errmessage{Please edit your file (probably \jobname.tex) and remove
any duplicate ``\string\input'' lines} \fi

\def\renewcommand{\newcommand}	       
\edef\cite{\the\catcode`@}%
\catcode`@ = 11
\let\@oldatcatcode = \cite
\chardef\@letter = 11
\chardef\@other = 12
%
%
%
%
\def\@innerdef#1#2{\edef#1{\expandafter\noexpand\csname #2\endcsname}}%
%
%
\@innerdef\@innernewcount{newcount}%
\@innerdef\@innernewdimen{newdimen}%
\@innerdef\@innernewif{newif}%
\@innerdef\@innernewwrite{newwrite}%
%
%
%
\def\@gobble#1{}%
%
%
%
\ifx\inputlineno\@undefined
   \let\@linenumber = \empty 
\else
   \def\@linenumber{\the\inputlineno:\space}%
\fi
%
%
%
\def\@futurenonspacelet#1{\def\cs{#1}%
   \afterassignment\@stepone\let\@nexttoken=
}%
\begingroup 
\def\\{\global\let\@stoken= }%
\\ 
\endgroup
\def\@stepone{\expandafter\futurelet\cs\@steptwo}%
\def\@steptwo{\expandafter\ifx\cs\@stoken\let\@@next=\@stepthree
   \else\let\@@next=\@nexttoken\fi \@@next}%
\def\@stepthree{\afterassignment\@stepone\let\@@next= }%
%
%
%
\def\@getoptionalarg#1{%
   \let\@optionaltemp = #1%
   \let\@optionalnext = \relax
   \@futurenonspacelet\@optionalnext\@bracketcheck
}%
%
%
\def\@bracketcheck{%
   \ifx [\@optionalnext
      \expandafter\@@getoptionalarg
   \else
      \let\@optionalarg = \empty
      \expandafter\@optionaltemp
   \fi
}%
\def\@@getoptionalarg[#1]{%
   \def\@optionalarg{#1}%
   \@optionaltemp
}%
%
%
%
\def\@nnil{\@nil}%
\def\@fornoop#1\@@#2#3{}%
\def\@for#1:=#2\do#3{%
   \edef\@fortmp{#2}%
   \ifx\@fortmp\empty \else
      \expandafter\@forloop#2,\@nil,\@nil\@@#1{#3}%
   \fi
}%
\def\@forloop#1,#2,#3\@@#4#5{\def#4{#1}\ifx #4\@nnil \else
       #5\def#4{#2}\ifx #4\@nnil \else#5\@iforloop #3\@@#4{#5}\fi\fi
}%
\def\@iforloop#1,#2\@@#3#4{\def#3{#1}\ifx #3\@nnil
       \let\@nextwhile=\@fornoop \else
      #4\relax\let\@nextwhile=\@iforloop\fi\@nextwhile#2\@@#3{#4}%
}%
%
%
%
\@innernewif\if@fileexists
\def\@testfileexistence{\@getoptionalarg\@finishtestfileexistence}%
\def\@finishtestfileexistence#1{%
   \begingroup
      \def\extension{#1}%
      \immediate\openin0 =
         \ifx\@optionalarg\empty\jobname\else\@optionalarg\fi
         \ifx\extension\empty \else .#1\fi
         \space
      \ifeof 0
         \global\@fileexistsfalse
      \else
         \global\@fileexiststrue
      \fi
      \immediate\closein0
   \endgroup
}%
%
%
%
%
\def\bibliographystyle#1{%
   \@readauxfile
   \@writeaux{\string\bibstyle{#1}}%
}%
\let\bibstyle = \@gobble
%
%
\let\bblfilebasename = \jobname
\def\bibliography#1{%
   \@readauxfile
   \@writeaux{\string\bibdata{#1}}%
   \@testfileexistence[\bblfilebasename]{bbl}%
   \if@fileexists
      \nobreak
      \@readbblfile
   \fi
}%
\let\bibdata = \@gobble
%
%
\def\nocite#1{%
   \@readauxfile
   \@writeaux{\string\citation{#1}}%
}%
\@innernewif\if@notfirstcitation
%
%
\def\cite{\@getoptionalarg\@cite}%
%
%
\def\@cite#1{%
   \let\@citenotetext = \@optionalarg
   \printcitestart
   \nocite{#1}%
   \@notfirstcitationfalse
   \@for \@citation :=#1\do
   {%
      \expandafter\@onecitation\@citation\@@
   }%
   \ifx\empty\@citenotetext\else
      \printcitenote{\@citenotetext}%
   \fi
   \printcitefinish
}%
\def\@onecitation#1\@@{%
   \if@notfirstcitation
      \printbetweencitations
   \fi
   \expandafter \ifx \csname\@citelabel{#1}\endcsname \relax
      \if@citewarning
         \message{\@linenumber Undefined citation `#1'.}%
      \fi
      \expandafter\gdef\csname\@citelabel{#1}\endcsname{%
\strut
\vadjust{\vskip-\dp\strutbox
\vbox to 0pt{\vss\parindent0cm \leftskip=\hsize 
\advance\leftskip3mm
\advance\hsize 4cm\strut\openup-4pt 
\rightskip 0cm plus 1cm minus 0.5cm ?  #1 ?\strut}}
         {\tt
            \escapechar = -1
            \nobreak\hskip0pt
            \expandafter\string\csname#1\endcsname
            \nobreak\hskip0pt
         }%
      }%
   \fi
   \csname\@citelabel{#1}\endcsname
   \@notfirstcitationtrue
}%
%
%
\def\@citelabel#1{b@#1}%
%
%
\def\@citedef#1#2{\expandafter\gdef\csname\@citelabel{#1}\endcsname{#2}}%
%
%
%
\def\@readbblfile{%
   \ifx\@itemnum\@undefined
      \@innernewcount\@itemnum
   \fi
   \begingroup
      \def\begin##1##2{%
         \setbox0 = \hbox{\biblabelcontents{##2}}%
         \biblabelwidth = \wd0
      }%
      \def\end##1{}
      %
      %
      \@itemnum = 0
      \def\bibitem{\@getoptionalarg\@bibitem}%
      \def\@bibitem{%
         \ifx\@optionalarg\empty
            \expandafter\@numberedbibitem
         \else
            \expandafter\@alphabibitem
         \fi
      }%
      \def\@alphabibitem##1{%
         \expandafter \xdef\csname\@citelabel{##1}\endcsname {\@optionalarg}%
         \ifx\biblabelprecontents\@undefined
            \let\biblabelprecontents = \relax
         \fi
         \ifx\biblabelpostcontents\@undefined
            \let\biblabelpostcontents = \hss
         \fi
         \@finishbibitem{##1}%
      }%
      \def\@numberedbibitem##1{%
         \advance\@itemnum by 1
         \expandafter \xdef\csname\@citelabel{##1}\endcsname{\number\@itemnum}%
         \ifx\biblabelprecontents\@undefined
            \let\biblabelprecontents = \hss
         \fi
         \ifx\biblabelpostcontents\@undefined
            \let\biblabelpostcontents = \relax
         \fi
         \@finishbibitem{##1}%
      }%
      \def\@finishbibitem##1{%
         \biblabelprint{\csname\@citelabel{##1}\endcsname}%
         \@writeaux{\string\@citedef{##1}{\csname\@citelabel{##1}\endcsname}}%
         \ignorespaces
      }%
      %
      %
      \let\em = \bblem
      \let\newblock = \bblnewblock
      \let\sc = \bblsc
      \frenchspacing
      \clubpenalty = 4000 \widowpenalty = 4000
      \tolerance = 10000 \hfuzz = .5pt
      \everypar = {\hangindent = \biblabelwidth
                      \advance\hangindent by \biblabelextraspace}%
      \bblrm
      \parskip = 1.5ex plus .5ex minus .5ex
      \biblabelextraspace = .5em
      \bblhook
      \input \bblfilebasename.bbl
   \endgroup
}%
%
%
\@innernewdimen\biblabelwidth
\@innernewdimen\biblabelextraspace
%
%
%
\def\biblabelprint#1{%
   \noindent
   \hbox to \biblabelwidth{%
      \biblabelprecontents
      \biblabelcontents{#1}%
      \biblabelpostcontents
   }%
   \kern\biblabelextraspace
}%
%
%
%
\def\biblabelcontents#1{{\bblrm [#1]}}%
%
%
\def\bblrm{\rm}%
%
%
\def\bblem{\it}%
%
%
\def\bblsc{\ifx\@scfont\@undefined
              \font\@scfont = cmcsc10
           \fi
           \@scfont
}%
%
%
\def\bblnewblock{\hskip .11em plus .33em minus .07em }%
%
%
\let\bblhook = \empty
%
%
%
\def\printcitestart{[}
\def\printcitefinish{]}
\def\printbetweencitations{, }
\def\printcitenote#1{, #1}
%
%
%
\let\citation = \@gobble
%
%
%
\@innernewcount\@numparams
%
%
\def\newcommand#1{%
   \def\@commandname{#1}%
   \@getoptionalarg\@continuenewcommand
}%
%
%
\def\@continuenewcommand{%
   \@numparams = \ifx\@optionalarg\empty 0\else\@optionalarg \fi \relax
   \@newcommand
}%
%
%
\def\@newcommand#1{%
   \def\@startdef{\expandafter\edef\@commandname}%
   \ifnum\@numparams=0
      \let\@paramdef = \empty
   \else
      \ifnum\@numparams>9
         \errmessage{\the\@numparams\space is too many parameters}%
      \else
         \ifnum\@numparams<0
            \errmessage{\the\@numparams\space is too few parameters}%
         \else
            \edef\@paramdef{%
               \ifcase\@numparams
                  \empty  No arguments.
               \or ####1%
               \or ####1####2%
               \or ####1####2####3%
               \or ####1####2####3####4%
               \or ####1####2####3####4####5%
               \or ####1####2####3####4####5####6%
               \or ####1####2####3####4####5####6####7%
               \or ####1####2####3####4####5####6####7####8%
               \or ####1####2####3####4####5####6####7####8####9%
               \fi
            }%
         \fi
      \fi
   \fi
   \expandafter\@startdef\@paramdef{#1}%
}%
%
%
%
%
\def\@readauxfile{%
   \if@auxfiledone \else 
      \global\@auxfiledonetrue
      \@testfileexistence{aux}%
      \if@fileexists
         \begingroup
            \endlinechar = -1
            \catcode`@ = 11
            \input \jobname.aux
         \endgroup
      \else
         \message{\@undefinedmessage}%
         \global\@citewarningfalse
      \fi
      \immediate\openout\@auxfile = \jobname.aux
   \fi
}%
%
%
\newif\if@auxfiledone
\ifx\noauxfile\@undefined \else \@auxfiledonetrue\fi
%
%
%
%
\@innernewwrite\@auxfile
\def\@writeaux#1{\ifx\noauxfile\@undefined \write\@auxfile{#1}\fi}%
%
%
%
\ifx\@undefinedmessage\@undefined
   \def\@undefinedmessage{No .aux file; I won't give you warnings about
                          undefined citations.}%
\fi
%
%
\@innernewif\if@citewarning
\ifx\noauxfile\@undefined \@citewarningtrue\fi
%
%
%
\catcode`@ = \@oldatcatcode


\def\widestnumber#1#2{}

\def\rm{\fam0 \tenrm}

\def\fakesubhead#1\endsubhead{\bigskip\noindent{\bf#1}\par}


%
%
%

%

\font\textrsfs=rsfs10
\font\scriptrsfs=rsfs7
\font\scriptscriptrsfs=rsfs5

\newfam\rsfsfam
\textfont\rsfsfam=\textrsfs
\scriptfont\rsfsfam=\scriptrsfs
\scriptscriptfont\rsfsfam=\scriptscriptrsfs

\edef\oldcatcodeofat{\the\catcode`\@}
\catcode`\@11

\def\Cal@@#1{\noaccents@ \fam \rsfsfam #1}

\catcode`\@\oldcatcodeofat

\newpage

\head {Annotated Content} \endhead  \resetall
\bn
I\S0 \ub{Introduction}
\mr
\item "{{}}"  [We review background and some 
definitions and theorems on abstract elementary classes.]
\ermn
I\S1 \ub{The Framework}
\mr
\item "{{}}"  [We define types, stability in $\lambda,{\Cal S}(M)$ and
$E_\mu$: equivalence relations on types all whose restrictions to models
of cardinality $\le \mu$ are equal.  We recall that categoricity in $\lambda$
implies stability in $\mu \in [LS({\frak K}),\lambda)$.]
\ermn
I\S2 \ub{Variant of Saturation}
\mr
\item "{{}}"  [We define $<^\ell_{\mu,\alpha}$ and ``$N$ is
$(\mu,\kappa)$-saturated over $M$" and show universality and uniqueness.]
\ermn
I\S3 \ub{Splitting}
\mr
\item "{{}}"  [We note that stability in $\mu$ implies that there are 
not so many $\mu$-splittings.]
\ermn
I\S4 \ub{Indiscernibility and E.M. models}
\mr
\item "{{}}"  [We define strong splitting and dividing, and connect them to
the order property and unstability.]
\ermn
I\S5 \ub{Rank and Superstability}
\mr
\item "{{}}"  [We define one variant of superstability; in particular
categoricity implies it.]
\ermn
I\S6 \ub{Existence of many non-splitting}
\mr
\item "{{}}"  [We prove (e.g. for ${\frak K}$ categorical in $\lambda = 
\text{ cf}(\lambda)$) that if $M_0 <^1_{\mu,\kappa} M_1 \le_{\frak K} 
N \in {\frak K}_{< \lambda}$ and $p \in {\Cal S}(M)$ does not $\mu$-split 
over $M_0$, then $p$ can be extended 
to $q \in {\Cal S}(N)$ which does not $\mu$-split over $M_0$. \nl
(Note: up to $E_\mu$-equivalence the extension is unique).  Secondly, if 
$\langle M_i:i \le \delta \rangle$ is $\le^1_{\mu,\kappa}$-increasing 
continuous in $K_\mu$ and $p \in {\Cal S}(M_\delta)$ then for some $i$ we 
have: $p$ does not $\mu$-split over $M_i$.]
\ermn
I\S7 \ub{More on Splitting}
\mr
\item "{{}}"  [We connect non-splitting to rank and to dividing.]
\ermn
II\S8 \ub{Existence of nice $\Phi$}
\mr
\item "{{}}"  [We try to successively extend the $\Phi$ we use which is proper
for linear orders such that we have as many definable automorphisms as
possible.  We also relook at omitting types theorems over larger model (so
only restrictions will appear).]
\ermn
II\S9 \ub{Small Pieces are Enough and Categoricity}
\mr
\item "{{}}"  [The main claim is that for some not too large $\chi$, if
$p_1,p_2 \in {\Cal S}(M)$ are $E_\chi$-equivalent, $\|M\| < \lambda$ where
$K$ is categorical in $\lambda$ we have $p_1\,E_\chi\,p_2 \Leftrightarrow
p_1 = p_2$. \nl
Lastly, we derive that categoricity is downward closed for successor
cardinals large enough above LS$({\frak K})$.]
\endroster
\bn

\head {\S0 Introduction} \endhead  \resetall
\bigskip

We try to find something on

$$
\text{Cat}_K = \{\lambda:K \text{ categorical in } \lambda\}
$$
\medskip

\noindent
for ${\frak K}$ an abstract elementary class with amalgamation 
(see \scite{0.A} below). \newline
The Los conjecture = Morley theorem deals with the case where $K$ is the 
class of models of a countable first order theory $T$.  See \cite{Sh:c} 
for more on
first order theories.  What for $T$ a theory in an infinitary language?  
(For a theory $T$, $K$ is the class $K_T = \{M:M \models T\}$ we may 
write Cat$_T$ instead of Cat$_{K_T} =$ Cat$_K$).  
Keisler gets what can be gotten from Morley's
proof on $\psi \in L_{\aleph_1,\aleph_0}$.  Then see \cite{Sh:48} on 
categoricity in $\aleph_1$ for $\psi \in L_{\aleph_1,\aleph_0}$ and even 
$\psi \in L_{\aleph_1,\aleph_0}(Q)$, and \cite{Sh:87a}, \cite{Sh:87b} 
on the behaviour in the $\aleph_n$'s.
Makkai Shelah \cite{MaSh:285} proved: if $T \subseteq L_{\kappa,\aleph_0},
\kappa$ a compact cardinal then Cat$_T \cap \{\mu^+:\mu \ge 
\beth_{(2^{\kappa +|T|})^+}\}$ is empty or is $\{\mu^+:\mu \ge 
\beth_{(2^{\kappa+|T|})^+}\}$ (it relies on some developments from 
\cite{Sh:300} but is self-contained).
\medskip

It was then reasonable to deal with weakening the requirement on $\kappa$ to
measurability.  Kolman Shelah \cite{KlSh:362} proved that if
$\mu \in \text{ Cat}_T$, \underbar{then} (after cosmetic changes), for the
right $\le_T$ the class $\{M:M \models T,\|M\| < \lambda\}$ has 
amalgamation and joint embedding property.  This is continued in 
\cite{Sh:472} which gets results on categoricity parallel to the one in 
\cite{MaSh:285} for the ``downward" implication.
\medskip

In \cite{Sh:88} we deal with abstract elementary classes (they include models 
of $T \subseteq L_{\kappa,\aleph_0}$, see \scite{0.A}), prove a 
representation theorem (see \scite{0.E} below), and investigate categoricity
in $\aleph_1$ (and having models in $\aleph_2$, limit models, realizing and
materializing types).  Unfortunately, we do not have anything interesting 
to say here on this context.  So we add amalgamation and the joint embedding
properties thus getting to the framework of Jonsson \cite{J} 
(they are the ones needed to
construct homogeneous universal models).  So this context is more narrow than
the ones discussed above, but we do not use large cardinals.
We concentrate here, for categoricity on $\lambda$, on the case ``$\lambda$ is
regular".   See for later works \cite{Sh:576}, \cite{Sh:600} and
\cite{ShVi:635}. \newline
We quote the basics from \cite{Sh:88} (or \cite{Sh:576}).

We thank Andres Villaveces and Rami Grossberg for much help.
\bigskip

\definition{\stag{0.A} Definition}  ${\frak K} = (K,\le_{\frak K})$ 
is an abstract elementary class \ub{if} for some vocabulary 
$\tau = \tau(K) = \tau({\frak K}),K$ is a class of $\tau(K)$-models, and 
the following axioms hold.
\enddefinition
\bigskip

$Ax0$: The holding of $M \in K,N 
\le_{\frak K} M$ depends on $N,M$ only up to isomorphism i.e. 
$[M \in K, M \cong N \Rightarrow N \in K]$, and [if $N \le_{\frak K} M$ 
and $f$ is an isomorphism from $M$ onto the $\tau$-model $M'$ mapping $N$ 
onto $N'$ \underbar{then} $N' \le_{\frak K} M'$].
\medskip

$AxI$: If $M \le_{\frak K} N$ then $M \subseteq N$ (i.e. $M$ is a submodel
of $N$).
\medskip

$AxII$: $M_0 \le_{\frak K} M_1 \le_{\frak K} M_2$ implies $M_0 \le_{\frak K}
M_2$ and $M \le_{\frak K} M$ for $M \in K$.
\medskip

$AxIII$: If $\lambda$ is a regular cardinal, $M_i$ (for $i < \lambda$) is a
$\le_{\frak K}$-increasing (i.e. $i < j < \lambda$ implies $M_i \le_{\frak K}
M_j$) and continuous (i.e. for limit ordinal $\delta < \lambda$ we have \nl
$M_\delta = \dsize \bigcup_{i < \delta} M_i$) \underbar{then} $M_0
\le_{\frak K} \dsize \bigcup_{i < \lambda} M_i \in {\frak K}$.
\medskip

$AxIV$:  If $\lambda$ is a regular cardinal, $M_i (i < \lambda)$ is
$\le_{\frak K}$-increasing continuous and $M_i \le_{\frak K} N$
\underbar{then} $\dsize \bigcup_{i < \lambda} M_i \le_{\frak K} N$.
\medskip

$AxV$:  If $M_0 \subseteq M_1$ and $M_\ell \le_{\frak K} N$ for $\ell = 0,1$,
\ub{then} $M_0 \le_{\frak K} M_1$.
\medskip

$AxVI$:  $LS({\frak K})$ exists \footnote{We normally assume $M \in {\frak K}
\Rightarrow \|M\| \ge LS({\frak K})$, here there is no loss in it.  It is
also natural to assume $|\tau({\frak K})| \le LS({\frak K})$ which just means
increasing $LS({\frak K})$.}; see below Definition \scite{0.C}.
\bigskip

\definition{\stag{0.B} Definition}  1) $K_\mu =: \{M \in K:\|M\| = \mu\}$.
\newline
2) We say $h$ is a $\le_{\frak K}$-embedding of $M$ into $N$ is for some 
$M' \le_{\frak K} N,h$ is an isomorphism from $M$ onto $M'$.
\enddefinition
\bigskip

\definition{\stag{0.C} Definition}  1) We say that $\mu$ is a Skolem Lowenheim
number of ${\frak K}$ \ub{if} \nl
$\mu \ge \aleph_0$ and:
\medskip
\roster
\item "{$(*)^\mu_K$}"  for every $M \in K,A \subseteq M,|A| \le \mu$ there
is $M',A \subseteq M' \le_{\frak K} M$ and \nl
$\|M'\| \le \mu$.
\endroster
\medskip

\noindent
2) $LS'({\frak K}) = \text{ Min}\{\mu:\mu 
\text{ is a Skolem Lowenheim number of } {\frak K}\}$. \nl
3) $LS({\frak K}) = LS'({\frak K}) + |\tau(K)|$. 
\enddefinition
\bigskip

\proclaim{\stag{0.D} Claim}  1) If $I$ is a directed partial order, $M_t \in
K$ for $t \in I$ and \nl
$s <_I t \Rightarrow M_s \le_{\frak K} M_t$ \ub{then}
\medskip
\roster
\item "{$(a)$}"  $M_s \le_{\frak K} \dsize \bigcup_{t \in I} M_t \in K$ 
for every $s \in I$
\smallskip
\noindent
\item "{$(b)$}"  if $(\forall t \in I)[M_t \le_{\frak K} N]$ then 
$\dsize \bigcup_{t \in I} M_t \le_{\frak K} N$.
\endroster
\medskip

\noindent
2) If $A \subseteq M \in K,|A| + LS'({\frak K}) \le \mu \le \|M\|$, 
\underbar{then} there is $M_1 \le_{\frak K} M$ such that \newline
$\|M_1\| = \mu$ and $A \subseteq M_1$. \newline
3) If $I$ is a directed partial order, $M_t \le N_t \in K$ for $t \in I$
and $s \le_I t \Rightarrow M_s \le_{\frak K} M_t \and N_s \le_{\frak K} N_t$ 
\ub{then} $\dsize \bigcup_t M_t \le_{\frak K} \dsize \bigcup_t N_t$.
\endproclaim
\bigskip

\proclaim{\stag{0.E} Claim}  Let ${\frak K}$ be an abstract elementary class.
There are $\tau^+,\Gamma$ such that:
\medskip
\roster
\item "{$(a)$}"  $\tau^+$ is a vocabulary extending $\tau(K)$ of cardinality
$LS({\frak K})$
\smallskip
\noindent
\item "{$(b)$}"  $\Gamma$ is a set of quantifier free types in $\tau^+$
(each is an $m$-type for some $m < \omega$)
\smallskip
\noindent
\item "{$(c)$}"  $M \in K$ iff for some $\tau^+$-model $M^+$ omitting every
$p \in \Gamma$ we have \newline
$M = M^+ \restriction \tau$
\smallskip
\noindent
\item "{$(d)$}"  $M \le_{\frak K} N$ \underbar{iff} there are $\tau^+$-models
$M^+,N^+$ omitting every $p \in \Gamma$ such that $M^+ \subseteq N^+,M =
M^+ \restriction \tau(K)$ and $N = N^+ \restriction \tau(K)$. \nl
We can replace $M \le_{\frak K} N$ by $M_i \le_{\frak K} N$ for a 
family $\{M_i:i \in I\}$ (getting $\tau^+$-expansions $M^+_i,N^+$ of $M_i,N$
respectively, such that $M^+_i \subseteq N^+$ and $M^+_i,N^+$ omit every
$p \in \Gamma$ for every $i \in I$) \ub{if} for any
$\bar a \in {}^{\omega >}M$ for some $i,\bar a \in {}^{\omega >}(M_i)$ and
$\bar a \in {}^{\omega >}(M_j) \Rightarrow M_i \subseteq M_j$
\sn
\item "{$(e)$}"  if $M \le_{\frak K} N$ and $M^+$ is an expansion of $M$ to
a $\tau^+$-model omitting every $p \in \Gamma$ \ub{then} we can find a
$\tau^+$-expansion of $N$ omitting every $p \in \Gamma$ such that $M^+
\subseteq N^+$.
\endroster
\endproclaim
\bigskip

\proclaim{\stag{0.F} Claim}  Assume ${\frak K}$ has a member of 
cardinality $\ge \beth_{(2^{LS({\frak K})})^+}$ (here and elsewhere we 
can weaken this to: 
has a model of cardinality $\ge \beth_\alpha$ for every 
$\alpha < (2^{LS({\frak K})})^+$).  
\underbar{Then} there is $\Phi$ proper for linear orders 
(see \cite[Ch.VII,\S2]{Sh:c}) such that:
\medskip
\roster
\item "{$(a)$}"  $|\tau(\Phi)| = LS({\frak K})$
\smallskip
\noindent
\item "{$(b)$}"  for linear orders $I \subseteq J$ we have \newline
$EM_\tau(I,\Phi) \le_{\frak K} EM(J,\Phi)(\in K)$.
\smallskip
\noindent
\item "{$(c)$}"  $EM_\tau(I,\Phi)$ has cardinality $|I| + LS({\frak K})$ 
(so ${\frak K}$ has a model in every cardinality $\ge LS({\frak K}))$.
\endroster
\endproclaim
\newpage

\head {PART 1 \\
\S1 The Framework} \endhead  \resetall
\bigskip

\demo{\stag{1.0} Hypothesis}
\roster
\item "{$(a)$}"  ${\frak K} = (K,\le_{\frak K})$  an abstract elementary 
class (\scite{0.A}) so \nl
$K_\lambda  = \{M \in  K:\|M\| = \lambda \}$
\smallskip
\noindent
\item "{$(b)$}"  ${\frak K}$  has amalgamation 
and the joint embedding property      
\smallskip
\noindent
\item "{$(c)$}"  $K$  has members of arbitrarily large cardinality,
equivalently: $K$ has a member of cardinality at least
$\beth_{(2^{LS({\frak K})})^+}$.
\endroster
\enddemo
\bigskip

\demo{\stag{1.0A} Convention}  1) So there is a monster ${\frak C}$ (see
\cite[Ch.I,\S1]{Sh:a} = \cite[Ch.I,\S1]{Sh:c}).
\enddemo
\bigskip

\definition{\stag{1.1} Definition}  We say $K$ (or ${\frak K}$) is 
categorical in $\lambda$ \ub{if} it has
one and only one model of cardinality $\lambda$, up to isomorphism.
\enddefinition
\bigskip

\definition{\stag{1.2} Definition}  1) We can define tp$(\bar a,M,N)$ (when
$M \le_{\frak K} N$ and $\bar a \subseteq N$), 
as $(\bar a,M,N)/E$ where $E$ is the following equivalence relation:
$(\bar a^1,M^1,N^1)\,E\,(\bar a^2,M^2,N^2)$ \underbar{iff} $M^\ell
\le_{\frak K} N^\ell,\bar a^\ell \in {}^\alpha(N^\ell)$ (for some $\alpha$)
and $M^1 = M^2$ and there is $N \in K$ satisfying $M^1 = M^2 \le_{\frak K}
N$ and $\le_{\frak K}$-embedding $f^\ell:N^\ell \rightarrow N$ over $M^\ell$
(i.e. $f \restriction M^\ell$ is the identity) for $\ell=1,2$ and 
$f^1(\bar a^1) = f^2(\bar a^2)$. \newline
2)  We omit $N$ when $N = {\frak C}$ (see \scite{1.0A}) and may then write 
$\frac{\bar a}{M} = \text{ tp}(\bar a,M,{\frak C})$.  
We can define $N$ is $\kappa$-saturated 
(when $\kappa > LS({\frak K})$) by: if $M \le_{\frak K} N,\|M\| < \kappa$
and $p \in {\Cal S}^{< \omega}(M)$ (see below) \ub{then} $p$ is realized 
in $M$, i.e. for some $\bar a \subseteq N,p = \text{ tp}(\bar a,M,N)$. \nl
3) ${\Cal S}^\alpha(M) = \{\text{tp}(\bar a,M,N):\bar a \in {}^\alpha N,
M \le_{\frak K} N\}$; we define $p \restriction M$ when \nl
$M \le_{\frak K} N 
\and p \in {\Cal S}(N)$ as tp$(\bar a,M,N_1)$ when $N \le_{\frak K} N_1,
p = \text{ tp}(\bar a,N,N_1)$.  Let \nl
$p \le q$ mean $p \in {\Cal S}(M),q \in
{\Cal S}(N),p = q \restriction M$; see \cite[Ch.II]{Sh:300} or 
\cite[\S0]{Sh:576}. \newline
4) ${\Cal S}(M) = {\Cal S}^1(M)$ (could just as well use ${\Cal S}^{< \omega}
(M) = \dbcu_{n < \omega} {\Cal S}^n(M)$). \nl
5) If $M_0 \le_{\frak K} M_1$ and $p_\ell \in S^\alpha(M_\ell)$ for
$\ell =1,2$, \ub{then} $p_0 = p_1 \restriction M_0$ means that for some
$\bar a,N$ we have $M_1 \le_{\frak K} N$ and $\bar a \in {}^\alpha N$ and
$p_\ell = \text{ tp}(\bar a,M_\ell,N)$ for $\ell=1,2$.
\enddefinition
\bigskip

\definition{\stag{1.3} Definition}  Let ${\frak K}$ stable in $\lambda$ mean: 
$\|M\| \le \lambda \Rightarrow |{\Cal S}(M)| \le \lambda$ and $\lambda \ge 
LS({\frak K})$.
\enddefinition
\bigskip

\demo{\stag{1.4} Convention}  If not said otherwise, $\Phi$ is as in
\scite{0.F}.
\enddemo
\bigskip

\proclaim{\stag{1.5} Claim}  If $K$ is categorical in $\lambda$ and
$\lambda \ge LS({\frak K})$, \underbar{then} 
\medskip
\roster
\item "{$(a)$}"  ${\frak K}$ is stable in every $\mu$ which satisfifes
$LS({\frak K}) \le \mu < \lambda$, hence
\smallskip
\noindent
\item "{$(b)$}"  the model $M \in K_\lambda$ is cf$(\lambda)$-saturated (if
cf$(\lambda) > LS({\frak K})$).
\endroster
\endproclaim
\bigskip

\demo{Proof}  Like \cite{KlSh:362}.
\enddemo
\bigskip

\definition{\stag{1.6} Definition}  $E_\mu$ is the following relation,
$$
\align
p\,\,E_\mu\,\,q \text{ iff } &\text{for some } M \in K,m < \omega 
\text{ we have} \\
  &p,q \in {\Cal S}^m(M) \text{ and } 
[N \le_{\frak K} M \and \|N\| \le \mu \Rightarrow p \restriction N = q
\restriction N].
\endalign
$$
\medskip

\noindent
Obviously it is an equivalence relation.
\enddefinition
\bigskip

\remark{\stag{1.7} Remark}  1) In previous contexts $E_{LS({\frak K})}$ is 
equality, e.g. the axioms of $NF$ in \newline
\cite[Ch.II,\S1]{Sh:300} implies it; but here we do not know ---
this is the main difficulty. \newline
We may look at this as our bad luck, or inversely, a place to encounter some
of the difficulty of dealing with $L_{\mu,\omega}$ (in which our context is
included). \nl
2) In the cases we shall deal with 
we can define ``$M \in K_{LS({\frak K})}$" is saturated.
\endremark
\bigskip

\proclaim{\stag{1.8} Claim}  1) There is no maximal member in $K$, in fact for
every $M \in K$ there is $N,M <_{\frak K} N \in K,\|N\| \le \|M\| + 
LS({\frak K})$, hence for every $\lambda \ge \|M\| + LS({\frak K})$ there is
$N \in K_\lambda$ such that $M <_{\frak K} N \in K_\lambda$. \nl
2) If $p_2 \in {\Cal S}^\alpha(M_2)$ and $M_1 \le_{\frak K} M_2 \in K$
\ub{then} for one and only one $p_1 \in {\Cal S}^\alpha(M_1)$ we have
$p_1 = p_2 \restriction M_1$. \nl
3) If $p_1 \in {\Cal S}^\alpha(M_1)$ and $M_1 \le_{\frak K} M_2 \in K$
\ub{then} for some $p_2 \in {\Cal S}^\alpha(M_2)$ we have
$p_1 = p_2 \restriction M_1$. \nl 
4) If $M_1 \le_{\frak K} M_2 \le_{\frak K} M_3$ and
$p_\ell \in {\Cal S}^\alpha(M_\ell)$ for $\ell =1,2,3$ \ub{then}
$p_3 \restriction M_2 = p_2 \and p_2 \restriction M_1 = p_1 \Rightarrow
p_3 \restriction M_1 = p_2$. \nl
\endproclaim
\bigskip

\demo{Proof}  1) Immediate by clause (c) of the hypothesis \scite{1.0} and
claim \scite{0.F}. \nl
2) Straightforward. \nl
3) By amalgamation. \nl
4) Check.  \hfill$\square_{\scite{1.8}}$
\enddemo
\bigskip

\proclaim{\stag{1.9} Claim}  If $\langle M_i:i \le \omega \rangle$ is
$\le_{\frak K}$-increasing continuous and $p_n \in {\Cal S}^\alpha(M_n)$
and $p_n = p_{n+1} \restriction M_n$ for $n < \omega$, \ub{then} there is
$p_\omega \in {\Cal S}^\alpha(M_\omega)$ such that $n < \omega \Rightarrow
p_\omega \restriction _n = p_n$.
\endproclaim
\bigskip

\demo{Proof}  Straight chasing diagrams.
\enddemo
\bigskip

\remark{\stag{1.10} Remark}  In \scite{1.9} we do not claim uniqueness and
not existence replacing $\omega$ for $\delta$ of uncountable cofinality.
In general not true [Saharon add].
\endremark
\newpage

\head {\S2 Variant of Saturated} \endhead  \resetall
\bigskip

\definition{\stag{2.1} Definition}  Assuming ${\frak K}$ stable in 
$\mu$ and $\alpha$
is an ordinal $< \mu^+$, $\mu^+ \times \alpha$ means ordinal product. 
\newline
1)  $M <^\circ_{\mu,\alpha} N$ if: $M \in K_\mu,N \in K_\mu,M 
\le_{\frak K} N$ and there is a $\le_{\frak K}$-increasing sequence 
$\bar M = \langle M_i:i \le \mu \times \alpha \rangle$ which is continuous, 
$M_0 = M,M_{\mu \times \alpha} \le_{\frak K} N$ 
and every $p \in {\Cal S}^1(M_i)$ is realized in $M_{i+1}$. \newline
2)  We say $M <^1_{\mu,\alpha} N$ iff $M \in K_\mu,N \in K_\mu,
M \le_{\frak K} N$ and there is a $\le_{\frak K}$-increasing sequence 
$\bar M = \langle M_i:i \le \mu \times \alpha 
\rangle,M_0 = M,M_{\mu \times \alpha} = N$ and every $p \in {\Cal S}^1(M_i)$  
is realized in $M_{i+1}$.  \newline
3) If $\alpha = 1$, we may omit it.
\enddefinition
\bigskip

\proclaim{\stag{2.2} Lemma}  Assume ${\frak K}$ stable in 
$\mu$ and $\alpha < \mu^+$.
\newline
0) If $\ell \in \{0,1\}$ and $\alpha_1 < \alpha_2 < \mu^+$ and there is
$b \subseteq \alpha_2$ such that otp$(b) = \alpha_1$ \newline
and $[\ell = 1 \Rightarrow b$ unbounded in $\alpha_2]$ \ub{then} 
$<^\ell_{\mu,\alpha_2} \subseteq <^\ell_{\mu,\alpha_1}$. \newline
1) If $M \in  K_\mu$, \ub{then} for some $N$ we have $M <^\circ_{\mu,\alpha} 
N$ and for some $N,M <^1_{\mu,\alpha} N$. \newline
2) $\quad$ (a) \,\, If  $M \in K_\mu,M \le_{\frak K} M' \le^\ell_{\mu,\alpha} 
N$ \ub{then} $M \le^\ell_{\mu,\alpha} N$. \newline

$\quad$ (b) \,\, If  $M \in K_\mu,M \le_{\frak K} M' \le^\circ_{\mu,\alpha}
N' \le_{\frak K} N \in K_\mu$ \ub{then} $M \le^\circ_{\mu,\alpha} N$. \newline
3)  If $\langle M_i:i < \alpha \rangle$ is $\le_{\frak K}$-increasing sequence
in $K_\mu,M_i \le^\circ_\mu M_{i+1}$ and $\alpha < \mu^+$ is a
limit ordinal, \underbar{then} $M_0 \le^1_{\mu,\alpha} 
\dsize \bigcup_{i < \alpha} M_i$. \newline
4)  If $M \le^\circ_\mu N$  \underbar{then}:
\roster
\item "{$(a)$}"  any  $M' \in K_\mu$ can be $\le_{\frak K}$-embedded into $N$
(here we can waive $\|M\| = \mu$)
\smallskip
\noindent
\item "{$(b)$}"  If  $M' \le_{\frak K} N' \in K_{\le \mu}$, $h$ is a
$\le_{\frak K}$-embedding of $M'$ into $M$ then $h$ can be extended to
a $\le_{\frak K}$-embedding of $N'$ into $N$.
\endroster
\medskip

\noindent
5)  If  $M^\ell \le^1_{\mu,\kappa} N^\ell$ for  $\ell = 1,2$,  $h$ 
an isomorphism from  $M^1$ into [onto]  $M^2$ then  $h$  can be extended
to an isomorphism from  $N^1$ into [onto]  $N^2$.\newline 
6)  If  $M \le^1_{\mu,\kappa} N^\ell$ for $\ell = 1,2$ \ub{then} 
$N^1 \cong N^2$ (even over $M$).\newline
7)  If  $M \le^\circ_{\mu,\kappa} N$, $M \le_{\frak K} M' \in K_\mu$
\underbar{then} $M'$ can be $<_{\frak K}$-embedded into $N$ over $M$.\newline 
8)  If  $\mu \ge \kappa > LS({\frak K})$ and $M <^1_{\mu,\kappa} N$ 
\ub{then} $N$ is cf$(\kappa)$-saturated.
\endproclaim
\bigskip

\demo{Proof}  See \cite[Ch.II,3.10,p.319]{Sh:300} and around, we shall explain
and prove part (8) below.
\enddemo
\bigskip

\noindent
\underbar{\stag{2.2A} Discussion}:  There (in \cite[Ch.II,3.6]{Sh:300})
the main point was that for $\kappa > LS({\frak K})$, the notions 
``$\kappa$-homogeneous universal" and $\kappa$-saturation (i.e. 
every ``small" 1-type is realized) are equivalent.
\medskip

Not hard, still \cite[Ch.II,3.6]{Sh:300} was a surprise to some.  
In first order the equivalence saturated $\equiv$ homogeneous universal for
$\prec$ seemed a posteriori natural as the homogeneity used was anyhow for
sequences of elements realizing the same first order formulas so 
(forgetting about the models) to some extent this seemed natural; 
i.e. asking this for any type of 1-element was very natural.

But here, types of 1-element are really meaningful only over a model.  So it
seems that if over any small submodel every type of 1-element is realized (say
in ${\frak A}$) and we want to embed $N \ge_{\frak K} N_0,N_0 \le_{\frak K} 
{\frak A}$ into ${\frak A}$ over $N_0$, we encounter the following problem:
we cannot continue this as after $\omega$ stages, as we get a set which is 
not a model (if $LS({\frak K}) > \aleph_0$ this absolutely necessarily fails;
and if $LS({\frak K}) = \aleph_0$ at best the situation is 
as in \cite{Sh:87a}).

This explains a natural preconception making you not believe; i.e.
psychological barrier to prove. It does not mean that the proof is hard.
\bigskip

\remark{\stag{2.2M} Remark}  Note that $\le^1_{\mu,\kappa},\kappa$ regular are 
the interesting ones.
Still $\le^0_{\mu,\kappa}$ is enough for universality (\scite{2.2}(4)) and is
natural, $\le^1_{\mu,\kappa}$ is natural for uniqueness.  BUT
$<^1_{\mu,\aleph_0} = <^1_{\mu,\aleph_1}$ can be proved only under
categoricity (or something like superstability assumptions).  LOOK at first
order $T$ stable in $\mu$.  Then, $M <^1_{\mu,\kappa} N$ is equivalent to:

$$
\|M\| = \|N\| = \mu,M,N \models T
$$
\medskip

\noindent
and there is $\langle M_i:i \le \kappa \rangle$ which is 
$\prec$-increasing continuous such that

$$
M_0 = M \qquad \qquad M_\kappa = N
$$

$$
(M_{i+1},c)_{c \in M_i} \text{ is saturated}.
$$
\medskip

\noindent
\underbar{Question}:  Now, is $N$ saturated when 
$M <^1_{\mu,\kappa} N$?
\medskip

\noindent
\underbar{Answer}:  It is iff cf$(\kappa) \ge \kappa_r(T)$.  See
\cite[Ch.III,\S3]{Sh:c}.
\smallskip

\noindent
See on limit and superlimit models in \cite{Sh:88}.
\endremark
\bn
\underbar{Before we prove \scite{2.2}(8), recall}
\definition{\stag{2.2B} Definition}  $M \in K$ is $\kappa$-saturated if
$\kappa > LS({\frak K})$ and: \nl
$N \le_{\frak K} M,\|N\| < \kappa,
p \in {\Cal S}^1(N) \Rightarrow p$ realized in $M$.
\enddefinition
\bigskip

\demo{Proof of \scite{2.2}(8)} 
\sn
\underbar{Statement}:  If $M <^1_{\mu,\kappa} N$ ($\kappa$ regular) 
\ub{then} $N$ is $\kappa$-saturated.
\smallskip
\noindent
\underbar{Note}:  if $\kappa \le LS({\frak K})$ the conclusion is essentially
empty, but there is no need for the assumption ``$\kappa > LS({\frak K})$".
\enddemo
\bigskip

\demo{Proof}  Let $\bar M = \langle M_i:i \le \mu \times \kappa \rangle$
witness $M \le^1_{\mu,\kappa} N$ so $M_0 = M,M_{\mu \times \kappa} = N,
M_i \le_{\frak K}$-increasing continuous and every $p \in {\Cal S}(M_i)$ 
is realized in $M_{i+1}$.

Assume
\roster
\item "{$(*)$}"  $N' \le_{\frak K} N,\|N'\| < \kappa,p \in {\Cal S}(N')$.
\endroster
\medskip

\noindent
We should prove that ``$p$ is realized in $N$".  But 
$\langle M_i:i \le \mu \times \kappa \rangle$ is increasing continuous

$$
\text{cf}(\mu \times \kappa) = \kappa > \|N'\|
$$
\medskip

\noindent
so $N' \le_{\frak K} M_{\mu \times \kappa} = 
\dsize \bigcup_{i < \mu \times \kappa}
M_i$ implies there is $i(*) < \mu \times \kappa$, such that $N' \subseteq
M_{i(*)}$ hence by Axiom V we have $N' \le_{\frak K} M_{i(*)}$.  So $p$ 
has (by amalgamation!) an extension $p^* \in {\Cal S}(M_{i(*)})$ and 
$p^*$ is realized in 
$M_{i(*)+1}$ so in $M_{\mu \times \kappa} = N$. \hfill$\square_{\scite{2.2}}$
\enddemo
\bigskip

\noindent
\underbar{Comment}:  Hence length $\mu$ (instead of $\mu \times \kappa$)
suffices.

But for the uniqueness it does not.  See \scite{2.2}(4) + (5).
\bigskip

\noindent
\underbar{Comment}:  The definition of $\le^0_{\mu,\kappa},\le^1_{\mu,\kappa}$
is also essentially taken from \newline
\cite[Ch.II,3.10]{Sh:300}.  We need the intermediate steps to construct models 
so we have to have $\mu$ of them in order to deal with all the elements.
\bigskip

\proclaim{\stag{2.3} Claim}  If $K$ is categorical in $\lambda,M \in K_\lambda$
and cf$(\lambda) > \mu$ \ub{then}: \nl
if $N <_{\frak K} M \in K_\lambda,
N \in K_\mu,N' <_{\frak K} M,h$ an isomorphism from $N$ onto $N'$, 
\underbar{then} $h$ can be extended to an automorphism of $M$.
\endproclaim
\bigskip

\demo{Proof}  By \scite{1.2} we have $LS({\frak K}) \le \mu < \lambda
\Rightarrow {\frak K}$ stable in $\lambda$.  We can find 
$\langle M_i:i < \lambda \rangle$ which is
$<_{\frak K}$-increasing continuous, 
$\|M_i\| = |i| + LS({\frak K})$, \nl
$M_i <^1_{|i|+LS({\frak K}),
|i|+LS({\frak K})}M_{i+1}$.
By the categoricity assumption without loss of generality $M =
\dsize \bigcup_{i < \lambda} M_i$.  As cf$(\lambda) > \mu$ for some $i_0 <
\lambda$ we have $N,N' \prec M_{i_0}$.

By \scite{2.2} we can build an automorphism. \hfill$\square_{\scite{2.3}}$
\enddemo
\bigskip

\definition{\stag{2.4} Definition}  For $\mu \ge LS({\frak K})$, we say 
$N \in K_\mu$ is $(\mu,\kappa)$-saturated \ub{if} for some $M$ we have 
$M <^1_{\mu,\kappa} N$ (so $\kappa$ is $\le \mu$, 
normally regular).
\enddefinition
\bigskip

\proclaim{\stag{2.4A} Claim}  1) The $(\mu,\kappa)$-saturated model is unique 
(even over $M$) if it exists at all. \newline
2) If $M$ is $(\mu,\kappa)$-saturated, $\kappa = \text{{\rm cf\/}}(\kappa)$,
$\text{{\rm cf\/}}(\kappa) > LS({\frak K})$ \ub{then} $M$ is 
$\kappa$-saturated. \newline
3) If $M$ is $(\mu,\kappa)$-saturated for every $\kappa = \text{{\rm cf\/}}
(\kappa) \le \mu$ and $\mu > LS({\frak K})$ \ub{then} $M$ is $\mu$-saturated.
\endproclaim
\bigskip

\noindent
\underbar{Discussion}:  It is natural to define saturated as
$\|M\|$-saturated.  ($I$ may have confusions using the other being
$(\mu,\kappa)$-saturated for every regular $\kappa \le \mu$.)  This is
particularly reasonable when the cardinal is regular, e.g. if $K$ categorical
in $\lambda,\lambda = \text{ cf}(\lambda)$ the model in $K_\lambda$ is
$\lambda$-saturated.

Part of the program is to prove that all the definitions are equivalent.

For now in Definition \scite{2.4} we are not sure that such a model exists.
\newpage

\head {\S3 Splitting} \endhead  \resetall
\bn
Whereas non-forking is very nice in \cite{Sh:c}, in more general 
contexts, non first order, it
is not clear whether we have so good a notion, hence we go back to earlier
notions from \cite{Sh:3}, like splitting.  It still gives for many cases
$p \in {\Cal S}(M)$, a ``definition" of $p$ over some ``small"
$N \le_{\frak K} M$.  We need $\mu$-splitting because $E_{LS({\frak K})}$ is
not known to be equality (see \scite{1.6}).
\bigskip

\demo{\stag{3.0} Context}  Inside the monster model ${\frak C}$.
\enddemo
\bigskip

\definition{\stag{3.1} Definition}  $p \in {\Cal S}(M)$ does $\mu$-split over  
$N \le_{\frak K} M$ if:
\medskip
\roster
\item "{{}}"  $\|N\| \le \mu$, and there are $N_1,N_2,h$ such that:
\newline
$\|N_1\| = \|N_2\| \le \mu$ and $N \le_{\frak K} N_\ell \le_{\frak K} M$, for
$\ell =1,2$ \nl
$h$ an elementary mapping from $N_1$ onto $N_2$ over $N$ such that
\newline
the types $p \restriction N_2$ and $h(p\restriction N_1)$ are contradictory
and $N \le_{\frak K} N_\ell \le_{\frak K} M$.
\endroster
\enddefinition
\bigskip

\proclaim{\stag{3.2} Claim}  1) Assume ${\frak K}$ is stable in $\mu,
\mu \ge LS({\frak K})$.  
If $M \in {\frak K}_{\ge \mu}$ and $p \in {\Cal S}^1(M)$, \ub{then} for some 
$N_0 \subseteq M,\|N_0\| = \mu,p$ does not $\mu$-split over $N_0$ 
(see Definition \scite{3.1}). \newline
2) Moreover, if $2^\kappa > \mu,\langle M_i:i \le \kappa + 1 \rangle$ is
$<_{\frak K}$-increasing, $\bar a \in {}^m(M_{\kappa +1})$, \newline
tp$(\bar a,M_{i+1},M_{\kappa +1})$ does $(\le \mu)$-split over $M_i$,
\underbar{then} ${\frak K}$ is not stable in $\mu$.
\endproclaim
\bigskip

\demo{Proof of \scite{3.2}}  1) If not, we can choose by induction on  
$i < \mu$ \,\, $N_i,N^1_i,N^2_i,h_i$ such that:  
\medskip
\roster
\item "{$(a)$}"  $\langle N_i:i \le \mu \rangle$  is increasing continuous,
$N_i <_{\frak K} M$,  $\|N_i\| = \mu$ 
\smallskip
\noindent
\item "{$(b)$}"  $N_i \le_{\frak K} N^\ell_i \le_{\frak K} N_{i+1}$
\smallskip
\noindent
\item "{$(c)$}"  $h_i$ is an elementary mapping from $N^1_i$ onto $N^2_i$ 
over $N_i$,
\smallskip
\noindent
\item "{$(d)$}"  $p \restriction N^2_i,h_i(p \restriction N^1_i)$  
are contradictory, equivalently distinct (we could have defined them for  
$i < \mu^+)$.
\ermn
Let $\chi = \text{ Min}\{ \chi:2^\chi > \mu \}$ so $2^{<\chi} \le \mu$.
Now contradict stability in $\mu$ as in part (2). \newline
2) Similar to \cite[Ch,I,\S2]{Sh:a} or \cite[Ch.I,\S2]{Sh:c} (by using models),
but we give details.  
Without loss of generality $M_i \in K_{\le \mu}$ for $i \le \kappa +1$.
For each $i < \kappa$ let $N_{i,1},N_{i,2}$ be such that $M_i \le_{\frak K}
N_{i,\ell} \le_{\frak K} M_{i+1},g_i$ an isomorphism from $N_{i,1}$ 
onto $N_{i,2}$ over $M_i$ and tp$(\bar a,N_{i,2}) \ne g_i(\text{tp}
(\bar a,N_{i,1}))$.  Without loss of generality
$2^{< \kappa} \le \mu$.  We define by induction on $\alpha \le \kappa$ a
model $M^*_\alpha$ and for each $\eta \in {}^\alpha 2$, a mapping $h_\eta$
such that:
\mr
\item "{$(a)$}"  $M^*_\alpha \in K_\mu$ is $\le_{\frak K}$-increasing
continuous
\sn
\item "{$(b)$}"  for $\eta \in {}^\alpha 2,h_\eta$ is a
$\le_{\frak K}$-embedding of $M_\alpha$ into $M^*_\alpha$
\sn
\item "{$(c)$}"  if $\beta < \alpha,\eta \in {}^\alpha 2$, \ub{then}
$h_{\eta \restriction \beta} \subseteq h_\eta$
\sn
\item "{$(d)$}"  if $\alpha = \beta +1,\nu \in {}^\beta 2$, \ub{then}
$h_{\nu \char 94 <0>}(N_{i,1}) = h_{\nu \char 94 <1>}(N_{i,2})$.
\ermn
There is no problem to carry the definition (we are using amalgamation only
in $K_{\le \mu}$ and if we start with $M_0 \in K_\mu$ only in $K_\mu$).
Now for each $\eta \in {}^\kappa 2$ we can find $M^*_\eta \in K_\mu,
M^*_\kappa \le_{\frak K} M^*_\eta$ and $\le_{\frak K}$-embedding $h^+_\eta$
of $M_{\kappa +1}$ into $M^*_\eta$ extending $h_\eta = \dbcu_{\alpha < \kappa}
h_{\eta \restriction \alpha}$.  Now $\{\text{tp}(h^+_\eta(\bar a),M^*_\kappa,
M^*_\eta):\eta \in {}^\kappa 2\}$ is a family of $2^\kappa > \mu$ distinct
members of ${\Cal S}^m(M^*_\kappa)$ and recall $M^*_\kappa \in K_\mu$ so we
are done.   \hfill$\square_{\scite{3.2}}$
\enddemo
\bigskip

\demo{\stag{3.3} Conclusion}  If $p \in {\Cal S}^m(M),M$ is $\mu^+$-saturated,
$\kappa = \text{ cf}(\kappa) \le \mu$, \underbar{then} for some  
$N_0 <^\circ_{\mu,\kappa} N_1 \le_{\frak K} M$, (so $\|N_1\| = \mu$) we have: 
\newline
$p$ is the $E_\mu$-unique extension of $p \restriction N_1$ which does not 
$\mu$-split over $N_0$.
\enddemo
\newpage

\head {\S4 Indiscernibles and E.M. Models} \endhead  \resetall
\bigskip

\definition{\stag{4.1} Definition}  Let $h_i:Y \rightarrow {\frak C}$ for
$i < i^*$. \newline
1)  $\langle h_i:i < i^* \rangle$ is an indiscernible sequence (of character
$< \kappa$) (over $A$) if \underbar{for every} $g$, a partial one to one order
preserving map from  $i^*$ to $i^*$ (of cardinality $< \kappa$)
\underbar{there is} $f \in AUT({\frak C})$, such that

$$
g(i) = j \Rightarrow  h_j \circ  h^{-1}_i \subseteq f
$$

$$
\text{(and \, id}_A \subseteq f).
$$
\medskip

\noindent
2)  $\langle h_i:i < i^* \rangle$ is an indiscernible set (of character
$\kappa$) (over $A$) if:
for every $g$ partial one to one map from $i^*$ to $i^*$ (with
$|\text{Dom }g| \le \kappa$) \underbar{there is} $f \in  AUT({\frak C})$, 
such that

$$
g(i) = j \Rightarrow h_j \circ h^{-1}_i \subseteq f
$$

$$
\text{(and \, id}_A \subseteq f).
$$
\medskip

\noindent
3)  $\langle h_i:i < i^\ast \rangle$ is a strictly indiscernible sequence, if
$i^* \ge \omega$ and for some $\Phi$, proper for linear orders (see 
\cite[Ch.VII]{Sh:a} or \cite[Ch.VII]{Sh:c}) in vocabulary $\tau_1 = \tau(\Phi)$
extending $\tau(K)$, there is $M^1 = EM^1(i^*,\Phi)$ such that $M^1$ is the 
Skolem Hull of $\{x_i:i < i^*\}$, and a sequence of unary terms  
$\langle \sigma_t:t \in Y \rangle$ such that:

$$
\sigma_t(x_i) = h_i(t) \text{ for } i < i^\ast, t \in  Y
$$

$$
M^1\restriction \tau(K) <_{\frak K} {\frak C}.
$$
\mn
4) Let $h_i:Y_i \rightarrow {\frak C}$ for $i < i^*$ we say that
$\langle h_i:i < i^* \rangle$ has characteristic $\sigma$ if:
\mr
\item "{$(*)$}"  if $h'_i:Y_i \rightarrow {\frak C}$ for $i < i^*$ and for
every $u \in [i^*]^{< \sigma}$ there is an automorphism $f_u$ of ${\frak C}$
such that $f_u \restriction A = \text{ id}_A$ and $i \in u \Rightarrow f_u
\circ h_i = h'_i$, \ub{then} there is an automorphism $f$ of ${\frak C}$ such
that $f \restriction A = \text{ id}_A$ and 
$i < i^* \rightarrow f \circ h_i=h'_i$.
\endroster
\enddefinition
\bigskip

\demo{\stag{4.1A} Notation}  We can replace $h_i$ by the sequence
$\langle h_i(t):t \in  Y \rangle$.
\enddemo
\bigskip

\definition{\stag{4.1B} Definition}  1) ${\frak K}$ has the 
$(\kappa,\theta)$-order property \ub{if} for every $\alpha$ there are 
$A \subseteq {\frak C}$ and $\langle \bar a_i:i < \alpha \rangle$, where 
$\bar a_i \in {}^\kappa{\frak C}$ and $|A| \le \theta$ such that:
\medskip
\roster
\item "{$(*)$}"  if $i_0 < j_0 < \alpha,i_1 < j_1 < \alpha$ then for no
$f \in AUT({\frak C})$ do we have \newline
$f \restriction A = \text{ id}_A,
f(\bar a_{i_0} \char 94 \bar a_{j_0}) = \bar a_{j_1} \char 94 \bar a_{i_1}$.
\endroster
\medskip

\noindent
If $A = \emptyset$ i.e. $\theta = 0$, we write ``$\kappa$-order property". \nl
2) ${\frak K}$ has the $(\kappa_1,\kappa_2,\theta)$ order property \ub{if}
for every $\alpha$ there are $A \subseteq {\frak C}$ satisfying $|A| \le 
\theta,\langle \bar a_i:i < \alpha \rangle$ where $\bar a_i \in
{}^{\kappa_1}{\frak C}$ and $\langle \bar b_i:i < \alpha \rangle$ 
where $\bar b_i \in {}^{\kappa_2}{\frak C}$ such that
\mr
\item "{$(*)$}"  if $i_0 < j_0 < \alpha,i_1 < j_1 < \alpha$, \ub{then} for
no $f \in \text{ AUT}({\frak C})$ do we have \nl
$f \restriction A =
\text{ id}_A,f(\bar a_{i_0}) = \bar a_{j_1},f(\bar b_{j_0}) = \bar b_{i_1}$.
\endroster
\enddefinition
\bigskip

\demo{\stag{4.1C} Observation}  So we have obvious monotonicity properties
and if $\theta
\le \kappa$ we can let $A = \emptyset$; so the $(\kappa,\theta)$-order 
property implies the $(\kappa + \theta)$-order property.
\enddemo
\bigskip

\proclaim{\stag{4.1D} Claim}  1) Any strictly indiscernible sequence (over $A$)
is an indiscernible sequence (over $A$). \nl
2) Any indiscernible set (over $A$) is an indiscernible set (over $A$).
\endproclaim
\bigskip

\proclaim{\stag{4.2} Claim}  1) If $\mu \ge LS({\frak K}) + |Y|$ and 
$h^\theta_i:Y \rightarrow {\frak C}$, for $i < \theta < \beth_{(2^\mu)^+}$ 
(e.g. $h^\theta_i = h_i$) \ub{then} we can find $\langle h'_j:j < i^* 
\rangle$, a strictly indiscernible sequence, with $h'_j:Y \rightarrow 
{\frak C}$ such that: 
\medskip
\roster
\item "{$(*)$}"  for every $n < \omega,j_1 < \dots < j_n < i^*$ for
arbitrarily large $\theta < \beth_{(2^\mu)^+}$ we can find 
$i_1 < \dots < i_n < \theta$ and $f \in AUT({\frak C})$ such that  
$h'_{j_\ell} \circ (h^\theta_{i_\ell})^{-1} \subseteq f$.
\endroster
\medskip

\noindent
2)  If in part (1) for each $\theta$, the sequence 
$\langle h^\theta_j:j < \theta \rangle$ is an indiscernible sequence of
character $\aleph_0$, in $(*)$ any  
$i_1 < \dots < i_n < i^*$ will do. \newline
3) In Definition \scite{4.1B} we can restrict $\alpha$ to
$\alpha < \beth_{(2^{\kappa + \theta + LS({\frak K})})^+}$ and get an 
equivalent version. \newline
4) In Definition \scite{4.1B} we can demand $\langle \bar a \char 94
\bar a_i:i < \alpha \rangle$ is strictly indiscernible (where $\bar a$ lists
$A$) and get an equivalent version. \nl
5) If $\mu \ge LS({\frak K}) + |Y|,N \le_{\frak K} {\frak C}$ and
$h^\theta_i:Y \rightarrow N$ for $i < \theta < \beth_{(2^\mu)^+}$ and $N^1$
is an expansion of $N$ with $|\tau(N^1)| \le \mu$, \ub{then} for some
expansion $N^2$ of $N^1$ with $|\tau(N^2)| \le \mu$ and $\Psi$ we have:
\mr
\item "{$(a)$}"  $\tau(\Psi) = \tau(N^2)$ 
\sn
\item "{$(b)$}"  for linear orders $I \subseteq J$ we have \nl
$EM_{\tau({\frak K})}(I,\Psi) \le_{\frak K} EM_{\tau({\frak K})}(I,\Psi) 
\in K$ \nl
and the skeleton of $EM_{\tau({\frak K})}(I,\Psi)$ is 
$\langle \bar a_t:t \in I \rangle,\bar a_t = \langle a_{t,y}:y \in Y \rangle$
\sn
\item "{$(c)$}"  for every $n < \omega$ for arbitrarily large $\theta <
\beth_{(2^\mu)^+}$ for some $i_0 < \dots i_{n-1} < \theta$, for every linear
order $I$ and $t_0 < \dots < t_{n-1}$ in $I$, letting $J = \{t_0,\dotsc,
t_{n-1}\}$ there is an isomorphism $g$ from $EM(J,\Psi) \subseteq
EM(I,\Psi)$ (those are $\tau(N^2)$-models) onto the submodel of
$N^2$ generated by $\dbcu_{\ell < n} \text{Rang}(h^\theta_{i_\ell})$ such
that \nl
$h^\theta_{i_\ell}(y) = g(a_{t,y})$.
\endroster
\endproclaim
\bigskip

\demo{Proof}  As in \cite[Ch.VII,\S5]{Sh:c} and \cite{Sh:88} [Saharon read],
see \scite{X1.3B} for a similar somewhat more complicated proof.
\enddemo
\bigskip

\proclaim{\stag{4.3} Lemma}  1) If there is a strictly indiscernible sequence
which is not an indiscernible set of character $\aleph_0$ called 
$\langle \bar a^i:i < \omega \rangle$, \ub{then} ${\frak K}$ has the 
$|\ell g(a^i)|$-order property. \newline
2) If there is $\langle \bar a^i:i < i^* \rangle$ is a strictly indiscernible
sequence over $A$ of character $\theta^+$ but is not an indiscernible set over
$A$ of character $\theta^+$ and $i^* \ge \theta^+$, \ub{then} ${\frak K}$ has
the $(\ell g(\bar a^0),|A| + \theta \times \ell g(\bar a^0)$-order property.
\endproclaim
\bigskip

\remark{Remark}  Permutation of infinite sets is a more complicated issue.
\endremark
\bigskip

\proclaim{\stag{4.3A} Claim}  1) If ${\frak K}$ has the $\kappa$-order 
property \ub{then}:

$$
I(\chi,\kappa) = 2^\chi \text{ for every } \chi > (\kappa + LS({\frak K}))^+
$$
\medskip

\noindent
(and other strong non-structure properties). \nl
2)  If ${\frak K}$ has the $(\kappa_1,\kappa_2,\theta)$-order property and
$\chi \ge \kappa = \kappa_1 + \kappa_2 + \theta$ \ub{then} for some
$M \in K_\chi$, we have $|{\Cal S}^{\kappa_2}(M)/E_\kappa| > \chi$. 
\endproclaim
\bigskip

\demo{Proof}  1) By \cite[Ch.III,\S3]{Sh:e} (preliminary version appears in
\cite[Ch.III,\S3]{Sh:300}) (note the version on e.g.
$\triangle(L_{\lambda^+,\omega})$). \newline
2) Straight.  \hfill$\square_{\scite{4.3A}}$
\enddemo
\bigskip

\definition{\stag{4.4} Definition}  1) Suppose $M \le_{\frak K} N$ and 
$p \in {\Cal S}^m(N)$.  Then $p$ divides over $M$ if there are elementary 
maps $\langle h_i:i < \bar \kappa \rangle$, $\text{Dom}(h_i) = N$,  
$h_i \restriction M = \text{ id}_M$,
$\langle h_i:i < \bar \kappa \rangle$ is a strictly indiscernible sequence and
$\{h_i(p):i < \bar \kappa \}$ is contradictory i.e. no element (in some
${\frak C}',{\frak C} <_{\frak K} {\frak C}'$)  realizing all of them; recall
$\bar \kappa$ is the cardinality of ${\frak C}$.  Let $\mu$-divides mean no
elements realize $\ge \mu$ of them.  \newline
2)  $\bold \kappa_\mu({\frak K})$ [or $\kappa^*_\mu({\frak K})$] is the set 
of regular $\kappa$ such that for some $\le_{\frak K}$-increasing continuous
$\langle M_i:i \le \kappa + 1 \rangle$ in $K_\mu$ and 
$b \in M_{\kappa +1}$ for every $i < \kappa$ we have:  
tp$(b,M_\kappa,M_{\kappa +1})$ [or tp$(b,M_{i+1},M_{\kappa +1})$] divides 
over $M_i$; so $\kappa \le \mu$. \newline
3) $\kappa_{\mu,\theta}({\frak K})$ [or $\kappa^*_{\mu,\theta}({\frak K})$] 
is the set of regular $\kappa$ such that
for some $\le_{\frak K}$-increasing continuous sequence 
$\langle M_i:i \le \kappa +1 \rangle$ in $K_\theta$ and 
$b \in M_{\kappa +1}$ for every $i < \kappa$ we have: 
tp$(b,M_\kappa,M_{\kappa +1})$ [or tp$(b,M_{i+1},M_{\kappa +1})$],
$\mu$-divides over $M_i$, so $\kappa \le \theta$ (see Definition \scite{4.6} 
below).
\enddefinition
\bigskip

\remark{\stag{4.4A} Remark}  1) A natural question: is there a 
parallel to forking? \nl
2)  Note the difference between $\kappa_\mu({\frak K})$ and
$\kappa^*_\mu({\frak K})$.  Note that now the ``local character" is apparently
lost.  
\endremark
\bigskip

\demo{\stag{4.5} Fact}  1) In Definition \scite{4.4}(1) we can equivalently 
demand: no element realizing $\ge \beth_{{(2^\chi)}^+}$ of them, where
$\chi  = \|N\|$. \nl
2) If $\kappa \in \kappa^*_\mu({\frak K}),\theta = \text{ cf}(\theta) \le
\kappa$ \ub{then} $\theta \in \kappa^*_\mu({\frak K})$ and similarly of
$\kappa^*_{\mu,\theta}({\frak K})$. \nl
3) $\kappa^*_\mu({\frak K}) \subseteq \kappa_\mu({\frak K})$ similarly
$\kappa^*_{\mu,\theta}({\frak K}) \subseteq \kappa_{\mu,\theta}({\frak K})$.
\enddemo
\bigskip

\definition{\stag{4.6} Definition}  Suppose $M \le_{\frak K} N,p \in 
{\Cal S}(N),M \in K_{\le \mu},\mu \ge LS({\frak K})$. \newline
1) We say $p$ does $\mu$-strongly splits over $M$, \ub{if} there are   
$\langle \bar a^i:i < \omega \rangle$  such that:
\roster
\item "{(i)}"  $\bar a^i \in {}^{\gamma \ge}{\frak C}$  for $i < \omega$,
$\gamma < \mu^+$, $\langle \bar a^i:i < \omega \rangle$ is strictly
indiscernible over $M$
\smallskip
\noindent
\item "{(ii)}"  for no $b$ realizing $p$ do we have tp$(\bar a^0 \char 94
 \langle b \rangle,M,{\frak C}) = \text{ tp}(\bar a^1 \char 94 \langle b
\rangle,M,{\frak C})$.
\endroster
\medskip
 
\noindent
2)  We say $p$ explicitly $\mu $-strongly splits over $M$ if in addition  
$\bar a^0 \cup \bar a^1 \subseteq  N$. \newline
3) Omitting $\mu$ means any $\mu$ (equivalently $\mu = \|N\|$).
\enddefinition
\bigskip

\proclaim{\stag{4.7} Claim}  1) Strongly splitting implies
dividing with models of cardinality $\le \mu$ if $(*)_\mu$ holds where  
$(*)_\mu = (*)_{\mu,\aleph_0,\aleph_0}$ and 
\medskip
\roster
\item "{$(*)_{\mu,\theta,\sigma}$}"  If  
$\langle \bar a^i:i < i^* \rangle$ is a strictly indiscernible sequence,  
$\bar a^i \in {}^\mu{\frak C},\bar b \in {}^{\sigma >}{\frak C}$, \ub{then} 
for some  $u \subseteq i^*$,  $|u| < \theta $  and the isomorphism type of
$({\frak C},\bar a^i \char94 \bar b)$ for all $i \in i^* \backslash u$
is the same.
\endroster
\endproclaim
\bigskip

\proclaim{\stag{4.8} Claim}  1) Let $\mu(*) = \mu + \sigma + LS({\frak K})$.
Assume $\langle \bar a^i:i < i^* \rangle$ and
$\bar b$ form a counterexample to $(*)_{\mu,\theta,\sigma}$ of
\scite{4.7} and $\theta \ge \beth_{(2^{\mu(*)})^+}$
\underbar{then} ${\frak K}$ has the $\mu(*)$-order property.
\newline
2) We can also conclude that for $\chi \ge \mu + LS({\frak K})$, for some
$M \in K_\chi$ we have $|{\Cal S}^{\ell g(\bar b)}(M)| > \chi$. \nl
3) If we have ``$\theta < \beth_{(2^{\mu(*)})^+}$" we can still get that
for every $\chi \ge \mu + \sigma + LS({\frak K}) + \theta$ for some 
$M \in K_\chi$, we have $|{\Cal S}^{\ell g(\bar b)}(M)| \ge \chi^\theta$. \nl
4) In part (1) it suffices to have such an example for every $\theta <
\beth_{(2^{\mu(*)})^+}$, of course, for fixed $\mu(*)$. 
\endproclaim
\bigskip

\demo{Proof}  Straight, using \scite{4.9} below.
\enddemo
\bigskip

\proclaim{\stag{4.9} Claim}  Assume 
$M = EM(I,\Phi),LS({\frak K}) + \ell g(\bar a_i) \le \mu,\mu \ge
|\alpha| + LS({\frak K})$ and $M \le_{\frak K} N,\bar b \in {}^\alpha N$ and
\medskip
\roster
\item "{$(*)$}"  for no $J \subseteq I,|J| < \beth_{(2^\mu)^+}$ do we have
for all $t,s \in I \backslash J$, \newline
tp$(\bar a_t \char 94 \bar b,\emptyset,N) = \text{ tp}
(\bar a_s \char 94 \bar b,\emptyset,N)$.
\endroster
\medskip

\noindent
\underbar{Then}
\medskip
\roster
\item "{$(A)$}"  we can find $\Phi'$ proper for linear orders and a 
formula $\varphi$ (not necessarily
first order, but $\pm \varphi$ is preserved by $\le_{\frak K}$-embeddings)
such that for any linear order $I'$ \newline
$M = EM(I',\Phi'),\bar a_t = \bar a^t \char 94 \bar b_t,\ell g(\bar a^t) \le
\mu,\ell g(\bar b_t) = \alpha$ and \newline
$M \models \varphi[\bar a^t,\bar b_s] \Leftrightarrow t < s$ \newline
(if $\alpha < \omega$, this is half the finitary order property)
\smallskip
\noindent
\item "{$(B)$}"  this implies instability in every $\mu' \ge \mu$ if
$\alpha < \omega$
\smallskip
\noindent
\item "{$(C)$}"  this implies the $(\mu + |\alpha|)$-order property and even
the $(\mu,|\alpha|,0)$-order property
\sn
\item "{$(D)$}"  if $\bar b \in {}^\alpha M$ then $``|J| < \mu^+"$ or just
$``|J| < |\alpha|^+ + \aleph_0"$ in $(*)$ suffices
\sn
\item "{$(E)$}"  if $\chi \ge \mu$, for some $M \in K_\chi$, \ub{then}
$|{\Cal S}^\alpha (M)| > \chi$ moreover $|{\Cal S}^\alpha(M)/E_\mu| > \chi$.
\endroster
\endproclaim
\bigskip

\demo{Proof}  As we can increase $I$, \wilog \,the linear order $I$ is dense
with no first or last element and is $(\beth_{(2^\mu)^+})^+$-strongly
saturated, see Definition \scite{4.11} below.  
So for some $p$ and some interval $I_0$ of $I$, the set
$Y_0 = \{t \in I_0:\text{tp}(\bar a_t \char 94 \bar b,\emptyset,N) = p\}$
is a dense subset of $I_0$.  Also for some $q \in {\Cal S}^\alpha(M)
\backslash \{p\}$, the set $Y_1 = \{t \in I:
\text{tp}(\bar a_t \char 94 \bar b,\emptyset,N) = q\}$ has cardinality 
$\ge \beth_{(2^\mu)^+}$ and let $Y'_1 \subseteq Y_1$ have cardinality
$\beth_{(2^\mu)^+}$.  As we can shrink $I_0$ \wilog \,$I_0$ is disjoint from
$Y'_1$ and as we can shrink $Y_1$ \wilog \,$(\forall s \in Y'_1)(\forall t \in
I_0)(s <^I t)$ or $(\forall s \in Y'_q)(\forall t \in I_0)(t <^I s)$. \nl
By the Erd\"os-Rado theorem, for every $\theta < \beth_{(2^\mu)^+}$ there are
$s^\theta_\alpha \in Y'_1$ for $\alpha < \theta$ such that $\langle
s^\theta_\alpha:\alpha < \theta \rangle$ is strictly increasing or strictly
decreasing; without loss of generality the case does not depend on $\theta$, 
so as we can invert $I$ without loss of generality it is increasing.  
Let $t^*_\alpha \in Y'_1$ for $\alpha < \beth_{(2^\mu)^+}$ be strictly 
increasing.  Hence (try $(p_1,p_2) = (p,q)$ and
$(p_1,p_2) = (q,p)$, one will work)
\mr
\item "{$(*)$}"  we can find $p_1 \ne p_2$ such that
\sn
\item "{$(**)$}"   for every $\theta < \beth_{(2^\mu)^+}$ there is an 
increasing sequence $\langle t^\theta_\alpha:\alpha < \theta + \theta 
\rangle$ of members of $I$ such that

$$
\text{(i)} \qquad \alpha < \theta = 
\text{ tp}(\bar a_{t^\theta_\alpha} \char 94 \bar b,\emptyset,N) = p_0
$$

$$
\text{(ii)} \qquad \theta \le \alpha < \theta + \theta \Rightarrow \text{ tp}
(\bar a_{t^\theta_\alpha} \char 94 \bar b,\emptyset,N) = p_1.
$$
[Note that we could have replaced ``increasing" by

$$
\text{(iii)} \qquad \alpha < \beta < \theta \Rightarrow t^\theta_\alpha <_I
t^\theta_\beta <_I t^\theta_{\theta + \alpha} <_I t^\theta_{\theta + \beta}.
$$
\ermn
Why?  Let $I_1 = \{t \in I:(\forall \alpha < \beth_{(2^\mu)^+}) \, t^*_\alpha
< t\}$, so every $A \subseteq I_1$ of cardinality $\le \beth_{(2^\mu)^+}$ has
a bound from below, so for some $q_1 \in {\Cal S}^\alpha(M)$ the set
$I_2 = \{t \in I_1:\text{tp}(\bar a_t \char 94 \bar b,\emptyset,N)=q_1\}$ is
unbounded from below in $I_1$.  If $q_1 \ne p$ then $q_1,p$ can serve as
$p_1,p_2$, so assume $q_1=p$, so $q,q_1$ can serve as $p_1,p_2$.]
\mn
Now we apply \scite{4.2}(5) with $h^\theta_i$ listing
$\bar a^\theta_\alpha \char 94 \bar a^\theta_{\theta + \alpha} \char 94
\bar b$ and letting $N^1$ be $EM(I,\Phi)$ (so $\tau(N^1) = \tau(\Phi))$ and 
we get $\Psi$ as there.  Now for any linear order $I^*$, look at 
$EM(I^*,\Psi)$ and its skeleton $\langle \bar a^*_t:t \in I^* \rangle$.  
Clearly $\bar a^*_t = \bar a^1_t \char 94 \bar a^2_t \char 94 \bar b^*$, 
and letting $M^*$ be the submodel of $EM_{\tau(\Phi)}(I^*,\Psi)$ generated 
by $\{a^1_t,a^2_t:t \in I^*\} \cup \bar b$, it is isomorphic to 
$EM(I^* + I^*,\Psi)$, so \wilog \, $M = M^* \restriction \tau({\frak K}) 
\le_{\frak K} {\frak C}$, so
$\text{tp}(\bar a^1_t \char 94 \bar b,\emptyset,M) = p_1,
\text{tp}(\bar a^2_t \char 94 \bar b,\emptyset,M) = p_2$.  Now for any $\chi$
we can choose $I^* = I^*_\chi$ such that $\bold D = \{J:J \text{ an initial
segment of } I^* \text{ and } J \cong I^* \text{ and } I^* \backslash J
\text{ is isomorphic to } I^*\}$ has cardinality $> \chi$.

So we have proved clause (E) and clause (B), by easy manipulations we get
clause (A) and so (C).

We are left with clause (D).  Clearly there is $\bar t = \langle t_i:i <
i^* \rangle$ satisfying $i^* < |\alpha|^+ + \aleph_0$ such that $\bar b = 
\langle b_\beta:\beta < \alpha \rangle,b_\beta = 
\tau_\beta(\bar a_{t_{i(\beta,0)}},\dotsc,
\bar a_{t_{i(\beta,n(\beta)-1)}})$ where $i(\beta,\ell) < i^*,\tau_\beta$ a
$\tau(\Phi)$-term. \nl
Let $J = \{t_i:i < i^*\}$ so by the version of $(*)$ used in clause (D),
necessarily for some $s_1,s_2 \in I \backslash J$ we have:
  
$$
p_1 \ne p_2 \text{ where}
$$

$$
p_1 = \text{ tp}(\bar a_{s_1} \char 94 \bar b,\emptyset,N)
$$

$$
p_2 = \text{ tp}(\bar a_{s_2} \char 94 \bar b,\emptyset,N)
$$
\mn
Clearly $s_1 \ne s_2$.  By renaming \wilog \,$s_1 <^I s_2$ and $0 = i_0 \le
i_1 \le i_2 \le i_3 = i^*$ and $t_i <^I s_1 \Leftrightarrow i < i_1$ and
$s_1 <^I t_i <^I s_2 \Leftrightarrow i_1 \le i < i_2$ and $s_2 <^I t_i
\Leftrightarrow i_2 < i < i_3$. \nl
As $I$ is $(\beth_{(2^\mu)^+})^+$-strongly saturated we can increase $J$ so
renaming without loss of generality 
$i(\beta,\ell) \notin \{i_1,i_2\}$, and replace
$t_{i_1},t_{i_2}$ by $s_1,s_2$.  So for every linear order $I'$ we can define
a linear order $I^*$ with a set of elements

$$
\{t_i:i < i_1 \text{ or } i_2 < i < i^*\} \cup \{(s,i):s \in I',i_1 \le i <
i_2\}
$$
\mn
linearly ordered by:

$$
t_{j_1} < t_{j_2} \text{ \ub{if} } j_1 < j_2 < i_1
$$

$$
t_{j_1} < t_{j_2} \text{ \ub{if} } i_2 < j_1 < j_2 < i^*
$$

$$
\align
t_{j_1} < (s',j') < (s'',j'') < t_{j_2} \text{ \ub{if} } &j_1 < i_1,
i_2 < j_2 < i^*, \\
  &s',s'' \in I',j',j'' \in [i_1,i_2] \\
  &(s' <^{I'} s'') \vee (s' = s'' \and j' < j'').
\endalign
$$
\mn
In $M = EM(I^*,\Phi)$ define, for $s \in I'$

$$
\bar c_{s,i} \text{ is } \bar a_{t_i} \text{ if } i < i_1 \vee i > i_2,
$$

$$
\bar c_{s,i} = \bar a_{(s,i)} \text{ if } i \in [i_1,i_2]
$$

$$
\bar b_s = \langle \tau_\beta(\bar c_{s,i(\beta,0)},\bar c_{s,i(\beta,1)},
\dotsc,\bar c_{s,i(\beta,n(\beta)-1)}):\beta < \alpha \rangle.
$$
\mn
Easily

$$
s' <^{I'} s'' \Rightarrow \text{ tp}(\bar a_{(s',i_1)} \char 94 
\bar b_{s''},\emptyset,M) = p_1
$$

$$
s'' \le^{I'} s' \Rightarrow \text{ tp}(\bar a_{(s',i_1)} \char 94 
\bar b_{s''},\emptyset,M) = p_2.
$$
\mn
By easy manipulations we can finish. \hfill$\square_{\scite{4.9}}$
\enddemo
\bigskip

\proclaim{\stag{4.10} Claim}  Assume $K$ is categorical in $\lambda$ and
\medskip
\roster
\item "{$(a)$}"  $1 \le \kappa$ and $LS({\frak K}) < \theta =
\text{ cf}(\theta) \le \lambda$ and \newline
$(\forall \alpha < \theta)(|\alpha|^\kappa < \theta)$
\smallskip
\noindent
\item "{$(b)$}"  $\bar a_i \in {}^\kappa{\frak C}$ for $i < \theta$.
\endroster
\medskip

\noindent
\ub{Then} for some $W \subseteq \theta$ of cardinality $\theta$, the
sequence $\langle \bar a_i:i \in W \rangle$ is strictly indiscernible.
\endproclaim
\bigskip

\demo{Proof of \scite{4.10}}  Let $M' \prec {\frak C},\|M'\| = \theta$ and
$\alpha < \theta \Rightarrow \bar a_\alpha \subseteq M'$.  
There is $M'',M' \prec M'' \prec {\frak C},\|M''\| = \lambda$.  So 
$M'' \cong EM(\lambda,\Phi)$ and without loss of
generality equality holds.  So there is $u \subseteq \lambda,|u| \le \theta$
such that $M' \subseteq EM(u,\Phi)$.  So without loss of generality $M' = 
EM(u,\Phi)$.
So $a_\alpha \in EM(v_\alpha,\Phi)$ for some 
$v_\alpha \subseteq u,|v_\alpha| \le \kappa$. \nl
Without loss of generality: otp$(v_\alpha) = j^*$, so for $\alpha < \beta$,
OP$_{u_\alpha,u_\beta}$ the order preserving map from $v_\beta$ onto
$v_\alpha$ induces $f_{\alpha,\beta}:EM(u_\beta,\Phi) 
\overset \text{iso}\to{\underset \text{onto}\to \longrightarrow}
EM(u_\alpha,\Phi)$, and without loss of generality 
$f_{\alpha,\beta}(\bar a_\beta) = \bar a_\alpha$.

Now as $u$ is well ordered and the assumption (a), (or see below) for some 
$w \in [\theta]^\theta$ the sequence $\langle v_\alpha:\alpha \in w \rangle$
\ub{is indiscernible} in the linear order sense (make them a sequence).  Now
we can create the right $\Phi$. 
\mn
[Why?  Let $u_\alpha = \{\gamma_{\alpha,j}:j < j^*\}$ where
$\gamma_{\alpha,j}$ increases with $j$.  For $\alpha < \theta$, let \nl
$A_\alpha = \{\gamma_{\beta,j}:
\beta < \alpha,j < j^*\} \cup \{\dbcu_{\beta < \alpha,j} 
\gamma_{\beta,j} + 1\}$.
Let $\gamma^*_{\beta,j} = \text{ Min}\{\gamma \in A_\alpha:\gamma_{\beta,j}
\ge \gamma\}$ and for each $\alpha \in S^*_0 = \{\delta < \theta:
\text{cf}(\delta) > \kappa\}$ let $h(\delta) = \text{ Min}\{\beta < \delta:
\gamma^*_{\delta,j} \in A_\beta\}$ (defining $\langle A_\beta:\beta \le
\delta \rangle$ as increasing continuous, cf$(\delta) > \kappa \ge |j^*|$
and $\gamma^*_{\delta,j} \in A_\delta$ by definition).
\mn
By Fodor's lemma for some stationary $S_1 \subseteq S_0,h \restriction S_1$
is constantly $\beta^*$.  As \nl
$(\forall \alpha < \theta)(|\alpha|^\kappa <
\theta = \text{ cf}(\theta))$ for some $S_2 \subseteq S_1$ for each
$j < j^*$ and for all $\delta \in S_2$, the truth value of
$``\gamma_{\delta,j} \in A_\delta"$ (e.g. $\gamma_{\delta,j} = \gamma^*
_{\delta,j}$) is the same and $\langle \gamma^*_{\delta,j}:\delta \in S_2
\rangle$ is constant.  Now $\langle u_\delta:\delta \in S_2 \rangle$ is as
required.  See more \cite[\S7]{Sh:620}.] \hfill$\square_{\scite{4.10}}$
\enddemo
\bigskip

\definition{\stag{4.11} Definition}  A model $M$ is $\lambda$-strongly 
saturated if:
\medskip
\roster
\item "{$(a)$}"  $\lambda$-saturated
\smallskip
\noindent
\item "{$(b)$}"  strongly $\lambda$-homogeneous: if $f$ is a partial
elementary mapping from $M$ to $M$, $|\text{Dom}(f)| < \lambda$ \newline
then $(\exists g \in AUT(M))(f \subseteq g)$.
\endroster
\medskip

\noindent
Note: if $\mu = \mu^{< \lambda},I$ a linear order of cardinality $\le \mu$,
then there is a $\lambda$-strongly saturated dense linear order
$J,I \subseteq J$.
\enddefinition
\bigskip

\remark{Remark}  We can even get a uniform bound on $|J|$ (which only depends
on $\mu$).
\endremark
\newpage

\head {\S5 Rank and Superstability} \endhead  \resetall
\bigskip

\definition{\stag{5.1} Definition}  For $M \in K_\mu,p \in {\Cal S}^m(M)$
we define $R(p)$ an ordinal or $\infty$ as follows: 
$R(p) \ge \alpha$ iff for every  
$\beta < \alpha$ there are $M^+,M \le_{\frak K} M^+ \in K_\mu,
p \subseteq p^+ \in {\Cal S}^1(M^+),R(p^+) \ge \beta \and$
[$p^+$\,$\mu$-strongly splits over $M$].  In case of doubt we write $R_\mu$.
This is well defined and has the obvious properties: 
\mr
\item "{$(a)$}"  monotonicity,
\sn
\item "{$(b)$}"  if $M \in K_\mu,p \in {\Cal S}^m(M)$ and Rk$(p) \ge \alpha$
\ub{then} for some $N,q$ satisfying $M \le_{\frak K} N \in K_\mu$ and
$q \in {\Cal S}^m(N)$ we have: $q \restriction M = p$ and Rk$(q) = \alpha$
\sn
\item "{$(c)$}"  automorphisms of ${\frak C}$ preserve everything
\sn
\item "{$(d)$}"  the set of values is 
$[0,\alpha)$ or $[0,\alpha) \cup \{\infty\}$ for some $\alpha < 
(2^\mu)^+$, etc.
\endroster
\enddefinition
\bigskip

\definition{\stag{5.2} Definition}  We say ${\frak K}$ is $(\mu,1)$-superstable
if \newline
$M \in K_\mu \and p \in {\Cal S}(M) \Rightarrow R(p) < \infty$ \, 
$\biggl( \text{equivalently } < (2^\mu)^+ \biggr)$.
\enddefinition
\bigskip

\proclaim{\stag{5.3} Claim}  If $(*)_\mu$ from \scite{4.7} above fails, 
\underbar{then} ($\mu,1$)-superstability fails.
\endproclaim
\bigskip

\demo{Proof}  Straight.
\enddemo
\bigskip

\proclaim{\stag{5.4} Claim}  If ${\frak K}$ is not 
($\mu,1$)-superstable, \ub{then} 
there are a sequence \newline
$\langle M_i:i \le \omega + 1 \rangle$ which is $<_{\frak K}$-increasing 
continuous in $K_\mu$ and $m < \omega$ and \nl
$\bar a \in {}^m(M_{\omega +1})$ such that 
$(\forall i < \omega) \bigl[ {\bar a \over M_{i+1}}$ does 
$\mu$-strongly split over $M_i \bigr]$. \nl
Also the inverse holds.
\endproclaim
\bigskip

\demo{Proof}  As usual.
\enddemo
\bigskip

\proclaim{\stag{5.5} Claim}  1) If ${\frak K}$ is not $(\mu,1)$-superstable
\underbar{then}  $K$ is unstable in every  $\chi$  such that
$\chi^{\aleph_0} > \chi  + \mu  + 2^{\aleph _0}$. \newline
2)  If $\kappa \in \kappa^*_\mu({\frak K})$ and $\chi^\kappa > \chi \ge
LS({\frak K})$, \ub{then} ${\frak K}$ is not $\chi$-stable, even 
modulo $E_\mu$. \nl
3) If $\kappa \in \kappa_\mu({\frak K})$ and $\chi^\kappa > \chi =
\chi^\kappa \ge LS({\frak K})$ or just there is a tree with $\chi$ nodes and
$> \chi \, \kappa$-branches and $\chi \ge LS({\frak K})$, \ub{then}
${\frak K}$ is not $\chi$ stable even modulo $E_\mu$.
\endproclaim
\bigskip

\remark{Remark}  We intend to deal with the following elsewhere; we need
stable amalgamation
\mr
\item "{$(*)$}"  if $\kappa \in \bold \kappa_\mu({\frak K})$, 
cf$(\chi) = \kappa,\dsize \bigwedge_{\lambda < \chi} \lambda^\mu \le \chi$,\nl
\underbar{then} ${\frak K}$ is not $\chi$-stable.
\endroster
\endremark
\bigskip

\remark{\stag{5.5A} Remark}  1) In (1) this implies 
$I(LS({\frak K})^{+(\omega(\alpha_0+\alpha)+n)},K) \ge |\alpha|$ when
$\mu = \aleph_{\alpha_0}$.  We conjecture that \cite{GrSh:238} can be
generalized to the content of (1) with cardinals which exists by ZFC. \nl
2)  Note that for FO stable theory $T,{\frak K} = \text{ MOD}(T)$, for 
$\kappa$ regular we have
$(*)^\kappa_1 \Leftrightarrow (*)^\kappa_2$ where 
\medskip
\roster
\item "{$(*)^\kappa_1$}"  for any increasing chain $\langle M_i:i < \kappa 
\rangle$ of $\lambda$-saturated models of length \nl
$\kappa$, the union $\dsize \bigcup_{i < \kappa} M_i$ is $\lambda$-saturated,
\smallskip
\noindent
\item "{$(*)^\kappa_2$}"  $\kappa \in \kappa_r({\frak K})$.
\endroster
\medskip

\noindent
In \cite{Sh:e}, $(*)^\kappa_2$ is changed to
\medskip
\roster
\item "{$(**)$}"  $\kappa < \bold\kappa_r(T)$ \newline
(really $\bold\kappa_r({\frak K}$) (i.e. $\kappa_r(T)$) is a \ub{set} of
regular cardinals)).
\endroster
\medskip

\noindent
\relax From this point of view, FO theory $T$ is a degenerated case: $\bold \kappa_r
(T)$ is an initial segment so naturally we write the first regular not in
it.  This is a point where \cite{Sh:300} opens our eyes. \nl
3) In fact in \scite{5.5} not only do we get $\|M\| = \chi,|{\Cal S}(M)| >
\chi$ but also $|{\Cal S}(M)/E_\mu| > \chi$. \nl
4)  Let me try to explain the proof of \scite{5.5}, of course, being
influenced by the first order case.  If the class is superstable, one of the
consequences of not having the appropriate order property is that
(see \scite{4.9}) for a strictly indiscernible sequence $\langle \bar a_t:
t \in I \rangle$ over $A$ each $\bar a_t$ of length at most $\mu$ and 
$\bar b$, singleton for simplicity, for all except few of the $\bar a_t$'s, 
the type
of $\bar a_t \cong \bar b$ realizes the same type.  Of course, we can get
better theorems generalizing the ones for first order theories: we can use
$\kappa \notin \kappa_\mu({\frak C})$ and/or demand that after adding to 
$A,\bar c$ and few of the $\bar a_t$'s the rest is strictly indiscernible 
over the new $A$, but this is not used in \scite{5.5}.  Now if ${\frak C}$ 
is $(\mu,1)$-superstable the
number of exceptions is finite, however, the inverse is not true: for some
non $(\mu,1)$-superstable class ${\frak C}$ 
still the number of exceptions in such
situations is finite.  In the proof of \scite{5.5}(1) this property is used
as a dividing line.
\endremark
\bigskip

\demo{Proof} 1)
\enddemo
\bn
$\bold{Case \,I}$   There are $M,N,p,\langle \bar a_i:
i < i^* \rangle$ as in \scite{4.7}$(*)_\mu$ and $\bar c$,
(in fact $\ell g(\bar c) = 1$) such that $\bar c$ realizes $h_i(p)$
for infinitely many $i$'s and fails to realize $h_i(p)$ for infinitely many
$i$'s.
 
Let $I$ be a $\beth(\chi + \beth_{(2^\mu)^+})^+$-strongly saturated 
dense linear order (see Definition \scite{4.11}) such that even 
if we omit  $\le \beth_{(2^\mu)^+}$ members,
it remains so.  By the strict indiscernibility we can find 
$\langle \bar a_t:t \in I \rangle,c$ as above. 

So there is $u \subseteq I,|u| < \beth_{(2^\mu)^+}$ such that
$q = \text{ tp}(\bar a_t \char 94 \bar c,\emptyset,{\frak C})$ is the same
for all $t \in I \backslash u$; without loss of generality $q = \text{ tp}
(\bar a_t \char 94 \bar c,\emptyset,{\frak C}) \Leftrightarrow t \in I
\backslash u$, so $u$ is infinite.  So we can find $i_n \in i^* \cap u$
such that $i_n < i_{n+1}$.  Let $I' = I \backslash (u \backslash
\{i_n:n < \omega\})$, so that $I'$ is still $\chi^+$-strongly saturated.  
Hence for every $J \subseteq I'$ of order type $\omega$ for some $c_J(\in
{\frak C})$ we have

$$
t \in I' \backslash J \Rightarrow \text{ tp}(\bar a_t \char 94 \bar c_J,
\emptyset,{\frak C}) = q
$$

$$
t \in J \Rightarrow \text{ tp}(\bar a_t \char 94 \bar c_J,\emptyset,{\frak C})
\ne q.
$$
\medskip

\noindent
This clearly suffices.
\bn
$\bold{Case \, II}$  Not Case I. \nl
As in \cite{Sh:3} (the finitely many finite exceptions do not matter) or
see part (2).  \nl
2) If $\chi < 2^\kappa$ the conclusion follows from \scite{3.2}(2).
Possibly decreasing $\kappa$ (allowable as $\kappa \in \kappa^*_\mu
({\frak K})$ rather than $\kappa \in \kappa_\mu({\frak K})$ is assumed) we
can find a tree ${\Cal T} \subseteq {}^{\kappa \ge}\chi$, so closed under 
initial segments such that $|{\Cal T} \cap {}^{\kappa >} \chi| \le \chi$ but 
$|{\Cal T} \cap {}^\kappa \chi| > \chi$.  (The cardinal arithmetic 
assumption is needed just for this).
Let $\langle M_i:i \le \kappa +1 \rangle,c \in M_{\kappa +1}$ exemplify
$\kappa \in \kappa^*_\mu({\frak K})$ and let ${\Cal T}' = {\Cal T} \cup 
\{ \eta \char 94 \langle 0 \rangle:\eta \in {}^\kappa\text{Ord}
\text{ and } i < \kappa \Rightarrow \eta \restriction i \in {\Cal T}\}$.  Now 
we can by induction on $i \le \kappa +1$ choose $\langle h_\eta:\eta \in 
{\Cal T}' \cap {}^i \chi \rangle$, such that:
\mr
\item "{$(a)$}"  $h_\eta$ is a 
$\le_{\frak K}$-embedding from $M_{\ell g(\eta)}$ into ${\frak C}$
\sn
\item "{$(b)$}"  $j < \ell g(\eta) \Rightarrow h_{\eta \restriction j}
\subseteq h_\eta$
\sn
\item "{$(c)$}"  if $i = j+1,\nu \in {\Cal T} \cap {}^j \chi$, \ub{then} 
$\langle h_\eta(M_i):\eta \in \text{ Suc}_T(\nu) \rangle$ 
is strictly indiscernible, and can be extended to a sequence of length
$\bar \kappa$ such that $\langle h_\eta(p \restriction M_i):\eta \in
\text{Suc}_I(\nu)\rangle$ is contradictory (i.e. as in Definition 
\scite{4.4}(1)).
\ermn
There is no problem to do this.  Let $M \le_{\frak K} {\frak C}$ be of
cardinality $\chi$ and include $\bigcup\{h_\eta(M_i):i < \kappa
\text{ and } \eta \in 
{\Cal T} \cap {}^i \chi\}$ hence it includes also $h_\eta(M_\kappa)$ if 
$\eta \in {\Cal T} \cap {}^\kappa \chi$ as 
$M_\kappa = \dbcu_{i < \kappa} M_i$. \nl
For $\eta \in {\Cal T} \cap {}^\kappa \chi$ let $c_\eta = 
h_{\eta \char 94 <0>}(c)$ and $M_\eta = h_\eta(M_i)$ when $\eta \in
{\Cal T} \cap {}^i\text{Ord}$ and $i \le \kappa +1$, so by \scite{4.9} 
clearly (by clause (C))
\mr
\item "{$(*)$}"  if $i < \kappa,\eta \in {\Cal T} \cap {}^i \chi$, and
$\eta \triangleleft \eta_1 \in {\Cal T} \cap {}^\kappa \chi$, \ub{then} \nl
$\{\rho \in \text{ Suc}_T(\eta):\text{for some } \rho_1,\rho \triangleleft
\rho_1 \in {\Cal T} \cap {}^\kappa \chi \text{ and }$ \nl

$\qquad \qquad \qquad c_{\rho_1} \text{ realizes tp}(c_{\eta_1},
h_{\eta_1 \restriction (i+1)}(m_{i+1}))\}$ \nl
has cardinality $< \beth_{(2^{\mu + LS({\frak K})})^+}$.
\ermn
Next define an equivalence relation $\bold e$ on ${\Cal T} \cap 
{}^\kappa \chi$:

$$
\eta_1 \, \bold e \, \eta_2 \text{ iff tp}(c_{\eta_1},M) =
\text{ tp}(c_{\eta_2},M).
$$
\mn
or just

$$
\eta_1 \bold e \eta_2 \text{ iff } (\forall \nu)[\nu \in {\Cal T} \Rightarrow
\text{ tp}(c_{\eta_1},M_\nu) = \text{ tp}(c_{\eta_2},M_\nu)].
$$
\mn
Now if for some $\eta \in {\Cal T} \cap {}^\kappa \chi,|\eta / \bold e| >
\beth_{(2^{\mu+LS({\frak K})})^+}$ then for some $\eta^* \in {\Cal T} \cap
{}^{\kappa >}\chi$, we have

$$
\{\nu \restriction (\ell g(\eta^* +1)):\nu \in \eta / \bold e\}
\text{ has cardinality } > \beth_{(2^{\mu+LS({\frak K})})^+}
$$
\mn
which contradicts $(*)$; so if $\chi \ge \beth_{(2^{\mu+LS({\frak K})})^+}$,
we are done.

But if for some $\eta \in {\Cal T} \cap {}^{\kappa >} \chi$ 
the set in $(*)$ has
cardinality $\ge \kappa$, then we can continue as in case I of the proof of
part (1) replacing ``infinite" by ``of cardinality $\ge \kappa$", so 
assume this never happens.  So above if $|\eta / \bold e| > 2^\kappa$, we 
get again a contradiction.  So if $|{\Cal T} \cap {}^\kappa 
\chi| > 2^\kappa$, we conclude $|{\Cal T} \cap {}^\kappa \chi / \bold e| = 
|{\Cal T} \cap {}^\kappa \chi|$, so we are done.  We are left with the case 
$\chi < 2^\kappa$, covered in the beginning (note that for $\chi < 2^\kappa$ 
the interesting notion is splitting). \nl
3) Proof similar to part (2).    \hfill$\square_{\scite{5.5}}$ 
\bigskip

\proclaim{\stag{5.7} Claim}  If $\lambda > \mu^+,\mu \ge LS({\frak K},K)$,
${\frak K}$ is categorical in $\lambda$ \ub{then} \newline
1)  $K$ is $(\mu,1)$-superstable. \newline
2)  $\kappa^*_\mu({\frak K})$ is empty.
\endproclaim
\bigskip

\demo{Proof}  1) Assume the conclusion fails.  If 
$\lambda > \mu^{+\omega}$, we can use \scite{5.5} + \scite{1.5} so
\wilog \, cf$(\lambda) > LS({\frak K})$.

By \scite{1.5} if $M \in K_\lambda$ then $M$ is cf$(\lambda)$-saturated.
On the other hand from the Definition of $(\mu,1)$-superstable we get
a non-$\mu^+$-saturated model.

Let $\chi = \beth_{(2^\lambda)^+}$.  Assume ${\frak K}$ is not $(\mu,1)$-
superstable so we can find in $K_\mu$ an increasing continuous sequence
$\langle M_i:i \le \kappa +1 \rangle$ and $c \in M_{\omega+1}$ such that
$p_{n+1} = \text{ tp}(c,M_{n+1},M_{\omega +1})$ $\mu$-strongly splits 
over $M_n$ for $n < \omega$.  For each $n < \omega$ let $\langle \bar a^n_i:
i < \omega \rangle$ be a strictly indiscernible sequence over $M_n$ 
exemplifying $p_{n+1}$ does $\mu$-strongly splits over $M_n$ (see Definition 
\scite{4.6}).  So we can
define $\bar a^n_i \in {\frak C}$ for $i \in [\omega,\chi)$ such that
$\langle \bar a^n_i:i < \chi \rangle$ is strictly indiscernible 
over $M_n$.
Let ${\Cal T}_n = \{\eta \in {}^{2n}\chi:\eta(2m) < \eta(2m+1) \text{ for }
m < n\}$.  For $n < \omega,i < j < \chi$ let $h^n_{i,j} \in \text{ AUT}
({\frak C})$ be such that $h^n_{i,j} \restriction M_n = \text{ id},
h^n_{i,j}(\bar a^n_0 \char 94 \bar a^n_1) = \bar a^n_i \char 94 \bar a^n_j$.
Now we choose by induction on $n < \omega,\langle f_\eta:\eta \in {\Cal T}_n 
\rangle,\langle g_\eta:\eta \in {\Cal T}_n \rangle,\langle a^n_i:i < \chi,
\eta \in {\Cal T}_n \rangle$ such that:
\mr
\item "{$(a)$}"  $f_\eta,g_\eta$ are restrictions of automorphisms of
${\frak C}$
\sn
\item "{$(b)$}"  Dom$(f_\eta) = M_n$
\sn
\item "{$(c)$}"  $g_\eta \in \text{ AUT}({\frak C})$
\sn
\item "{$(d)$}"  $\bar a^n_i = g_\eta(\bar a^n_i)$ if $\eta \in T_n$
\sn
\item "{$(e)$}"  $f_{<>} = \text{ id}_{M_0}$,
\sn
\item "{$(f)$}"  $f_\eta \subseteq g_\eta$
\sn
\item "{$(g)$}"  if $\eta \in {}^{2n}\chi,m < n$ then
$f_{\eta \restriction (2m)} \subseteq f_\eta$
\sn
\item "{$(h)$}"  if $\eta \in {}^{2n}\chi$ and $i < j < \chi$ then
$f_{\eta \char 94 <i,j>} = (g_\eta \circ h^n_{i,j}) \restriction M_{n+1}$.
\ermn
There is no problem to carry the induction.  Now choose by induction on
$n,M^*_n,\eta_n,i_n,j_n$ such that
\mr
\item "{$(\alpha)$}"  $i_n < j_n < \chi$ and $\eta_n = \langle i_0,j_0,
\dotsc,i_{n-1},j_{n-1} \rangle$ so $\eta _n \in {\Cal T}_n$
\sn
\item "{$(\beta)$}"  $M_n \in K_\lambda,M^*_n <^1_{\mu,\omega} M^*_{n+1}$
\sn
\item "{$(\gamma)$}"  Rang$(f_{\eta_n}) \subseteq M_n$
\sn
\item "{$(\delta)$}"  $\bar a^{\eta_n}_i,\bar a^{\eta_n}_j$ realizes the
same type over $M_n$
\sn
\item "{$(\varepsilon)$}"  $\bar a^{\eta_n}_i,\bar a^{\eta_n}_j \subseteq
M^*_{n+1}$.
\ermn
There is no problem to carry the induction (using the theorem on existence
of strictly indiscernibles to choose $i_n < j_n$).

So $\dbcu_{n < \omega} f_{\eta_n}$ can be extended to $f \in \text{ AUT}
({\frak C})$.  Let $c^* = f(c),M^*_\omega = \dbcu_n M^*_{\eta_n},
M^*_{\omega+1} = f(M_{\omega+1})$.  Clearly tp$(c,M^*_{n+1},M^*_{\omega+1})$
does $\mu$-split over $M_n$ hence $M_\omega$ is not $\mu^+$-saturated (as
cf$(\lambda) > \mu$) (see \scite{5.7a}); contradiction. \nl
2) Follows.  \hfill$\square_{\scite{5.7}}$
\enddemo
\bigskip

\proclaim{\stag{5.7a} Claim}  If $\mu \ge LS({\frak K}),\langle M_i:i \le
\delta \rangle$ is $\le_{\frak K}$-increasing continuous, \newline
$p \in {\Cal S}^{\le \mu}(M_\delta),p$\,\,$\mu$-strongly splits over 
$M_i$ for all $i$ (or just $\mu$-splits over $M_i$) and $\delta < \mu^+$ 
\underbar{then} $M_\delta$ is not $\mu^+$-saturated.
\endproclaim
\bigskip

\demo{Proof}  Straight.
\enddemo
\bigskip

\proclaim{\stag{5.7B} Claim}  Assume there is a Ramsey cardinal
$> \mu +LS({\frak K})$.  If ${\frak K}$ is not $(\mu,1)$-superstable,
\ub{then} for every $\chi > \mu + LS({\frak K})$ there are $2^\chi$ pairwise
non-isomorphic models in ${\frak K}_\chi$.
\endproclaim
\bigskip

\demo{Proof}  By \cite{GrSh:238} for $\chi$ regular; together with
\cite{Sh:e} for all $\chi$.
\enddemo
\bigskip

\proclaim{\stag{5.8} Lemma}  1) If for some $M,|{\Cal S}(M)/E_\mu| > \chi 
\ge \|M\| + \beth_{(2^\mu)^+}$ and $\mu \ge LS({\frak K})$  \underbar{then} 
${\frak K}$ is not $(\mu,1)$-superstable. \newline
2)  If $\chi^\kappa \ge |{\Cal S}(M)/E_\mu | > \chi^{<\kappa} \ge \chi \ge 
\|M\| + \beth_{(2^\mu)^+},\mu \ge LS({\frak K}) + \kappa$ \underbar{then}  
$\kappa \in \kappa^*_\mu({\frak K})$.
\endproclaim
\bigskip

\demo{Proof}  No new point when you remember the definition of 
$E_\mu$ (see \scite{1.6}).
\enddemo
\newpage

\head {\S6 Existence of Many Non-Splitting} \endhead  \resetall
\bigskip

\subhead {\stag{6.1} Question} \endsubhead   Suppose $\kappa + LS({\frak K})
\le \mu < \lambda$ and $\bar N = \langle N_i:i \le \delta \rangle$ is 
$<^1_{\mu,\kappa}$-increasing continuous (we mean for $i < j,j$ non-limit  
$N_i <^1_{\mu,\kappa} N_j$), $\delta < \mu^+$ and $p \in 
{\Cal S}^m(N_\delta)$.  Is there $\alpha < \delta$ such that for every 
$M \in {\frak K}_{\le \lambda},N_\delta \le_{\frak K} M,p$ has an extension 
$q \in {\Cal S}^m(M)$ which does not $\mu$-split over $N_\alpha$ (and so
in particular $p$ does not $\mu$-split over $N_\alpha$).  
\bigskip

\remark{\stag{6.1A} Remark}  If $p \restriction N_{\alpha +1}$ does not
$\mu$-split over $N_\alpha$, then 
$p \restriction N_{\alpha +1}$ has at most one extension  
$\text{mod } E_\mu$ which does not $\mu $-split over $N_\alpha$ because
$N_{\alpha +1} \in K_\mu$ is universal over $N_\alpha,N_{\alpha +1} 
\le_{\frak K} M \in K_\lambda$.  So in \scite{6.1} if $p$ does not 
$\mu$-split over $N_\alpha$, \ub{then} there is at most one $q/E_\mu$.
\endremark
\bigskip

\proclaim{\stag{6.2} Lemma}  Suppose $K$ is categorical in $\lambda,
\text{{\rm cf\/}}(\lambda) > \mu \ge LS({\frak K})$.  \ub{Then} the 
answer to question \scite{6.1} is yes.
\endproclaim
\bigskip

\remark{\stag{6.2A} Remark}  We intend later to deal with the case
$\lambda > \mu \ge \text{ cf}(\lambda) + LS({\frak K})$ as in 
\cite{KlSh:362}.
\endremark
\bigskip

\subhead {Notation} \endsubhead  $I \times \alpha$ is $I + I + \dots (\alpha$ 
times) (with the obvious meaning).
\bigskip

\demo{Proof}  Let $\Phi$ be proper for linear order, 
$|\tau(\Phi)| \le LS({\frak K}),
EM_\tau(I,\Phi) \in K$ (of power $|I| + \mu(K)$) where $I$ is a linear order, 
of course and $I \subseteq J \Rightarrow EM_\tau(I,\Phi) \le_{\frak K}
EM_\tau(J,\Phi)$.  So $EM_\tau(\lambda,\Phi)$ is $\mu^+$-saturated (by
\scite{1.5}).  Let $I^*$ be a linear order of power $\mu$ such that 
$I^* \times (\alpha + 1) \cong I^*$ for 
$\alpha < \mu^+$ and $I^* \times \omega \cong I^*$.  
By \scite{1.5} we know that $EM_\tau (I^* \times \lambda,
\Phi)$ is $\mu^+$-saturated. 

Now we choose by induction on $i$ an ordinal $\alpha_i < \mu^+$ and an 
isomorphism $h_i$ from $N_{1+i}$ onto $EM(I^* \times \alpha_i,\Phi)$,
both increasing with $i$ where $N_{i+1}$ is from \scite{6.1} and
cf$(\alpha_i) = \aleph_0$ for $i$ nonlimit. 

For $i = 0$, use the proof of the uniqueness of $N_1$ over $N_0$ 
(see \scite{2.3} and reference there); more specifically using the back 
and forth argument
we can find $J_0 \subseteq \lambda,|J_0| = \mu$ and isomorphism $h_0$ from  
$N_1 = N_{0+1}$ onto $EM(I^* \times J_0,\Phi) \subseteq (I^* \times \lambda,
\Phi)$.  Now let $J^0 = J_0 \cup \{\alpha < \lambda:(\forall \beta \in J_0)\,
\beta < \alpha\}$  so $J^0 \cong \lambda$ (note: $J_0$ is bounded in
$\lambda$ as cf$(\lambda) > \mu \ge |J_0|$) and also $EM_\tau(I^* \times J^0,
\Phi)$ is $\mu^+$-saturated (being isomorphic to $EM_\tau(I^* \times \lambda,
\Phi)$), so without loss of generality $J_0$ is some ordinal $\alpha_0 < 
\mu^+.$ 

So we have $h_0$.  The  continuation is similar. \newline
Now $h_\delta$ is defined $h_\delta:N_\delta \overset{\text{onto}} \to
\rightarrow EM_\tau(I^* \times \alpha_\delta,\Phi)$,  so as $EM_\tau
(I^* \times \lambda,\Phi)$ is $\mu^+$-saturated, $h_\delta(p)$ is realized 
say by $\bar a$, so let $\bar a = \bar \sigma(x_{(t_1,\gamma_1)},\dotsc,
x_{(t_n,\gamma_n)})$ where $\bar \sigma$ is a sequence of terms in
$\tau(\Phi)$ and $(t_\ell,\gamma_\ell)$ is increasing with $\ell$ (in 
$I^* \times \lambda$).  Let $\beta < \delta$ be such that:

$$
\bigl\{ \gamma_1,\dotsc,\gamma_n \bigr\} \cap \alpha_\delta  
\subseteq \alpha_\beta.
$$

\medskip
\noindent
Let
$$
\gamma'_\ell  = \cases \gamma_\ell &\text{ \underbar{if} } \quad 
\gamma_\ell < \alpha_\delta \\
 \lambda + \gamma_\ell &\text{ \underbar{if} } \quad 
\gamma_\ell \ge \alpha_\delta
\endcases
$$
\medskip
\noindent
Then in the model $N = EM_{\tau({\frak K})}(I^* \times \lambda + 
\lambda,\Phi)$, we shall show that the finite sequence $\bar \sigma^1 = 
\bar \sigma \bigl( x_{(t_1,\gamma'_1)},\dotsc,x_{(t_n\gamma'_n)} \bigr)$ 
realizes a type as required over $M = EM_{\tau({\frak K})}
(I^* \times \lambda,\Phi)$. \nl
Why?  Let $M_\gamma = EM_{\tau({\frak K})}(I^* \times \alpha_\gamma,\Phi)$
for $\gamma < \delta$.  Assume toward contradiction that 
\mr
\item "{$(*)$}"  tp$(\bar a',M,N)$ does $\mu$-split over $M_{\beta +1}$.
\ermn
Let $\bold {\bar c},\bold {\bar b} \in {}^\mu M$ realize the same type over 
$M_{\beta +1}$ but witness splitting.

We can find $w \subseteq \lambda,|w| \le \mu$ such that $\bold{\bar c},
\bold{\bar b} \subseteq EM(I^* \times w,\Phi)$.  Choose $\gamma$ such that

$$
\sup(w) < \gamma < \lambda.
$$
\mn
Let $M^- = EM_{\tau({\frak K})}(I^* \times 
(\alpha_\delta \cup w \cup [\gamma,\lambda)),\Phi) <_{\frak K} M$.
\sn
Let $N^- = EM_{\tau({\frak K})}(I^* \times 
(\alpha_\delta \cup w \cup [\gamma,\lambda) \cup
[\lambda,\lambda + \lambda)),\Phi) <_{\frak K} N$.
\sn
So still $\bold{\bar c},\bold{\bar b}$ witness that tp$(\bar a',M^-,N^-)$
does $\mu$-split over $M_{\beta +1}$.
\mn
There is an automorphism $f$ of the linear order
$I^* \times (\alpha_\delta \cup w \cup [\gamma,\lambda)) \cup [\lambda,
\lambda + \lambda))$ such that

$$
f \restriction (I^* \times \alpha_{\beta +1}) = \text{ the identity}
$$

$$
f \restriction (I^* \times [\gamma + 1,\lambda + \lambda)) =
\text{ the identity}
$$

$$
\text{Rang}\bigl(f \restriction (I^* \times w)\bigr) 
\subseteq I^* \times [\alpha_{\beta +1},\alpha_{\beta +2}).
$$
\mn
Now $f$ induces an automorphism of $N^-$ naturally called $\hat f$. \nl
So
 
$$
\hat f \restriction M_\beta = \text{ identity}
$$

$$
\hat f(\bar a') = \bar a'
$$

$$
\hat f(M^-) = M^-
$$
\mn
As $\hat f$ is an automorphism, $\hat f(\bold{\bar c}),\hat f(\bold{\bar b})$
witness that tp$(\hat f(\bar a'),\hat f(M^-),\hat f(N^-))$ does $\mu$-splits 
over $\hat f(M_{\alpha_{\beta +1}})$; i.e. tp$(\bar a',M^-,N^-)$ does 
$\mu$-splits over $M_{\alpha_{\beta +1}}$.  
So tp$(\bar a',M_{\alpha_{\beta+2}},N)$ does $\mu$-splits over 
$M_{\alpha_{\beta +1}}$.
\sn
Now choose $\alpha_\gamma < \mu^+$ for $\gamma \in (\delta,\mu^+]$, increasing
continuous by

$$
\alpha_{\delta +i} = \alpha_\delta + i
$$

$$
M_\gamma = EM_{\tau({\frak K})}(I^* \times \alpha_\gamma,\Phi).
$$
\mn
So $\langle M_\gamma:\gamma \le \mu \rangle$ is 
increasing continuous.
So for $\gamma_1 \in [\beta,\mu^+)$ there is 
$f \in AUT(I^* \times (\lambda + \lambda))$ such that

$$
f \restriction I^* \times \alpha_\beta = \text{ identity}
$$

$$
f \text{ takes } I^* \times [\alpha_\beta,\alpha_{\beta +1}) \text{ onto } 
I^* \times [\alpha_\beta,\alpha_{\gamma_1 + 1})
$$

$$
f \text{ takes } I^* \times [\alpha_{\beta +1},\alpha_{\beta +2}) 
\text{ onto } I^* \times \{\alpha_{\gamma_1 + 1}\}
$$

$$
f \text{ takes } I^* \times [\alpha_{\beta + 2},\alpha_{\gamma_1 + 2}) 
\text{ onto } I^* \times \{\alpha_{\gamma_1 + 2}\}
$$

$$
f \restriction I^* \times [\alpha_{\gamma_1 +2},\lambda + \lambda) = 
\text{ identity}.
$$
\mn
As before this shows (using obvious monotonicity of $\mu$-splitting)

$$
\text{tp}(\bar a^1,M_{\gamma_1 + 2}N)\,\mu\text{-splits over }
M_{\gamma_1 + 1}.
$$
\mn
So $\{\gamma < \mu:\text{tp}(\bar a',M_{\gamma +1},N)$ does $\mu$-split over
$M_\gamma\}$ has order type $\mu$, so without loss of generality is $\mu$.
By \scite{3.2}(2) we get a contradiction.
\hfill$\square_{\scite{6.2}}$
\enddemo
\bigskip

\proclaim{\stag{6.3} Theorem}  Suppose $K$ categorical in $\lambda$ and the
model in $K_\lambda$ is $\mu^+$-saturated (e.g. cf$(\lambda) > \mu$) and 
LS$({\frak K}) \le \mu < \lambda$. \newline
1)  $M <^1_{\mu,\kappa} N \Rightarrow N$ is saturated if
$LS({\frak K}) < \mu$. \newline
2)  If $\kappa_1,\kappa_2$ and for $\ell = 1,2$ we have
$M_\ell <^1_{\mu,\kappa_\ell} N_\ell$, \ub{then} $N_1 \cong N_2$. \newline
3) There is $M \in K_\mu$ which is saturated.
\endproclaim
\bigskip

\remark{\stag{6.3A} Remark}  1) The model we get by (2) we call
{\bf the saturated model} of ${\frak K}$ \newline
in $\mu$. \newline
2)  Formally --- we do not use \scite{6.2}. \nl
3)  By the same proof $M \le^1_{\mu,\kappa_\ell} N_\ell \Rightarrow N_1
\cong_M N_2$ and we call $N$ {\bf saturated over} $M$.
\endremark
\bigskip

\demo{Proof}  1) By the uniqueness proofs \scite{2.2} as 
$M <^1_{\mu,\kappa} N$ there are \newline
$\langle M_i:i \le \kappa \rangle,
M_i <^1_{\lambda,\kappa} M_{i+1},<_{\frak K}$-increasing continuous 
$M_0 = M,M_\kappa = N$ and as in the proof of \scite{6.2} without 
loss of generality $M_i = EM(\alpha_i,\Phi)$ where $\alpha_i < \mu^+$. 

To prove $N = N_\kappa$ is $\mu$-saturated suppose $p \in {\Cal S}^1(M^*),
M^* \le_{\frak K} N,\|M^*\| < \mu$; as we can extend $M^*$ (as long as its 
power is $< \mu$ and it is $<_{\frak K} N$), without loss of generality 
$M^* = EM(J,\Phi),J \subseteq \alpha_\kappa,|J| < \mu$. 

So for some $\gamma$ we have $[\gamma,\gamma + \omega) \cap J = \emptyset$
and $\gamma + \omega \le \alpha_\kappa$.  We can replace $[\gamma,\gamma 
+ \omega)$ by a copy of $\lambda$; this will make the model $\mu$-saturated
[alternatively, use $I^* \times$ ordinal as in a previous proof].
\newline
But easily this introduces no new types realized over $M^*$.  So $p$ is 
realized. \newline
2) Follows by part (1) or its proof.  \newline
3)  Follows from the proof of part 1).  Left to the reader.
\enddemo
\bigskip

\remark{Remark}  In part (1) we have used just 
cf$(\lambda) > \mu > LS({\frak K})$.  \hfill$\square_{\scite{6.3}}$
\endremark
\bigskip

\proclaim{\stag{6.4} Claim}  Assume $K$ categorical in $\lambda$,
cf$(\lambda) > \mu > LS({\frak K})$.  If $N_i \in K_\mu$ is saturated,
increasing with $i$ for $i < \delta$ and $\delta < \mu^+$ \ub{then}
$N = \dbcu_{i < \delta} N_i \in K_\mu$ is saturated.
\endproclaim
\bigskip

\demo{Proof}  We prove this by induction on $\delta$, so \wilog \,$\langle
N_i:i < \delta \rangle$ is not just $\le_{\frak K}$-increasing and also
contradicts the conclusion but is increasing continuous and each $N_i$
saturated.  Without loss of generality $\delta = \text{ cf}(\delta)$.
If cf$(\delta) = \mu$ the conclusion clearly holds so assume cf$(\delta) <
\mu$.  Let $M \le_{\frak K} N,\|M\| < \mu$ and $p \in {\Cal S}(M)$ be omitted 
in $N$ and let $\theta = \delta + \|M\| + LS({\frak K}) < \mu$, and let 
$p \le q \in {\Cal S}(N)$.  Now we can choose by induction on $i \le \delta,
M_i \le N_i$ and $M^+_i \le_{\frak K} N$ such that $M_i \in K_\theta,M^+_i 
\in K_\theta,M_i$ is $\le_{\frak K}$-increasing continuous and $M \cap N_i 
\subseteq M_i,j < i \Rightarrow M^+_j \cap N_i \subseteq M_{i+1}$ and 
$M_i <^1_{\theta,\omega} M_{i+1}$ and if $q$ does $\theta$-split over 
$M_i$ then $q \restriction M^+_i$ does $\theta$-split over $M_i$.

So by \scite{6.2}, \scite{6.3} we know that 
$M_\delta$ is saturated, and for some $i(*)
< \delta$ we have: $q \restriction M_\delta$ does not $\theta$-split over
$M_{i(*)}$.  But $M^+_{i(*)} \subseteq N = \dbcu_{i < \delta} N_i,
M^+_{i(*)} \cap N_j \subseteq M_{j+1}$ so $M^+_{i(*)} \subseteq M_\delta$.
So necessarily $q \in {\Cal S}(N)$ does not $\theta$-split over $M_{i(*)}$.

Now we choose by induction on $\alpha < \theta^+,M_{i(*),\alpha},b_\alpha,
f_\alpha$ such that: \nl
$M_{i(*),\alpha} \in K_\theta,M_{i(*)} \le_{\frak K} M_{i(*),\alpha}
\le_{\frak K} N_{i(*)},M_{i(*),\alpha}$ is $\le_{\frak K}$-increasing
continuous in $\alpha,b_\alpha \in N_{i(*)}$ realizes $q \restriction
M_{i(*),\alpha},f_\alpha$ is a function with domain $M_\delta$ and range
$\subseteq N_{i(*)}$ such that the sequences
$\bar c = \langle c:c \in M_\delta \rangle$ and $\bar c^\alpha =: 
\langle f_\alpha(c):c \in M_\delta \rangle$ realize the
same type over $M_{i(*),\alpha}$ and $\{b_\alpha\} \cup \text{ Rang}(f_\alpha)
\subseteq M_{i(*),\alpha + 1}$.  As $N_{i(*)}$ is saturated we can carry the
construction; if some $b_\alpha$ realizes $q \restriction M_\delta$ we are 
done (as $b_\alpha \in N$ realizes $p$).  Let $d \in {\frak C}$ realize $q$ so
\mr
\item "{$(*)_1$}"  $\alpha < \beta < \theta^+ \Rightarrow \bar c^\beta
\char 94 \langle b_\alpha \rangle$ does not realize tp$(\bar c \char 94
\langle d \rangle,M_{i(*)},{\frak C})$. \nl
[Why?  As $\bar c \char 94 \langle b_\alpha \rangle$ does not realize
tp$(\bar c \char 94 \langle d \rangle,M_{i(*)},{\frak C})$ because $d$
realizes $p \restriction \bar c$ whereas $b_\alpha$ does not realize
$p \restriction \bar c$.]
\ermn
On the other hand as $q$ does not $\theta$-split over $M_{i(*)}$ we have \nl
tp$(\bar c \char 94 \langle d \rangle,M_{i(*)},{\frak C}) = \text{ tp}
(\bar c^\alpha \char 94 \langle d \rangle,M_{i(*)},{\frak C})$ so by the
choice of $b_\beta$:
\mr
\item "{$(*)_2$}"  if $\alpha < \beta < \theta^+$ then $\bar c^\alpha
\char 94 \langle b_\beta \rangle$ realizes tp$(\bar c \char 94 \langle d
\rangle,M_{i(*)},{\frak C})$.
\ermn
We are almost done by \scite{4.9}. \nl
[Why only almost?  We would like to use the ``$\theta$-order property fail",
now if we could define $\langle \bar c^\beta \char 94 \langle b_\beta
\rangle:\text{ for } \beta < (2^\theta)^+ \rangle$ fine, 
but we have only $\alpha < \theta^+$, this is too short.]  Now we will 
refine the construction to make
$\langle \bar c^\beta \char 94 \langle b_\beta \rangle:\beta < \theta^+
\rangle$ strictly indiscernible which will be enough.
As $N_{i(*)}$ is saturated \wilog \,
$N_{i(*)} = EM_{\tau({\frak K})}(\mu,\Phi)$ and $M_{i(*)} = 
EM_{\tau({\frak K})}(\theta,\Phi)$ (using \scite{6.5} below).  
As before for some $\gamma < \theta^+$
there are sequences $\bar c',\bar b'$ in $EM_{\tau({\frak K})}(\mu + \gamma,
\Phi)$ realizing tp$(\bar c,N_{i(*)},{\frak C}),q \restriction N_{i(*)}$
respectively, here we use cf$(\lambda) > \mu$ rather than just cf$(\lambda)
\ge \mu$.  For each $\beta < \theta^+$ there is a canonical isomorphism
$g_\beta$ from $EM_{\tau(\Phi)}(\beta \cup [\mu,\mu + \gamma),\Phi)$ onto
$EM_{\tau(\Phi)}(\beta + \gamma,\Phi)$.  So \wilog \, $M_{i(*),\alpha} =
EM_{\tau({\frak K})}(\theta + \gamma_\alpha,\Phi),\bar c^\alpha =
g_{\theta + \gamma_\alpha}(\bar c'),b_\alpha = g_{\theta + \gamma_\alpha}
(b')$.  So $(*)_1 + (*)_2$ gives the order property. 
\hfill$\square_{\scite{6.4}}$
\enddemo
\bn
We really proved, in \scite{6.3} (from $\lambda$ categoricity):

\bigskip

\demo{\stag{6.5} Subfact}  Assume $K$ is categorical in $\lambda$. \nl
1) If $I \subseteq J$ are linear order, of power $< \text{ cf}(\lambda)$; 
\medskip
\roster
\item "{$(*)$}"   $t \in J\backslash I \Rightarrow
\biggl( \exists^{\aleph_0} s \in J \biggr)[s \sim_I t]$ where $s \sim_I t$  
means ``$s,t$ realize the same Dedekind cut",
\endroster
\medskip

\noindent
\underbar{then} every type over $EM_{\tau({\frak K})}(I,\Phi)$ 
is realized in $EM_{\tau({\frak K})}(J,\Phi)$. \newline
2)  Adding more Skolem functions we can omit $(*)$,
for a suitable $\Phi$ we can make even the extension $\mu$-saturated over
$EM_{\tau({\frak K})}(I,\Phi)$.
\enddemo
\bigskip

\demo{Proof}  Why?  Use the proof of \scite{6.3}(1).  \newline
Replace the cut of $t$ in $I$ by $\lambda$:  
we get cf$(\lambda)$-saturated model.  \hfill$\square_{\scite{6.5}}$
\enddemo
\newpage

\head {\S7 More on Splitting} \endhead  \resetall
\bigskip

\demo{\stag{7.0} Hypothesis}  As before $+$ conclusions of {\S6} for $\mu \in
[LS({\frak K}),\text{cf}(\lambda))$. \newline
So 
\medskip
\roster
\item "{$(*)(a)$}"  ${\frak K}$ has a saturated model in $\mu$.
\smallskip
\noindent
\item "{$(b)$}"  union of increasing chain of saturated models in $K_\mu$ of 
length $\le \mu$ \nl
is saturated.
\smallskip
\noindent
\item "{$(c)$}"  \underbar{if} $\langle M_i:i \le \delta \rangle$ increasing 
continuous in $K_\mu$, each $M_{i+1}$ saturated over $M_i$ (the previous one),
$p \in {\Cal S}(M_\delta)$ \underbar{then} for some $i < \delta$, $p$ does 
not $\mu$-split over $M_i$.
\endroster
\enddemo
\bigskip

\subhead {\stag{7.1} Conclusion} \endsubhead  If $p \in {\Cal S}^m(M)$ 
and $M \in K_\mu$ is saturated, \underbar{then} for some \nl
$M^- <^1_{\mu,\omega} M,M^- \in K_\mu$ is saturated
and $p$ does not $\mu $-split over $M^-$.
\bigskip

\demo{Proof}  We can find $\langle M_n:n \le \omega \rangle$ in
$K_\mu$, each $M_n$ saturated $M_n \le^1_{\mu,\omega} M_{n+1}$ and
$M_\omega = \dbcu_{n < \omega} M_n$ so as $M_\omega$ is saturated, 
without loss of generality $M_\omega = M$.  Now using $(*)(c)$ of 
\scite{7.0} some $M_n$ is O.K. as $M^-$. 
\hfill$\square_{\scite{7.1}}$
\enddemo
\bigskip

\subhead {\stag{7.2} Fact} \endsubhead  If $M_0 \le^1_{\mu,\omega} M_2 
\le^1_{\mu,\omega} M_3,p \in {\Cal S}^m(M_3)$,  $p$  does not $\mu $-split 
over $M_0$, \underbar{then} $R(p) = R(p\restriction M_2)$.
\bigskip

\demo{Proof}  We can find (by uniqueness) $M_1 \in K_\mu$ such that $M_0 
\le^1_{\mu,\omega} M_1 \le^1_{\mu,\omega} M_2$ and we can find $M_4 \in
K_\mu$ such that $M_3 \le^1_{\mu,\omega} M_4$. \nl
We can find an isomorphism $h_1$ from $M_3$ onto $M_2$ over $M_1$ (by the
uniqueness properties $<^1_{\mu,\omega}$).  By uniqueness there is an
automorphism $h$ of $M_4$ extending $h_1$.  Also by uniqueness there is
$q \in {\Cal S}(M_4)$ which does not $\mu$-split over $M_0$ and extend
$p \restriction M_1$.  As $p,p \restriction M_2$ does not $\mu$-split over
$M_0$ and have the same restriction to $M_1$ and $M_0 \le_{\mu,\omega} M_1$
clearly $p = q \restriction M_2$.  Consider $q$ and $h(q)$ both 
from ${\Cal S}(M_3)$, both do not $\mu$-split over $M_0$ and have the
same restriction to $M_1$; as $M_0 <^1_{\mu,\omega} M_1$ it follows that
$q = h(q)$. \nl
So $R(p \restriction M_1) = R(q \restriction M_1) = R(h(q \restriction M_2))
= R(q \restriction M_2) = R(p)$ as required. \nl
${{}}$ \hfill$\square_{\scite{7.2}}$
\enddemo
\bigskip

\proclaim{\stag{7.3} Claim}  [$K$ categorical in $\lambda$, cf$(\lambda) 
> \mu > LS({\frak K})$]. 

Suppose $m < \omega,M \in K_\mu$ is saturated, $p \in {\Cal S}^m(M),
M \le_{\frak K} N \in K_\mu,p \le q \in {\Cal S}^m(N),N$ saturated 
over $M,q$ not a stationarization of $p$
(i.e. for no $M^- <^\circ_{\mu,\omega} M$, $q$ does not $\mu$-split over  
$M^-$).  \underbar{Then} $q$ does $\mu$-divide over $M$.
\endproclaim
\bigskip

\demo{Proof}  By \scite{7.4} below and \scite{6.2} (just $p$ does not
$\mu$-split over some $N_m$ where \nl
$\langle N_\alpha:\alpha \le \omega \rangle$ witness
$N_0 <^1_{\mu,\omega} M$).
\enddemo
\bigskip

\proclaim{\stag{7.4} Claim}  [Assumptions of \scite{7.4}]
Assume $M_0 <^1_{\mu,\omega} M_1 <^1_{\mu,\omega}
M_2$ all saturated.  If $q \in {\Cal S}(M_2)$ does not $\mu$-split over $M_1$
and $q \restriction M_1$ does not $\mu$-split over $M_0$, \ub{then} $q$
does not $\mu$-split over $M_0$.
\endproclaim
\bigskip

\demo{Proof}  Let $M_3 \in K_\mu$ be such that $M_2 <^1_{\mu,\omega} M_3$ and
$c \in M_3$ realizes $q$.  Choose a linear order $I^*$ such that
$I^* \times (\mu + \omega^*) \cong I^* \cong I^* \times \mu$, remember that
on the product we do not use lexicographic order.  $I^*$ has no
first nor last element \nl
(see \cite[AP]{Sh:220}).

Let $I_0 = I^* \times \mu,I_1 = I_0 + I^* \times \Bbb Z,I_2 = I_1 + I^*
\times \Bbb Z,I_3 = I_2 + I^* \times \mu$.  Clearly \wilog \, 
$M_\ell = EM_{\tau({\frak K})}(\Phi,I_\ell)$, let 
$c = \tau(\bar a_{t_0},\dotsc,a_{t_k})$ so $t_0,\dotsc,t_k \in I_3$; let
$I_{1,n} = I_0 + I^* \times \{m:\Bbb Z \models m < n\}$ 
and $I_{2,n} = I_0 + I^* \times \{m:\Bbb Z \models m < n\}$ and 
$I_{0,\alpha} = \alpha \times I^*$.
So we can find a (negative) integer $n(*)$ small enough and $m(*) \in
\Bbb Z$ large enough such that $\{t_0,\dotsc,t_n\} \cap I_{2,n(*)+1} \subseteq
I_{1,m(*)-1}$.  Let $M_{1,n} = EM_{\tau({\frak K})}(I_{1,n},\Phi)$ 
and $M_{2,n} = EM_{\tau({\frak K})}(I_{2,n},\Phi)$.
Clearly $M_0 <^1_{\mu,\omega} M_{1,n} <^1_{\mu,\omega} M_1 <^1_{\mu,
\omega} M_{2,n} <^1_{\mu,\omega} M_2$.  Clearly (use automorphism of $I_3$)
\mr
\item "{$(*)_0$}"  $q \restriction M_{2,n}$ does not $\mu$-split over
$M_{1,m}$ if $\Bbb Z \models n < n(*),m(*) \le m \in \Bbb Z$.
\ermn
By \scite{7.2} with $q,M_1,M_{2,n},M_2,q$ here standing for 
$M_0,M_2,M_3,p$ there we get
\mr
\item "{$(*)_1$}"  $R(q) = R(q \restriction M_{2,n})$ if $n \in \Bbb Z$.
\ermn
Similarly
\mr
\item "{$(*)_2$}"  $R(q \restriction M_1) = R(q \restriction M_{1,m})$ 
if $m \in \Bbb Z$.
\ermn
By $(*)_0$ and \scite{7.2} we have
\mr
\item "{$(*)_3$}"  $R(q \restriction M_{2,n(*)}) = 
R(q \restriction M_{1,m(*)})$.
\ermn
Similarly we can find a successor ordinal 
$\alpha(*) < \mu$ and $k(*) \in \Bbb Z$ such that

$$
\{t_0,\dotsc,t_k\} \cap I_{1,k(*)+1} \subseteq I_{0,\alpha(*)-1}
$$
\mn
and then prove
\mr
\item "{$(*)_4$}"  $R(q \restriction M_0) = R(q \restriction M_{0,\alpha})$ 
if $\alpha(*) \le \alpha < \mu$
\sn
\item "{$(*)_5$}"  $R(q \restriction M_{1,\ell(*)}) = 
R(q \restriction M_{0,\alpha})$ if $\alpha(*) \le \alpha < \mu$.
\ermn
Together $R(q) = R(q \restriction M_0)$, hence $q$ does not $\mu$-split over
$M_0$ as required.  \hfill$\square_{\scite{7.4}}$
\enddemo
\newpage 

\head {PART II} \endhead \centerline {\S8 Existence\footnote{Done end of Oct.1988} of nice $\Phi$}
\resetall
\bigskip

We build $EM$ models, where ``equality of types over $A$ in the sense of the 
existence of automorphisms over $A$" behaves nicely.
\bigskip

\subhead {\stag{X1.0} Context} \endsubhead
\roster
\item "{(a)}"  ${\frak K}$ is an abstract elementary class with models of 
cardinality $\ge \beth_{(2^{LS({\frak K})})^+}$; it really suffices to assume:
\sn
\item "{(a)$'$}"  ${\frak K}$ is a class of $\tau(K)$-models, which is
$PC_{\kappa^+,\omega}$ with a model of cardinality  
$\ge \beth_{(2^{LS({\frak K})})^+}$.
\endroster
\bigskip

\definition{\stag{X1.1} Definition}  1) Let $\kappa \ge LS({\frak K})$, now
$\varUpsilon^{or}_\kappa = \varUpsilon^{or}_{\kappa,\tau}$ is the family of  
$\Phi$ proper for linear orders (see \cite[Ch.VII]{Sh:c}) such that:
\medskip
\roster
\item "{$(a)$}"  $|\tau(\Phi)| \le \kappa$
\smallskip
\noindent
\item "{$(b)$}"  $EM_{\tau({\frak K})}(I,\Phi) = 
EM(I,\Phi)\restriction \tau (K) \in K$
\smallskip
\noindent
\item "{$(c)$}"  $I \subseteq J \Rightarrow EM_{\tau({\frak K})}(I,\Phi) 
\le_{\frak K} EM_{\tau({\frak K})}(J,\Phi)$.
\endroster
\medskip

\noindent
2) $\varUpsilon^{or}$ is $\varUpsilon^{or}_{LS({\frak K})}$.
\enddefinition
\bigskip

\definition{\stag{X1.2} Definition}   We define partial orders
$\le^\oplus_\kappa$ and $\le^\otimes_\kappa$ on $\varUpsilon^{or}_\kappa$ 
(for $\kappa \ge LS({\frak K}))$: \newline
1)  $\Psi_1 \le^\oplus_\kappa \Psi_2$ \underbar{if} 
$\tau(\Psi_1) \subseteq \tau (\Psi_2)$ and $EM_{\tau({\frak K})}(I,\Psi_1) 
\le_{\frak K} EM_{\tau({\frak K})}(I,\Psi_2)$ and $EM(I,\Psi_1) \subseteq
EM_{\tau(\Psi_1)}(I,\Psi_2)$ and
$EM(I,\Psi_1) = EM_{\tau(\Psi_1)}(I,\Psi_1) \subseteq EM_{\tau(\Psi_1)}
(I,\Psi_2)$ for any linear order $I$. \newline
Again for $\kappa = LS({\frak K})$ we may drop the $\kappa$. \newline
2) For $\Phi_1,\Phi_2 \in \Upsilon^{\text{or}}_\kappa$, we say 
$\Phi_2$ is an inessential extension of $\Phi_1,\Phi_1 \le^{\text{ie}}_\kappa
\Phi_2$ if $\Phi_1 \le^\oplus_\kappa \Phi_2$ and for every linear order $I$,
we have

$$
EM_{\tau({\frak K})}(I,\Phi_1) = EM_{\tau({\frak K})}(I,\Phi_2).
$$

$$
\text{(note: there may be more functions in } \tau(\Phi_2)!)
$$
\mn
3) $\Phi_1 \le^\otimes_\kappa \Phi_2$ \underbar{iff}
there is $\Psi$ proper for linear order and producing linear orders such that:
\medskip
\roster
\item "{(a)}"  $\tau(\Psi)$  has cardinality  $\le \kappa$,
\smallskip
\noindent
\item "{(b)}"  $EM(I,\Psi)$  is a linear order which is an extension of $I$:  
in fact  \newline
$[t \in  I \Rightarrow  x_t = t]$
\smallskip
\noindent
\item "{(c)}"  $\Phi'_2 \le^{\text{ie}}_\kappa \Phi_2$ where
$\Phi'_2 = \Psi \circ \Phi_1$, i.e.
\endroster

$$
EM(I,\Phi'_2) = EM(EM(I,\Psi),\Phi_1).
$$
\medskip

\noindent
(So we allow further expansion by functions definable from earlier ones
(composition or even definition by cases), as long as the number is 
$\le \kappa$). 
\enddefinition
\bigskip

\proclaim{\stag{X1.3} Claim}  1)  
$(\varUpsilon^{or}_\kappa,\le^\otimes_\kappa)$
and $(\varUpsilon^{or}_\kappa,\le^\oplus)$ are partial orders (and 
$\le^\otimes_\kappa \subseteq \le^\oplus_\kappa$). \newline
2)  Moreover, if  $\langle \Phi_i:i < \delta \rangle$ is a 
$\le^\otimes_\kappa$-increasing sequence, $\delta < \kappa^+$, \ub{then} it 
has a $<^\otimes_\kappa$-l.u.b. \,  $\Phi$; $EM^1(I,\Phi) = \dsize 
\bigcup_{i<\delta} EM^1(I,\Phi_i)$.  \nl
3) Similarly for $<^\oplus_\kappa$.
\endproclaim
\bigskip

\proclaim{\stag{X1.3A} Lemma}  1) If $N \le_{\frak K} M,\|M\| \ge
\beth_{(2^\chi)^+},\chi \ge \|N\| + LS({\frak K})$, \ub{then} there is
$\Phi$ proper for linear order such that:
\mr
\item "{$(a)$}"  $EM_{\tau({\frak K})}(\emptyset,\Phi) = N$
\sn
\item "{$(b)$}"  $N \le_{\frak K} EM_{\tau({\frak K})}(I,\Phi)$, moreover \nl
$I \subseteq J \Rightarrow EM_{\tau({\frak K})}(I,\Phi) \le_{\frak K}
EM_{\tau({\frak K})}(I,\Phi)$
\sn
\item "{$(c)$}"  $EM_{\tau({\frak K})}(I,\Phi)$ omits every type 
$p \in {\Cal S}(N)$ which $M$ omits, moreover if $I$ is finite then 
$EM_{\tau({\frak K})}(I,\Phi)$ can be
$\le_{\frak K}$-embedded into $M$.
\endroster
\endproclaim
\bigskip

\demo{Proof}  Straight by \cite[1.7]{Sh:88} or deduce by \scite{4.2} or use
\scite{X1.3B} with $N_1 = N_0$.
\enddemo
\bigskip

\proclaim{\stag{X1.3B} Lemma}  Assume
\mr
\item "{$(a)$}"  $LS({\frak K}) \le \chi \le \lambda$
\sn
\item "{$(b)$}"  $N_0 \le_{\frak K} N_1 \le_{\frak K} M$
\sn
\item "{$(c)$}"  $\|N_0\| \le \chi,\|N_1\| = \lambda$ and $\|M\| \ge
\beth_{(2^\chi)^+}(\lambda)$
\sn
\item "{$(d)$}"  $\Gamma_0 = \{p^0_i:i < i^*_0\} \subseteq {\Cal S}(N_0)$
each $p^0_i$ omitted by $M$
\sn
\item "{$(e)$}"  $\Gamma_1 = \{p^1_i:i < i^*_1 \le \chi\} \subseteq
{\Cal S}(N_1)$ such that for no $i < i^*_i$ any $c \in M$ does $c$ realizes
$p'_i/E_\chi$ [i.e. realizes each $p'_i \restriction M,M \le_{\frak K} N_1,
M \in {\frak K}_{\le \chi}$].
\ermn
\ub{Then} we can find $\langle N'_\alpha:\alpha \le \omega \rangle,\Phi$ 
and $\langle q^1_i:i < i^*_1 \rangle$ such that
\mr
\item "{$(\alpha)$}"  $\Phi$ proper for linear order
\sn
\item "{$(\beta)$}"   $N'_\alpha \in {\frak K}_{\le \chi}$ is $\le_{\frak K}$-
increasing continuous (for $\alpha \le \omega$)
\sn
\item "{$(\gamma)$}"  $N'_0 = N_0$ and $N'_\alpha \le_{\frak K} N_1$
\sn
\item "{$(\delta)$}"  $q^1_i \in {\Cal S}(N'_\omega)$
\sn
\item "{$(\epsilon)$}"  $EM_{\tau({\frak K})}(\emptyset,\Phi)$ is
$N'_\omega$
\sn
\item "{$(\zeta)$}"  for linear order $I \subseteq J$ we have \nl
$EM_{\tau({\frak K})}(I,\Phi) \le_{\frak K} EM_{\tau({\frak K})}(J,\Phi)$
\sn
\item "{$(\eta)$}"  for each $n$, there is a 
$\le_{\frak K}$-increasing sequence
$\langle N_{n,m}:m < \omega \rangle$ with union 
$EM_{\tau({\frak K})}(n,\Phi)$ and a $\le_{\frak K}$-embedding $f_{n,m}$ of
$N_{n,m}$ into $M$ with range $N'_{n,m}$ such that
\sn
{\roster
\itemitem{ $(i)$ }  $N'_m = N'_{0,m}$,
\sn
\itemitem{ $(ii)$ }  $f_{n,m} \restriction N_0$ is the identity,
Rang$(f_{0,m}) \subseteq N_1$
\sn
\itemitem{ $(iii)$ }  $f_{n,m}(q^1_i \restriction N'_m) = p^1_i \restriction
\text{ Rang}(f_n)$ for $i < i^*_1$
\endroster}
\item "{$(\theta)$}"  $EM_{\tau({\frak K})}(I,\Phi)$ omits every $p^0_i$ for
$i < i^*_0$ and omits every $q^1_i$ in a strong sense: for every
$a \in EM_{\tau({\frak K})}(I,\Phi)$ for some $n$ we have \nl
$q^1_i \restriction N'_n \ne \text{ tp}(a,N'_n,EM_{\tau({\frak K})}(I,\Phi))$.
\endroster
\endproclaim
\bigskip

\remark{Remark}  1) So we really can replace $q^1_i$ by $\langle q^1_i
\restriction N'_n:n < \omega \rangle$, but for $\omega$-chains by chasing
arrows such limit $(q^1_i)$ exists. \nl
2) Clause $(\zeta)$ follows from Clause $(\eta)$.
\endremark
\bigskip

\demo{Proof}  By \cite[1.7]{Sh:88} (and see \scite{0.E}) we can 
find $\tau_1,\tau({\frak K}) \subseteq \tau_1,|\tau_1| \le \chi$ 
(here we can have $|\tau_1| \le 
LS({\frak K}) \le \chi$) and an expansion $M^+$ of $M$ to a 
$\tau_1$-model and a set $\Gamma$ of quantifier free types 
(so $|\Gamma| \le 2^{\aleph_0 + |LS({\frak K})|}$) such that:
\mr
\item "{$(A)$}"  $M^+$ omits every $p \in \Gamma$ and if $M^*$ is a
$\tau_1$-model omitting every $p \in \Gamma$ \ub{then} $M^* \restriction \tau
({\frak K}) \in K$ and $N^* \subseteq M^* \Rightarrow N^* \restriction
\tau({\frak K}) \le_{\frak K} M^* \restriction \tau({\frak K})$
\sn
\item "{$(B)$}"  for $\bar a \in {}^{\omega >}M$ we let $M^+_{\bar a} =
M^+ \restriction c \ell(\bar a,M^+)$ then $M^+_{\bar a} \restriction
\tau({\frak K}) \le_{\frak K} M^+ \restriction \tau({\frak K})$,
Rang$(\bar a) \subseteq \text{ Rang}(\bar b) \Rightarrow M^+_{\bar a}
\restriction \tau({\frak K}) \le_{\frak K} M^+_b \restriction \tau({\frak K})$
where \nl
$\bar a \in {}^{\omega >}(N_\ell) \Rightarrow |M^+_{\bar a}| \subseteq
N_\ell$.
\ermn
Note:  Further expansion of $M^+$ to $M^*$, as long as $|\tau(M^*)| \le \chi$
preserves $(A) + (B)$ so we can add
\mr
\item "{$(C)$}"  $N_0,M^+_{\langle \rangle}$ have the same universe \nl
and let for $\ell =0,1, \, 
M^+_{\bar a,\ell} = M^+_{\bar a} \restriction (|N_\ell| \cap |M^+_{\bar a}|),
\ell = 1,2$
\sn
\item "{$(D)$}"  $M^+_{\bar a,0} \restriction \tau({\frak K}) \le_{\frak K}
M^+_{\bar a,1} \restriction \tau({\frak K}) \le_{\frak K} M^+_{\bar a}
\restriction \tau({\frak K})$
\sn
\item "{$(E)$}"  for $i < i^*_1$, the type $p^1_i \restriction 
(M^+_{\bar a,1} \restriction \tau({\frak K}))$ is not realized in
$M^+_{\bar a} \restriction \tau({\frak K})$.
\ermn
Now we choose by induction on $n$, sequence $\langle f^n_\alpha:\alpha <
(2^\chi)^+ \rangle$ and $N'_n$ such that:
\mr
\widestnumber\item{$(iii)$}
\item "{$(i)$}"  $f^n_\alpha$ is a one-to-one function from $\beth_\alpha
(\lambda)$ into $M$
\sn
\item "{$(ii)$}"  $\langle f^n_\alpha(\zeta):\zeta < \beth_\alpha(\lambda)
\rangle$ is $n$-indiscernible in $M^+$
\sn
\item "{$(iii)$}"  moreover, if $\alpha,\beta < (2^\chi)^+$, and $m \le n$
and $\zeta_1 < \ldots < \zeta_m < \beth_\alpha(\lambda)$ and $\xi_1 < \ldots
< \xi_m < \beth_\beta(\lambda)$ \ub{then}: the sequences $\bar a =
\langle f^n_\alpha(\zeta_1),\dotsc,f^n_\alpha(\zeta_m) \rangle$, \nl
$\bar b = \langle f^m_\beta(\xi_1),\dotsc,f^m_\beta(\xi_m) \rangle$ 
realize the same quantifier free type in $M^+$ over $N^+_1$, so 
there is a natural isomorphism $g_{\bar b,\bar a}$
from $M^+_{\bar a}$ onto $M^+_{\bar b}$ (mapping $f_\alpha(\zeta_\ell)$ 
to $f_\beta(\xi_1)$), moreover
$$
i < i^*_1 \Rightarrow g_{\bar b,\bar a}(p^1_i \restriction (M^+_{\bar a,1}
\restriction \tau({\frak K}))) = p^1_i \restriction (M^+_{\bar b,1}
\restriction \tau({\frak K}))
$$
and
$$
N'_m = M^+_{\bar a,1} \restriction \tau({\frak K}).
$$
\ermn
The rest should be clear.  \hfill$\square_{\scite{X1.3B}}$
\enddemo
\bigskip

\proclaim{\stag{X1.4} Claim}  Suppose
\roster
\item "{(a)}"  $\Phi \in \varUpsilon^{or}_\kappa$
\smallskip
\noindent
\item "{(b)}"  $n < \omega$, $u$, $u_1$, $u_2$ are subsets of
$\{0,1,...,n-1\}$ and $\sigma_1(...,\bar x_\ell,...)_{\ell \in u_1}$, \newline
$\sigma_2(...,\bar x_\ell,...)_{\ell \in u_2}$ are $\tau(\Phi)$-terms. 
\smallskip
\noindent
\item "{(c)}"  for every $\alpha < (2^{LS({\frak K})})^+$ (or at least
$\beth_\alpha < \mu(\kappa)$ - see \cite[Ch.VII,\S4]{Sh:c} but for this we
should be careful as to omit only $\le LS({\frak K})$ types) there are linear 
orders $I \subseteq J,I$ \,\, $\aleph_0$-homogeneous 
inside \footnote{this means that every partial order preserving function 
$h$ from $I$ to $I$ can be extended to an automorphism of $J$.}
$J,I$ of cardinality $\ge \beth_\alpha$, such that for some (equivalently 
every) $t_0 < t_1 < \dots < t_{n-1}$ of $I$ we have: 
{\roster
\itemitem{ $\oplus$ }  for some automorphism $f$ of $EM_{\tau({\frak K})}
(J,\Phi)$, \newline
$f \restriction EM_\tau(I\backslash \{t_\ell:\ell < n,\ell \notin u\},\Phi)$ 
is the identity and \newline
$f\biggl( \sigma_1(...,\bar a_{t_\ell},...)
_{\ell \in u_1} \biggr) = \sigma_2(...,\bar a_{t_\ell},...)_{\ell \in u_2}$.
\endroster}
\endroster
\medskip

\noindent
THEN for some $\Phi'$,  $\Phi \le^\oplus_\kappa \Phi'$ and even $\Phi
\le^\otimes_\kappa \Phi'$ we have 
\medskip
\roster
\item "{$\otimes$}"  for every linear order $I$ and $t_0 < \dots < t_{n-1}$ 
from $I$, there is an automorphism $f$ of $EM_\tau(I,\Phi')$ such that:
{\roster
\itemitem{ $(\alpha)$ }  $f \restriction EM(I \backslash \{t_\ell:\ell < n,
\ell \notin u\},\Phi')$ is the identity and
\itemitem{ $(\beta)$ }  $f\biggl( \sigma_1(...,\bar a_{t_\ell },...)
_{\ell \in u_1} \biggr) = \sigma_2(...,\bar a_{t_\ell},...)_{\ell \in u_2}$
\itemitem{ $(\gamma)$ }  $f = F(-,\bar a_{t_0},\dotsc,\bar a_{t_{n-1}})$ 
for some $F \in \tau(\Phi')$.
\endroster}
\endroster
\endproclaim
\bigskip

\demo{Proof}  Expand $M = EM(J,\Phi)$ by the predicates 
$Q_1 = \{\bar a_t:t \in I\},Q_2 = \{\bar a_t:t \in J\}$ if $\alpha = \ell g
(\bar a_t)$ is finite, in any case we use $Q_{\ell,i,j} = 
\{(a_{t,i},a_{t,j}):t \in I\}$ for $\ell \in \{1,2\}$ and $i \le j < \alpha$;
and \wilog \, $t \ne s \Rightarrow a_{t,0} \ne a_{s,0}$ and we identify
$t \in J$ with $a_{t,0}$.
For $t_0 < \ldots < t_{n-1} \in I$, let
$f_{t_0,\dotsc,t_{n-1}} \in AUT(EM_{\tau({\frak K})}(J,\Phi))$ 
be as in $(\oplus)$ and let
$g_\ell$ (for $\ell < \omega$) be functions from $M$ into 
$\{ \bar a_t:t \in J\}$ such that $\forall x \in M,x = \sigma_x(g_0(x) \ldots 
g_{n-1}(x))$ such that $g_\ell(x) <^J g_{\ell +1}(x)$ if $\ell < n-1$ and
$g_\ell(x) = g_{\ell +1}(x)$ otherwise.  Lastly, let $P_\sigma = 
\{x:\sigma_x = \sigma\}$.

Let $F$ be an $(n+1)$-ary function, 
$F(\bar a_{t_0},\dotsc,\bar a_{t_{n-1}},b) = f_{t_0 \ldots t_{n-1}}(b)$ 
when defined.  The model we get we call $M^+$. 
Now use the omitting types theorem, e.g. \scite{X1.3A}.  So there is a model
$N^+$ and $\langle \bar b_n:n < \omega \rangle$ indiscernible in it such
that $N^+ \equiv M^+,N^+$ omits all types which $M^+$ omits, for every
$m< \omega$ for some $s_0 < \ldots < s_{n-1}$ from $I$ the type of
$\bar b_0 \char 94 \ldots \bar b_{n-1}$ in $N^+$ is equal to the type of
$\bar a_{s_0} \char 94 \ldots \char 94 \bar a_{s_{n-1}}$ in $M^+$.  Define
$\Phi'$ such that $EM(I^*,\Phi')$ is a $\tau(N^+)$-model generated by 
$\{\bar a_t:t \in I^*\}$ such that $t_0 < \ldots < t_{n-1} \in I^* 
\Rightarrow$ type of $\bar a_{t_0} \char 94
\ldots \char 94 \bar a_{t_{n-1}}$ in $EM(I^*,\Phi')$ is equal to type of
$\bar b_0 \char 94 \ldots \char 94 \bar b_{n-1}$ in $N^+$.
\medskip

\noindent
Why is $\Phi <^\otimes \Phi'$ and not just $\Phi <^\oplus \Phi'$? \newline
Here we use \footnote{if $\ell g(\bar a_t)$ is infinite, slightly more
complicated} $Q_1,Q_2$ in $M^+$ we have
\medskip
\roster
\item "{$(*)$}"  every $c \in M^+$ is in the $\tau_\Phi$-Skolem Hull of
$Q^{M^+}_2 = \{\bar a_t:t \in J\}$.
\endroster
\medskip

\noindent
So 
\roster
\item "{$(*)'$}"  $M^+$ omits the type
$$
p(x) = \bigl\{ \neg(\exists \bar y_0,\dotsc,\bar y_{n-1})
(\dsize \bigwedge_{\ell < n} Q_2(\bar y_\ell) \and x = \sigma(y_0,\dotsc,
y_n):\sigma \in \tau_\Phi) \bigr\}.
$$
${}$ \hfill$\square_{\scite{X1.4}}$
\endroster
\enddemo
\bigskip

\demo{\stag{X1.5} Conclusion}  For $\kappa \ge LS({\frak K})$ 
there is $\Phi^* \in  
\varUpsilon^{or}_\kappa$ (in fact for every $\Phi \in \varUpsilon^{or}_\kappa$
there is $\Phi^*,\Phi \le^\otimes_\kappa \Phi^* \in \varUpsilon^{or}_\kappa$)
satisfying:
\medskip
\roster
\item "{$(a)$}"  if  $\Phi^\ast$ satisfies the assumptions of \scite{X1.4} 
for some $I,J$ (playing the role of $\Phi$ there) \underbar{then} it 
satisfies its conclusion (i.e. playing the role of $\Phi'$ there)
\smallskip
\noindent
\item "{$(b)$}"  moreover if $\kappa \ge 2^{LS({\frak K})}$, for some 
$\chi(\Phi^*) < \mu(\kappa)$ (see \cite[Ch.VII,\S4]{Sh:c}), we can 
weaken the assumption 
$\alpha < \beth_{(2^\kappa)^+}$ to $\beth_\alpha \le \chi(\Phi^*)$
\smallskip
\noindent
\item "{$(c)$}"  moreover, in \scite{X1.4} we can omit ``$I$ is 
$\aleph_0$-homogeneous inside $J$"
\smallskip
\noindent
\item "{$(d)$}"  also we can replace clause $(\alpha)$ of $\otimes$ (of
\scite{X1.4}) by: $f$ extends some automorphism of $EM(I \backslash \{t_\ell:
\ell < n,\ell \notin u\},\Phi^*)$ definable as in clause $(\gamma)$ of
$\otimes$ of \scite{X1.4}.
\sn
\item "{$(e)$}"   we can deal similarly with automorphisms extending a given 
\newline
$f \restriction EM(I \restriction \{t_\ell:\ell < n\})$ and having finitely
many demands.
\endroster
\enddemo
\bigskip

\demo{Proof}  For (a) iterate \scite{X1.4}, by bookkeeping looking at all  
$\langle \sigma_1,\sigma_2,u,u_1,u_2 \rangle$ and use \scite{X1.3} for noting
that the iteration is possible.  Now (b) holds as cf$(\mu(\kappa)) > \kappa$, 
and the number of terms is $\le \kappa$.  For (c) we can let $\Psi$ be such 
that $EM(I,\Psi)$ is an $\aleph_0$-homogeneous linear order, $|\tau(\Psi)| = 
\aleph_0$ and use $\Psi \circ \Phi^*$.  The rest are easy, too. \nl
${{}}$  \hfill$\square_{\scite{X1.5}}$
\enddemo
\bigskip

\proclaim{\stag{X1.6} Lemma}  Let $\Phi^*$ be as in \scite{X1.5}, and $I$ be a
linear order of cardinality $\chi(\Phi^*)$ (where $\chi(\Phi^*)$ is from
\scite{X1.5}).  Assume $\sigma(\bar x_0,\dotsc,\bar x_{n-1})$ is a term in 
$\tau(\Phi^*)$, for $\ell = 1,2$ we have
$t^\ell_0 < \dots < t^\ell_{n-1}$ and $u \subseteq \{\ell:t^1_\ell = 
t^2_\ell \}$, and there is no automorphism $f$ of $EM(I,\Phi^*)$ such that
$f \restriction EM(I\backslash \{t^1_i,t^2_i:\ell < n,\ell \notin u\},\Phi^*)$
is the identity, and $f \bigl( \sigma \bigl( \bar a_{t^1_0},\dotsc,\bigr)
\bigr) = \sigma \bigl( \bar a_{t^2_0},\dotsc,\bigr)$. \newline
\underbar{Then}
\medskip
\roster
\item  for  $\chi > \chi(\Phi^*)^+$ we have $I(\chi,K) = 2^\chi$.
\smallskip
\noindent
\item  We have the $\chi(\Phi^*)$-order property in the sense of Definition
\scite{4.1B} (see more \cite[Ch.III,\S3]{Sh:300} or better 
\cite[Ch.III,\S3]{Sh:e}.)
\endroster
\endproclaim
\bigskip

\demo{Proof}  Without loss of generality $I$ is dense.

We can find $t^3_0 < \ldots < t^3_{n-1}$ such that

$$
\ell \in u \Rightarrow t^3_\ell = t^1_\ell,
$$

$$
\ell \notin u \Rightarrow t^3_\ell \notin \{t^1_m,t^2_m:m < n\}.
$$
\mn
Now
\mr
\item "{$\bigotimes_1$}"  there is no automorphism $f$ of $EM(I,\Phi^*)$
such that \nl
$f \restriction EM_\tau(I \backslash \{t^1_\ell,t^2_\ell,
t^3_\ell:\ell < n,\ell \notin u\},\Phi)$ is the identity and \nl
$f(\sigma(\bar a_{t^1_0},\ldots)) = \sigma(\bar a_{t^2_0},\ldots)$ \nl
[Why?  If there is, easily some $\Phi$ contradicts \scite{X1.5}(a)]
\sn
\item "{$\bigotimes_2$}"   for some $k \in \{1,2\}$, there is no automorphism
$f$ of $EM(I,\Phi)$ which is the identity of $EM(I \backslash \{t^1_\ell,
t^2_\ell,t^3_\ell:\ell < n,\ell \notin u\},\Phi)$ and 
$f(\sigma(\bar a_{t^k_0},\ldots)) = \sigma(\bar a_{t^3_0},\ldots)$ \nl
[Why?  If not such $f_1,f_2$ exists and $f^{-1}_2 \circ f_1$ contradict
$(*)_2$].
\sn
\item "{$\bigotimes_3$}"  for some $k \in \{1,2\}$ there is no automorphism
$f$ of $EM(I,\Phi)$ which is the identity on $EM(I \backslash \{t^k_\ell,
t^3_\ell:\ell < n,\ell \notin u\},\Phi)$ and 
$f(\sigma(\bar a_{t^k_0},\ldots)) = \sigma(a_{t^3_0},\ldots)$ \nl
[Why?  We negate a stronger demand than in $(*)_2$].
\ermn
By renaming we get that \wilog

$$
t^1_\ell = t^2_k \Rightarrow \ell = k \in u.
$$
\mn
By the transitivity of ``there is an automorphism" we can
assume that just for a singleton $\ell(*),t^1_{\ell(*)} \ne t^2_{\ell(*)}$.  
Now if we increase $u$, surely such isomorphism does not exist so without 
loss of generality $u = \{\ell < n:\ell \ne \ell(*)\}$ and 
$t^1_{\ell(*)} <_I t^2_{\ell(*)}$, by symmetry.  Let $I^0 = \{t \in I:t <_I 
t^1_{\ell(*)}\},I^1 = \{t \in I:t^1_{\ell(*)} \le_I t <_I t^2_{\ell(*)}\},
I^2 = \{t \in I:t^2_{\ell(*)} <_I t\}$ (yes: $<_I$ not $\le_I)$. \newline
Now for every linear order $J$ we can define $I(J)$ as follows:  $I(J)$ is
a linear order which is the sum  $I^0 + \dsize \sum_{t\in J}I^1_t(J) + I^2$,  
$I^1_t(J)$ is isomorphic to $I^1$, so let $f_t:I^1 \rightarrow I^1_t(J)$  
be such an isomorphism.  Let $\bold{\bar b}^t$ list  
$EM(I^0 + I^1_t(J) + I^2)$ (such that for $t,s,\left( \text{id}_{I^0} + 
f_sf^{-1}_t + \text{ id}_{I^2}\right)$  induces a
mapping from $\bold{\bar b}^t$ onto $\bold{\bar b}^s$).  Let $\bar c^t =
f_t\left( \sigma (t^1_0,\dotsc,t^1_{n-1})\right)$.  Now 
\medskip
\roster
\item "{$(*)_1$}"  if  $s_0 <_J r <_J s_1$ \ub{then} there is no
automorphism  $f$  of  $EM_\tau(I(J),\Phi^*)$ over $\bar{\bold b}^r$ 
mapping $\bar c^{s_0}$ to $\bar c^{s_1}$,
\smallskip
\noindent
\item "{$(*)_2$}"  if $J$ is $\aleph_0$-homogeneous (or just 2-transitive)
and  $r <_J s_0 \and r <_J s_1$ or $s_0 <_J r \and s_1 <_J r$ \underbar{then}
there is an automorphism $f$ of $EM_\tau(I(J),\Phi^*)$ over 
${\bold b}^r$ mapping $\bar c^{s_0}$ to $\bar c^{s_1}$.
\endroster
\medskip

\noindent
So by \cite[Ch.III,\S3]{Sh:e} (or earlier version \cite[Ch.III,\S3]{Sh:300}), we
have the order property for sequences of length $\chi(\Phi^*)$; the formula
appearing in the definition of the order is preserved by automorphisms of the
model; though it looks as second order, it does not matter.  So conclusion (2)
holds and (1) follows.  \hfill$\square_{\scite{X1.6}}$
\enddemo
\bigskip

\proclaim{\stag{X1.7} Claim}  Assume
\medskip
\roster
\item "{$(a)$}"  $K$ is categorical in $\lambda$
\smallskip
\noindent
\item "{$(b)$}"  the $M \in K_\lambda$ is $\chi^+$-saturated (holds if
cf$(\lambda) > \chi$)
\smallskip
\noindent
\item "{$(c)$}"  $\chi \ge LS({\frak K})$.
\endroster
\medskip

\noindent
\underbar{Then} every $M \in K$ of cardinality $\ge \beth_{(2^\chi)^+}$
(or just $\ge \beth_{\mu(\chi)}$ if $\chi \ge 2^{LS({\frak K})}$) is
$\chi^+$-saturated.
\endproclaim
\bigskip

\demo{Proof}  If $M$ is a counterexample, let $N \le_{\frak K} M,\|N\| \le
\chi$ and $p \in {\Cal S}(N)$ be omitted by $N$.  By the omitting type theorem
for abstract elementary classes (see \scite{X1.3A}, i.e. \cite{Sh:88}), we
get $M' \in K_\lambda,N \le_K M',M'$ omitting $p$ a contradiction.
\hfill$\square_{\scite{1.7}}$
\enddemo
\bigskip

\proclaim{\stag{X1.8} Claim}  Assume
\medskip
\roster
\item "{$(a)$}"  $LS({\frak K}) \le \chi$
\smallskip
\noindent
\item "{$(b)$}"  for every $\alpha < (2^\chi)^+$ there are $M_\alpha
<_{\frak K} N_\alpha$ (so $M_\alpha \ne N_\alpha$), $\|M_\alpha\| \ge
\beth_\alpha$ and $p \in {\Cal S}(M_\alpha)$ such that $c \in N_\alpha
\Rightarrow \neg pE_\chi \text{ tp}(c,M_\alpha,{\frak C})$.
\endroster
\medskip

\noindent
1) For every $\theta > \chi$ there are $M <_{\frak K} N$ in $K_\theta$
and $p \in {\Cal S}(M_\alpha)$ as in clause (b). \newline
2) Moreover, if $\Phi$ is proper for orders as usual, 
$|\tau(\Phi)| \le \chi,\beth_{(2^\chi)^+} \le \lambda,K$ 
categorical in $\lambda$ and $I$ a linear
order of cardinality $\theta$, \underbar{then} we can demand
$M = EM_{\tau({\frak K})}(I,\Phi)$.
\endproclaim
\bigskip

\demo{Proof}  Straight.
\enddemo
\newpage

\head {\S9 Small Pieces are Enough and Categoricity} \endhead  \resetall
\bigskip

\demo{\stag{X2.0} Context}
\roster
\widestnumber\item{(iii)}
\item "{(i)}"  ${\frak K}$  an abstract elementary class
\smallskip
\noindent
\item "{(ii)}"  $K$  categorical in $\lambda,\lambda > LS({\frak K})$
\smallskip
\noindent
\item "{(iii)}"  some ($\equiv$ any) $M \in K_\lambda$ is saturated
(if $\lambda$ is regular this holds) 
\sn
\item "{$(iv)$}"  $\Phi^*$ is as in \scite{X1.5}.
\endroster
\enddemo
\bigskip 

\noindent
Hence 
\demo{\stag{X2.1} Fact}  For $\mu \in [LS({\frak K}),\lambda)$, 
there is a saturated model of cardinality $\mu$, \nl
(why?  by \scite{6.3}(3)) and there is also $\Phi^* \in 
\varUpsilon^{or}_\mu$ as in \scite{X1.5}.
\enddemo
\bigskip

\proclaim{\stag{X2.2} Main Claim}  If $M \in K$ is a saturated model of 
cardinality $\chi$, \nl
$\chi(\Phi^*) < \chi < \text{ cf}(\lambda) \le \lambda$ \ub{then}
${\Cal S}(M)$ has character $\le \chi(\Phi^*)$, i.e. if $p_1 \ne p_2$ 
are in ${\Cal S}(M)$ then for some $N \le_{\frak K} M,N \in K_{\chi(\Phi^*)}$ 
we have $p_1 \restriction N \ne p_2 \restriction N$.
\endproclaim
\bigskip

\demo{Proof}  We can find $I \subseteq J,|I| = \chi,|J| = \lambda,M = 
EM_{\tau({\frak K})}(I,\Phi^*) \le_{\frak K} N^* = EM_{\tau({\frak K})}
(J,\Phi^*)$ and $I,J$ are $\aleph_0$-homogeneous.  So by \scite{6.4}: every 
$p \in {\Cal S}(M)$ is realized in $N^*$ and say $p$ is realized by
$\sigma_p(\bar a_{t_{p,0}},\bar a_{t_{p,1}},\dotsc,\bar a_{t_{p,n_p-1}})$
where $t_{p,0} < t_{p,1} < \ldots < t_{p,n_p-1}$.  If the conclusion fails,
then we can find $q \ne p$ in ${\Cal S}(M)$ such that
\medskip
\roster
\item "{$(*)$}"  $N \le_{\frak K} M,\|N\| \le \chi(\Phi^*) \Rightarrow p 
\restriction N = q \restriction N$.
\endroster
\medskip

\noindent
Choose $I' \subseteq J,|I'| = \chi(\Phi^*)$ such that $I \subseteq I'$ and
$\{t_{p,\ell}:\ell < n_p\} \subseteq I'$ and $\{t_{q,\ell}:\ell < n_q\} 
\subseteq I'$ and let $M' = EM_{\tau({\frak K})}(I' \cap I,\Phi^*) 
\le_{\frak K} M$.

So $p \restriction M' = q \restriction M'$; so \scite{X1.4} becomes
relevant (i.e. \scite{X1.5}(b)) considering the $\aleph_0$-homogeneity of
$J$) hence by the choice of $\Phi^*,p = q$ contradiction.
\hfill$\square_{\scite{X2.2}}$
\enddemo
\bigskip

\demo{\stag{X2.3} Conclusion}  Let $I$ be a directed partial order.  
Assume $M \in K_\chi$ is saturated, 
$\chi(\Phi^*) \le \chi < \lambda,\langle M_t:t \in I \rangle$
is a $\le_{\frak K}$-increasing family of $\le_{\frak K}$-submodels of $M$, 
each saturated and $[t < s \Rightarrow M_t$ saturated over $M_s]$ and 
$\|M_t\| \le \chi(\Phi^*)$,
\underbar{then} for every $p \in {\Cal S} 
\left( \dsize \bigcup_{t \in I} M_t \right)$
for some $t^* \in I$:
\medskip
\roster
\item "{$(*)$}"  $p$ does not $\mu$-split over $M_{t^*}$ \newline
(and even does not $\chi$-split over $M_{t^*}$).
\endroster
\enddemo
\bigskip

\demo{Proof}  Clear by the proof of \scite{X2.2}.
\enddemo
\bigskip

\proclaim{\stag{X2.3A} Claim}  If $T$ is categorical in $\lambda,LS({\frak K}) 
\le \chi(\Phi^*) \le \mu < \lambda$ and $\langle M_i:i < \delta \rangle$ an 
$<_{\frak K}$-increasing sequence of $\mu^+$-saturated models \ub{then}
$\dsize \bigcup_{i<\delta} M_i$ is 
$\mu^+$-saturated.
\endproclaim
\bigskip

\remark{Remark}  1) Hence this holds for limit cardinals $> LS({\frak K})$.
\nl
2) The addition compared to \scite{6.4} is the case cf$(\lambda) = \mu^+$, 
e.g. $\lambda = \mu^+$.
\endremark
\bigskip

\demo{Proof}  Let $M_\delta = \dsize \bigcup_{i < \delta} M_i$ and assume
$M_\delta$ is not $\mu^+$-saturated.  So there are $N \le_{\frak K} M_\delta$
of cardinality $\le \mu$ and $p \in {\Cal S}(N)$ which is not realized in
$M_\delta$.  Choose $p_1 \in {\Cal S}(M_\delta)$ such that $p_1 \restriction
N = p$. \nl
Without loss of generality $N$ is saturated.

Let $\chi = \chi(\Phi^*)$, \wilog \, $\delta = \text{ cf}(\delta)$. \nl
We claim
\mr
\item "{$\bigotimes$}"  there are $i(*) < \delta$ and 
$N^* \le_{\frak K} M_{i(*)}$
of cardinality $\chi$ such that $p$ does not $\chi$-split over $N^*$. \nl
Why?  Assume toward contradiction that this fails.  The proof of $\bigotimes$
splits to two cases.
\ermn
\ub{Case I}:  cf$(\delta) \le \chi$.

We can choose by induction on $i < \delta,N^0_i,N^1_i$ such that
\mr
\item "{$(a)$}"  $N^0_i \le_{\frak K} M_i$ has cardinality $\chi$
\sn
\item "{$(b)$}"  $N^0_i$ is increasingly continuous
\sn
\item "{$(c)$}"  $N^0_i <^1_{\chi,\omega} N^0_{i+1}$
\sn
\item "{$(d)$}"  $N^0_i \le_{\frak K} N^1_i \le_{\frak K} M_\delta$
\sn
\item "{$(e)$}"  $N^1_i$ has cardinality $\le \chi$
\sn
\item "{$(f)$}"  $p_1 \restriction N^1_i$ does $\chi$-split over $N^0_i$
\sn
\item "{$(g)$}"  for $\varepsilon,\zeta < i,N^1_\varepsilon \cap M_\zeta 
\subseteq N^0_i$.
\ermn
There is no problem to carry the definition and then $N^1_i \subseteq
\dbcu_{j < \delta} N^0_j$ and \nl
$\langle N^0_i:i < \delta \rangle,p_1 \restriction \dbcu_{i < \delta} N^0_i$ 
contradicts \scite{6.2}.
\bn
\ub{Case II}:  cf$(\delta) > \chi$.

Now by \scite{3.2}
\mr
\item "{$(*)$}"  for some $N^* \le_{\frak K} M_\delta$ of cardinality $\chi$
we have $p_1$ does not $\chi$-split over $N^*$.
\ermn
As $\delta = \text{ cf}(\delta) \ge \mu > \chi$, for some $i(*) < \delta$ we
have $N^* \le_{\frak K} M_{i(*)}$.  This ends the proof of $\bigotimes$. \nl
So $i(*),N^*$ are well defined.  Without loss of generality $N^*$ is
saturated.  Let $c \in {\frak C}$ realize $p_1$.  We
can find a $\le_{\frak K}$-embedding $h$ of $EM_{\tau({\frak K})}
(\mu^+ + \mu^+,\Phi)$ into ${\frak C}$ such that
\medskip
\roster
\item "{$(a)$}"  $N^*$ is the $h$-image of $EM_{\tau({\frak K})}(\chi,\Phi)$
\smallskip
\noindent
\item "{$(b)$}"  $h$ maps $EM_{\tau({\frak K})}(\mu^+,\Phi)$ onto 
$M' \le_{\frak K} M_{i(*)}$
\smallskip
\noindent
\item "{$(c)$}"  every $d \in N$ and $c$ belong to the range of $h$.
\endroster
\medskip

By renaming, $h$ is the identity, clearly for some $\alpha < \mu^+$ we
have \newline
$N \cup \{c\} \subseteq EM_\tau(\alpha \cup [\mu^+,\mu^+ + \alpha))$, so
some list $\bar b$ of the members of $N$ is \newline
$\bar \sigma(\ldots,\bar a_i,\dotsc,a_{\mu^+ + j})_{i < \alpha,j <\alpha}$
(assume $\alpha > \mu$ for simplicity) and \newline
$c = \sigma^*(\dotsc,\bar a_i,\dotsc,a_{\mu^+ + j},\ldots)_{i \in u,j \in w}$
($u,w \subseteq \mu^+$ finite as, of course, only finitely many $\bar a_i$'s 
are needed for the term $\sigma^*$).

For each $\gamma < \mu^+$ we define 
$\bar b^\gamma = \bar \sigma(\ldots,\bar a_i,\dotsc,
a_{(1 + \gamma)\alpha + j},\ldots)_{i < \alpha,j < \alpha}$ and
\newline
$c^\gamma = \sigma^*(\ldots,\bar a_i,\dotsc,a_{(1+\gamma)\alpha +j},
\ldots)_{i,j}$ and stipulate $\bar b^{\mu^+} = \bar b,c^{\mu^+} = c$ 
and let $q = \text{ tp}(\bar b \char 94 c,N^*,{\frak C})$.
Clearly $\langle \bar b^\gamma \char 94 c^\gamma:\gamma < \mu^+ \rangle 
\char 94 \bar b \char 94 c$ is a strictly indiscernible sequence over $N^*$
and $\subseteq M_\delta \cup \{c\}$, so also $\{\bar b^\gamma:\gamma \le
\mu\} \subseteq M_\delta$ is strictly indiscernible over $N$.  \nl
[Why?  Use $I \supseteq \mu^+ + \mu^+$ which is a strongly $\mu^{++}$
saturated dense linear order and use automorphisms of EM$(I,\Phi)$ induced
by an automorphism of $I$.]  \nl
As $c$ realizes $p_1$ clearly
tp$(c,M_\delta)$ does not $\chi$-split over $N^*$ hence by \scite{X2.2}
necessarily tp$(\bar b^\gamma \char 94 c,N^*,{\frak C})$ is the same for all
$\gamma \le \mu^+$, hence $\gamma < \mu^+ \Rightarrow \text{ tp}(\bar b^\gamma
\char 94 c^{\mu^+},N^*,{\frak C}) = q$, so by the indiscernibility
$\beta < \gamma \le \mu^+ \Rightarrow \text{ tp}(\bar b^\beta \char 94
c^\gamma,N^*,{\frak C}) = q$.

Similarly for some $q_1$,

$$
\beta < \gamma \le \mu^+ \Rightarrow \text{ tp}(\bar b^\gamma \char 94
c^\beta,N^*,{\frak C}) = q_1.
$$
\mn
If $q \ne q_1$, then tp$(c_0,\bar b^1,{\frak C}) \ne \text{ tp}(c_2,
\bar b^1,{\frak C})$, but Rang$(\bar b^\gamma)$ is a model $N^*_\gamma
\le_{\frak K} M_{i(*)},N^* \le_{\frak K} N^*_\gamma$, so by \scite{X1.6},
for some $v \subseteq \ell g(\bar b^\gamma)$ of cardinality $\chi$, \nl
tp$(c_0,\bar b^1 \restriction \nu,{\frak C}) \ne \text{ tp}(c_2,\bar b^1
\restriction v,{\frak C})$.  So clearly we get the $(\chi,\chi,1)$-order
property contradiction to \scite{4.9}.  \nl
Hence necessarily $\beta \le \mu^+ \and \gamma \le \mu^+ \and \beta \ne \gamma
\Rightarrow \text{ tp}(\bar b^\beta \char 94 c^\gamma,N^*,{\frak C}) = q$.
For $\beta = \mu^+,\gamma =0$ we get that $c^\gamma \in M_{i(*)} \le_K
M_\delta$ realizes tp$(c^\gamma,N,{\frak C}) = p_1 \restriction N$ as
desired.  \hfill$\square_{\scite{X2.3A}}$
\bn
We could have proved earlier \nl
\demo{\stag{X2.3B} Observation}  1) If $M$ is $\theta$-saturated, $\theta >
LS({\frak K})$  and $\theta < \lambda$
and $N \le_{\frak K} M,N \in K_{\le \theta}$ \ub{then} there is
$N',N \le_{\frak K} N' \le_{\frak K} M,N' \in K_\theta$ and every $p \in
{\Cal S}(N)$ realized in $M$ is realized in $N'$. \nl
2) Assume $\langle N_i:i \le \delta \rangle$ is $\le_{\frak K}$-increasingly
continuous, $\delta < \theta^+$ is divisible by $\theta,N_i \in 
K_{\le \theta},N_i \le_{\frak K} M,M$ is $\theta$-saturated,
and every $p \in {\Cal S}(N_i)$ realized in $M$ is realized in $N_{i+1}$
\ub{then} \mr
\item "{$(a)$}"  if $\delta = \theta \times \sigma$, LS$({\frak K}) <
\sigma = \text{ cf}(\sigma) \le \theta$, then $N_\delta$ is $\sigma$-saturated
\sn
\item "{$(b)$}"  if $\delta = \theta \times \theta,\theta > LS({\frak K})$,
\ub{then} $N_\delta$ is saturated.
\endroster
\enddemo
\bigskip

\proclaim{\stag{X2.4} Theorem}  (The Downward \L os theorem for $\lambda$ 
successors).

If $\lambda$ is successor $\ge \mu(\chi(\Phi^*)) = \mu < \chi < \lambda$,
\ub{then} $K$ is categorical in $\chi$.
\endproclaim
\bigskip

\remark{\stag{X2.4A} Remark}  1) We intend to try to 
prove in future work that also in 
$K_{>\lambda}$ we have categoricity and deal with $\lambda$ not successor.  
This calls for using ${\Cal P}^-(n)$-diagrams as in \cite{Sh:87a}, 
\cite{Sh:87b}, etc. \nl
2) We need some theory of orthogonality and regular
types parallel to \cite[Ch.V]{Sh:a} = \cite[Ch.V]{Sh:c}, as done in
\cite{Sh:h} and then \cite{MaSh:285} (which appeared) and then (without the
upward categoricity) \cite{KlSh:362}, \cite{Sh:472}.  Then the categoricity 
can be proved as in those papers.
\endremark
\bigskip

\demo{Proof}  Let $K' = \{M \in K:M \text{ is } \chi(\Phi^*)
\text{-saturated of cardinality } \ge \chi(\Phi^*)\}$.  So
\mr
\item "{$(*)_0$}"  there is $M \in K_\lambda$ which is $\lambda$-saturated
\nl
[why?  by \scite{2.3}, \scite{1.5}, as $\lambda$ is regular]
\sn
\item "{$(*)_1$}"  $K'$ is closed under $\le_{\frak K}$-increasing unions
\sn
\item "{$(*)_2$}"  if $\chi \ge \beth_{(2^{LS({\frak K})})^+}(\chi(\Phi^*))$
and $M \in K_\chi$ \ub{then} $M \in K'_\chi$, moreover $M$ is \nl
$\beth_{(2^{LS({\frak K})})^+}(\chi(\Phi^*))$-saturated \nl
[Why?  Otherwise by \scite{X1.3A} there is a non LS$({\frak K})^+$-saturated
$M \in K_\lambda$ contradicting $(*)_0$, or use \scite{X1.7}.  
For the ``Moreover" use \scite{X1.3B} instead of \scite{X1.3A}]
\sn
\item "{$(*)_3$}"  if $M \in K',p \in {\Cal S}(M)$ \ub{then} for some
$M_0 \le_{\frak K} M,M_0 \in K'_{\chi(\Phi^*)}$ and $p$ does not 
$\chi(\Phi^*)$-split over $M_0$ \nl
[why?  by \scite{3.2}, \scite{1.5}] 
\sn
\item "{$(*)_4$}"  \ub{Definition}: for $\chi \in [\chi(\Phi^*),\lambda)$ and
$M \in K'_\chi$ and $p \in {\Cal S}(M)$ we say $p$ is minimal if
\sn
{\roster
\itemitem{ $(a)$ }  $p$ is not algebraic which means $p$ is not realized by
any $c \in M$
\sn
\itemitem{ $(b)$ }  if $M \le_{\frak K} M' \in K'_\chi$ and $p_1,p_2 \in
{\Cal S}(M')$ are non-algebraic extending $p$, \ub{then} $p_1 = p_2$
\endroster}
\item "{$(*)_5$}"   \ub{Fact}:  if $M \in K'_\chi$ is saturated, $\chi \in
[\chi(\Phi^*),\lambda)$, \ub{then} some $p \in {\Cal S}(M)$ is minimal
\nl
[Why?  If not, we choose by induction on $\alpha \le \chi$ for every $\eta \in
{}^\alpha 2$ and triple $(M_\eta,N_\eta,a_\eta)$ and $h_{\eta,\eta
\restriction \beta}$ for $\beta \le \alpha$ such that:
\sn
{\roster
\itemitem{ $(a)$ }  $M_\eta <_{\frak K} N_\eta$ and $a_\eta \in N_\eta
\backslash M_\eta$
\sn
\itemitem{ $(b)$ }  $\langle M_{\eta \restriction \beta}:\beta \le \alpha
\rangle$ is $\le_{\frak K}$-increasingly continuous
\sn
\itemitem{ $(c)$ }  $M_{\eta \restriction \beta} <^1_{\mu,\omega}
M_{\eta \restriction (\beta +1)}$
\sn
\itemitem{ $(d)$ }  $h_{\eta,\eta \restriction \beta}$ is a $\le_{\frak K}$-
embedding of $N_{\eta \restriction \beta}$ into $N_\eta$ which is the identity
on $M_{\eta \restriction \beta}$ and maps $a_{\eta \restriction \beta}$ to
$a_\eta$
\sn
\itemitem{ $(e)$ }  if $\gamma \le \beta \le \alpha,\eta \in {}^\alpha 2$,
then $h_{\eta,\eta \restriction \gamma} = h_{\eta,\eta \restriction \beta}
\circ h_{\eta \restriction \beta,\eta \restriction \gamma}$
\sn
\itemitem{ $(f)$ }  $M_{\eta \char 94 \langle 0 \rangle} =
M_{\eta \char 94 \langle 1 \rangle}$ but \nl
tp$(a_{\eta \char 94 
\langle 0 \rangle},M_{\eta \char 94 \langle 0 \rangle},
N_{\eta \char 94 \langle 0 \rangle}) \ne \text{ tp}(a_{\eta \char 94
\langle 1 \rangle},M_{\eta \char 94 \langle 1 \rangle},
N_{\eta \char 94 \langle 1 \rangle})$
\sn
\itemitem{ $(g)$ }  $M_\eta <_{\frak K} {\frak C}$.
\endroster}
No problem to carry the definition and 
let $\kappa = \text{ Min}\{\kappa:2^\kappa > \chi\}$ and choose
$M <_{\frak K} {\frak C},\|M\| \le \chi,\eta \in {}^{\kappa >} 2 \Rightarrow
M_\eta \subseteq M$ hence $\eta \in {}^\kappa 2 \Rightarrow M_\eta
\subseteq M$ so $\{\text{tp}(a_\eta,M,{\frak C}):\eta \in {}^\kappa 2\}$ is
a subset of ${\Cal S}(M)$ of cardinality $2^\kappa > \chi$.  So then we 
can get a contradiction to stability in $\chi$].
\sn
\item "{$(*)_6$}"  Fix $M^* \in K'_{\chi(\Phi^*)}$ and minimal $p^* \in
{\Cal S}(M^*)$
\sn
\item "{$(*)_7$}"  if $M^* \le_{\frak K} M \in K'_{< \lambda}$, \ub{then}
$p^*$ has a non-algebraic extension $p \in {\Cal S}(M)$, moreover, if $M$
is saturated, it is unique and also $p$ is minimal \nl
[Why?  Existence by \scite{6.2}, uniqueness modulo $E_{\chi(\Phi^*)}$ follows
hence uniqueness by locality lemma \scite{X2.2}.  
Applying this to a saturated extension $M'$ of $M$ of cardinality $\|M\|$ we
get $p$ is minimal].
\ermn
Let $\lambda_1$ be the predecessor of $\lambda$.
\mr
\item "{$(*)_8$}"  there are no $M_1,M_2$ such that:
{\roster
\itemitem{ $(a)$ }  $M^* \le_{\frak K} M_1 \le_{\frak K} M_2$
\sn
\itemitem{ $(b)$ }  $M_1,M_2$ are saturated of cardinality $\lambda_1$
\sn
\itemitem{ $(c)$ }  $M_1 \ne M_2$
\sn
\itemitem{ $(d)$ }  no $c \in M_2 \backslash M_1$ realizes $p^*$
\endroster}
[Why?  If there are, we choose by induction on $\zeta < \lambda,N_\zeta \in
K_{\lambda_1}$ is $\le_{\frak K}$-increasingly continuous, each $N_\zeta$
is saturated, $N_0 = M_1,N_\zeta \ne N_{\zeta +1}$ and no $c \in
N_{\zeta +1} \backslash N_\zeta$ realizes $p^*$.  If we succeed, then
$N = \dbcu_{\zeta < \lambda} N_\zeta$ is in $K_\lambda$ (as $N_\zeta \ne
N_{\zeta +1}!$) but no $c \in N \backslash N_\zeta$ realizes $p^*$ \nl
(why?  as $\{\zeta:c \notin N_\zeta\}$ is an initial segment of 
$\lambda$, non-empty as $0$ is in so it has a last element $\zeta$, so 
$c \in N_{\zeta +1} \backslash N_\zeta$ so realizes $p^*$, contradiction); 
hence $N$ is not saturated, contradiction.  For $\zeta =0,
N_0 = M_1$ is okay by clause (b).  If $\zeta$ is limit $< \lambda$, let
$N_\zeta = \dbcu_{\varepsilon < \zeta} N_\varepsilon$, clearly $N_\zeta \in
K_{\lambda_1}$ and it is saturated by \scite{X2.3A}.  If $\zeta = \varepsilon
+1$, note that as $N_\varepsilon,M_1$ are saturated and in $K_{\lambda_1}$ and
$\le_{\frak K}$-extends $M^*$ which has smaller cardinality, there is an
isomorphism $f_\zeta$ from $M_1$ onto $N_\varepsilon$ which is the identity
on $M^*$.  We define $N_\zeta$ such that there is an isomorphism $f^+_\zeta$
from $M_2$ onto $N_\zeta$ extending $f_\zeta$.  
By assumption (b), $N_\zeta \in
K_{\lambda_1}$ is saturated by assumption (c), $N_\zeta \ne N_{\zeta +1}$, and
by assumption (d), no $c \in N_{\zeta +1} \backslash N_\zeta$ realizes $p^*$
(as $f \restriction M^* =$ the identity).  So as said above, we have derived
the desired contradiction].
\sn
\item "{$(*)_9$}"  if $M \in K'_{< \lambda}$ and $M^* \le M <_{\frak K} N,
M$ has cardinality $\ge \theta^* = \beth_{(2^{\chi(\Phi^*)})^+}$, 
\ub{then} some $c \in N \backslash M$ realizes $p^*$. \nl
[Why?  By $(*)_2$, $M,N$ are $\theta^*$-saturated.  So we can find saturated
$M' \le_{\frak K} M,N' \le_{\frak K} N$ of cardinality $\theta^*$ such that
$M' = N' \cap M,M^* \ne N'$ (why? by observation \scite{X2.3B}).  So still
no $c \in N' \backslash M'$ realizes $p^*$.  We would like to transfer the
appropriate omitting type theorem of this situation from $\theta^*$ to
$\lambda_1$; the least trivial point is preserving the saturation.  But this
can be expressed as: ``is isomorphic to $EM(I,\Phi)$ for some linear order 
$I$" for appropriate $\Phi$, and this is easily transferred].
\sn
\item "{$(*)_{10}$}"  if $M \in K'_{\le \lambda}$ has cardinality 
$\ge \theta^* = \beth_{(2^{\chi(\Phi^*)})^+}$ then it is 
$\theta^*$-saturated \nl
(so $\in K'_{\le \lambda}$) \nl
[why?  included in the proof of $(*)_9$]
\sn
\item "{$(*)_{11}$}"  if $M \in K'_{\le \lambda}$ has cardinality 
$\ge \theta^*$, \ub{then} $M$ is saturated \nl
[why?  Assume not; by $(*)_{10}$, $M$ is $\theta^*$-saturated let 
$\theta$ be such that $M$ is $\theta$-saturated but not 
$\theta^+$-saturated; by $(*)_{10}$, $\theta \ge \theta^*$, without loss
of generality $M^* \le_{\frak K} M$.  
Let $M_0 \le_{\frak K} M$ be such
that $M_0 \in K_\theta$ and some $q \in {\Cal S}(M_0)$ is omitted by $M$ and
\wilog \, $M^* \le_{\frak K} M_0$.
\nl
Now choose by induction on $i < \theta^+$ a triple $(N^0_i,N^1_i,f_i)$ such
that:
\sn
{\roster
\itemitem{ $(a)$ }  $N^0_i \le_{\frak K} N^1_i$ belong to $K_\theta$ and
are saturated
\sn
\itemitem{ $(b)$ }  $N^0_i$ is $\le_{\frak K}$-increasingly continuous
\sn
\itemitem{ $(c)$ }  $N^1_i$ is $\le_{\frak K}$-increasingly continuous
\sn
\itemitem { $(d)$ }  $N^0_0 = M_0$ and $d \in N^1_0$ realizes $q$
\sn
\itemitem{ $(e)$ }  $f_i$ is a $\le_{\frak K}$-embedding of $N^0_i$ into $M$
and $f_0 = \text{ id}_{M_0}$
\sn
\itemitem{ $(f)$ }  for each $i$, for some $c_i \in N^1_i \backslash N^0_i$
we have $c_i \in N^0_{i+1}$
\sn
\itemitem{ $(g)$ }  $f_i$ is increasing continuous. \nl
If we succeed, let $E = \{\delta < \theta^+:\delta$ limit and for every
$i < \delta$ and $c \in N^1_i$ we have $(\exists j < \theta^+)(c_j = c)
\rightarrow (\exists j < \delta)(c_j = c)\}$.  Clearly $E$ is a club of
$\theta^+$, and for each $\delta \in E,c_\delta$ belongs to $N^1_\delta =
\dbcu_{i < \delta} N^1_i$ so there is $i < \delta$ such that $c_\delta \in
N^1_i$, so for some $j < \delta,c = c_j$ so $c_\delta = c_j \in N^0_{j+1}
\le_{\frak K} N^0_\delta$, contradiction to clause (f). \nl
So we are stuck for some $\zeta$, now $\zeta \ne 0$ trivially.  Also
$\zeta$ not limit by \scite{X2.3A}, so $\zeta = \varepsilon + 1$.  Now if 
$N^0_\varepsilon = N^1_\varepsilon$, then $f_\varepsilon(d) \in M$ realizes
$q$ a contradiction, so $N^0_\varepsilon <_{\frak K} N^1_\varepsilon$.  Also
$f_\varepsilon(N^0_\varepsilon) <_{\frak K} M$ by cardinality consideration.
Now by $(*)_9$ some $c_\varepsilon \in N^1_\varepsilon \backslash 
N^0_\varepsilon$ realizes  $p^*$. \nl
We can find $N'_\zeta \le_{\frak K} M$ such that 
$f_\varepsilon(N^0_\varepsilon) <_{\frak K} N'_\zeta \in K_\theta,
N'_\zeta$ saturated (why? by \scite{X2.3B}). \nl
Again by $(*)_9$ we can find $c'_\zeta \in N'_\zeta \backslash f_\varepsilon
(N^0_\varepsilon)$ realizing $p^*$.  By $(*)_5$ clearly tp$(c'_\varepsilon,
f_\varepsilon(N^0_\varepsilon),M) = f_\varepsilon(\text{tp}(c_\varepsilon,
N^0_\varepsilon,N^1_\varepsilon))$ so we can find $N^1_\zeta \in K_\theta$
which is a $\le_{\frak K}$-extension of $N^1_\varepsilon$ and a
$\le_{\frak K}$-embedding $g_\varepsilon$ of $N'_\zeta$ into $N^1_\zeta$
which extends $f^{-1}_\varepsilon$ and maps $c'_\varepsilon$ to
$c_\varepsilon$.  Without loss of generality $N^1_\zeta$ is saturated.  Let
$N^0_\zeta = g_\varepsilon(N'_\zeta)$ and $N^1_\zeta,c_\varepsilon$ were
already defined. 
So we can carry the construction, contradiction, so $(*)_{11}$ holds].
\endroster}
\sn
\item "{$(*)_{12}$}"  $K_\lambda$ is categorical in every $\chi \in
[\beth_{(2^{\chi(\Phi^*)})^+},\lambda)$ \nl
[why?  by $(*)_{11}$ every model is saturated and the saturated model is
unique].  \nl
${{}}$ \hfill$\square_{\scite{X2.4}}$
\endroster
\enddemo
\newpage
\enddemo

\shlhetal
\newpage
    
REFERENCES.  
\bibliographystyle{lit-plain}
\bibliography{lista,listb,listx,listf,liste}

\enddocument

\bye